\documentclass[10pt]{article}
\usepackage[francais,english]{babel}
\usepackage{amsmath}
\usepackage{amsfonts}
\usepackage{amssymb}
\usepackage{amsthm}
\usepackage{graphics}
\usepackage{amscd}
\usepackage{fullpage}
\usepackage{epsfig}
\usepackage{color}
\newcommand{\R}{\mathbb{R}}

\newcommand{\C}{\mathbb{C}}

\newcommand{\U}{\mathbb{U}}
\newcommand{\E}{\mathbb{E}}
\newcommand{\isom}{\stackrel{\sim}{\longrightarrow}}
\newcommand{\mc}{\mathcal}
\newcommand{\mb}{\mathbb}

\newcommand{\eps}{\varepsilon}
\newcommand{\ind}{{\bf 1}}
\renewcommand{\P}{\mathbb{P}}
\newcommand{\sH}{{\mathcal{H}}}
\newcommand{\cc}{\mathbf{c}}

\renewcommand{\H}{\mathbb{H}}

\newcommand{\Lap}{\Delta\!}

\newcommand{\D}{\mathbb{D}}

\DeclareMathOperator{\End}{End}
\DeclareMathOperator{\Ker}{Ker}

\DeclareMathOperator{\Tr}{Tr}

\DeclareMathOperator{\SLE}{SLE}
\DeclareMathOperator{\FF}{FF}

\DeclareMathOperator{\wind}{wind}
\DeclareMathOperator{\Harm}{Harm}
\DeclareMathOperator{\dist}{dist}

\title{SLE and the free field: Partition functions and couplings}
\author{Julien Dub\'edat\footnote{Partially supported by NSF grant DMS0804314}}
\newtheorem{thm}{Theorem}[section]

\newtheorem{Thm}[thm]{Theorem}
\newtheorem{Def}[thm]{Definition}
\newtheorem{Prop}[thm]{Proposition}
\newtheorem{Lem}[thm]{Lemma}
\newtheorem{Cor}[thm]{Corollary}

\newtheorem{Rem}[thm]{Remark}

\begin{document}
\maketitle
\begin{abstract}
Schramm-Loewner Evolutions ($\SLE$) are random curves in planar simply connected domains; the massless (Euclidean) free field in such a domain is a random distribution. Both have conformal invariance properties in law. In the present article, some relations between the two objects are studied. We establish identities of partition functions between different versions of $\SLE$ and the free field with appropriate boundary conditions; this involves $\zeta$-regularization and the Polyakov-Alvarez conformal anomaly formula. We proceed with a construction of couplings of  $\SLE$ with the free field, showing that, in a precise sense, chordal $\SLE$ is the solution of a stochastic ``differential" equation driven by the free field. Existence, uniqueness in law, and pathwise uniqueness for these SDEs are proved for general $\kappa>0$.
\end{abstract}

\tableofcontents

\setcounter{section}{-1}
\section{Introduction}

In 2d statistical mechanics, various important models such as percolation or the Ising model are expected (or proved) to have, at criticality, a conformally invariant scaling limit. The general notion of conformal invariance underlaid the development of Conformal Field Theory. In 1999, Schramm (\cite{Sch99}) proposed a precise version of the notion of conformal invariance in distribution. This consists in considering the distribution of an isolated path (typically, an interface in the model) connecting two boundary points of a simply connected planar domain (in the chordal case). One obtains a collection of distributions $(\mu^{\SLE}_c)$ on simple paths indexed by configurations, viz. domains with two marked boundary points; the conformal invariance requirement reads $\varphi_*\mu^{\SLE}_c=\mu^{\SLE}_{\varphi(c)}$ for a conformal equivalence $\varphi:c\rightarrow\varphi(c)$ (in other words, $\mu^{\SLE}$ is a covariant functor on the groupoid of configurations). Under the conformal invariance requirement and an additional domain Markov property, the collection of measures is classified by a positive parameter $\kappa>0$ (\cite{Sch99}).

Another type of conformally invariant scaling limits involves distributions. In the case of dimers, a height function is associated to a configuration, following the definition of Thurston; Kenyon (\cite{Ken_domino_conformal,Ken_domino_GFF}) proved that in the case of the square lattice, with appropriate boundary conditions, this height converges in distribution to the massless free field. This is the Gaussian measure with covariance operator given by the Green kernel with Dirichlet boundary conditions. It can be seen as a random distribution (element of $C^\infty_0(D)'$) and is a basic object in constructive Field Theory (\cite{Simon_Pphi,Glimm_Jaffe}).

Temperley's bijection (see eg \cite{KPW}) relates dimer configurations (tilings) to uniform spanning trees; branches of these trees are distributed as loop-erased random walks. In this discrete setting, two types of invariance principle may be considered: a branch converges to $\SLE_2$, as proved by Lawler, Schramm, Werner (\cite{LSW_LERW}); the height function (at least in closely related set-ups) converges to a free field. Moreover, in the discrete setting, the height function determines the branches. The relation between the height of the tiling and the branches can be understood in terms of winding (of a curve running along the branches on the medial lattice, \cite{KPW}), as first conjectured by Benjamini. It was then proved that the scaling limit of (the Peano path of) the tree is $\SLE_8$ (\cite{LSW_LERW}). A question raised in \cite{Ken_domino_GFF} is whether the reconstruction of the tree from the height function, which is possible in the discrete set-up, can be carried out in the continuum. This will be answered affirmatively in Section 8.

In \cite{SS_freefield}, Schramm and Sheffield prove that the zero level line of a discrete Gaussian free field on the triangular lattice (with appropriate boundary conditions) converges in distribution to chordal $\SLE_4$, as the mesh goes to zero. Trivially, the discrete free field converges to the continuous massless free field. The relation between chordal $\SLE_4$ and the free field in the scaling limit, in particular in terms of couplings, is studied in details in the forthcoming \cite{SS_FF2}. A closely related situation is that of double domino tilings, that was conjectured by Kenyon to lie in the same universality class.

Work in progress relating the free field and $\SLE_\kappa$ for $\kappa\neq 4$ has been reported by Scott Sheffield, based partly on the ``winding" of $\SLE$ curves, seen as ``flow lines of $e^{ich}$", $h$ a free field, $c$ a parameter. A notion of ``local sets" of the free field, that applies to and extends the case of contour lines, has also been advanced.

In the examples of spanning trees and double domino tilings/discrete free field, two types of boundary conditions for fields appear: piecewise constant, with jumps at prescribed points; and a multiple of the winding of the boundary curve, again with jumps at prescribed points.

In the present article, we study relations between different variants of $\SLE$ and the free field with appropriate boundary conditions. The first  main result concerns {\em partition functions} of $\SLE$ and the free field. For the free field, the partition function is defined in a natural way from its Gaussian structure:
$${\mc Z}^{\FF}=\det(\Lap)^{-\frac 12}\exp(-\frac 12\langle m,m\rangle_{\sH^1})$$
where $m$ is the mean of the field. The Laplacian (in a bounded domain with Dirichlet condition on its smooth boundary) has a discrete spectrum going to infinity. In this situation, it is customary to resort to the $\zeta$-regularized ${\det}_\zeta(\Lap)$. Partition functions of $\SLE$ are defined in a way compatible with its absolute continuity properties; the form of the partition function is ${\mc Z}^{\SLE}=\det(\Lap)^{-\frac {\cc}2}$ times a conformally invariant tensor, where $\cc$ is the central charge, that depends on $\kappa$. For many variants of $\SLE$ (chordal, radial, multiple chordal and radial), we match the boundary condition (involving the winding of the boundary) of the free field and the $\SLE$ variant in such a way that the identity of partition functions (see Theorem \ref{PFident}):
$${\mc Z}^{\FF}={\mc Z}^{\SLE}$$
holds. These are functions on the configuration space (a configuration being now equipped with a Riemannian metric, not merely a complex structure), that are well defined up to a multiplicative constant. The apparent mismatch of exponents of the Laplacian determinant is resolved via the Polyakov-Alvarez conformal anomaly formula. We note that these partition functions are also relevant to Conformal Field Theory (as correlators of primary fields) and Virasoro representations, as detailed in the forthcoming \cite{Dub_Vir}. For earlier considerations on partition functions/CFT correlators in relation with $\SLE$, see \cite{FriKal,BB_review,KontSuh} and references therein.

When a field and an $\SLE$ are matched through their partition functions, one gets easily a ``local" coupling restricted to the $\SLE$ and field seen in disjoint subdomains. This plays a r\^ole closely analogous to that of local commutation of $\SLE$s considered in \cite{Dub_Comm} (here, the two ``commuting" objects are an $\SLE$ and a field, instead of two $\SLE$s). In the context of $\SLE$ reversibility, Zhan showed in \cite{DZ_revers} how to lift local couplings to global couplings. In \cite{Dub_dual}, it is shown how to extend this to the framework of local commutation, in which partition functions intervene naturally. We use similar techniques here to couple in a domain one $\SLE$ strand with a free field, in conjunction with Gaussian arguments and free field properties. One can also couple systems of commuting $\SLE$s with a free field, in such a way that the different $\SLE$ strands are independent conditionally on the field; the identity of partition functions is instrumental at this point.

In order to elucidate the nature of these results, we introduce a notion of stochastic ``differential" equation driven by the free field, by analogy with the classical framework of SDEs driven by linear Brownian motion. The relation between the $\SLE$ path and the free field does not involve a stochastic calculus, but is a condition that can be checked pathwise by an explicit construction. Informally, the field near the $\SLE$ trace converges to its boundary value given the position of the trace; some care has to be given to the fact that this boundary value is not defined on the trace (for $\kappa\neq 4$), which is rough. The equation reads:
$$\E(\phi_{|D\setminus K_t}|{\mc F}^{\FF}_{\partial K_t})=h((K^t_.))\quad \forall t\geq 0$$
where $\phi$ is the field, $K$ the Loewner chain and $h$ an harmonic function depending on the chain (playing the r\^ole of SDE coefficients $(\sigma,b)$). As in the case of SDEs, there are conditions of adaptness w.r.t. a filtration; the filtration is indexed here by a partially ordered set (for inclusion) of open subsets of the domain.

In this context, we prove that chordal $\SLE_\kappa$, for $\kappa>0$, is a solution of a stochastic equation driven by the free field, for which uniqueness in law holds (Theorem \ref{Thm:SDEweak}). The compatibility of the construction with various duality identities (reversibility for $\kappa=4$) leads to a proof of pathwise uniqueness for general $\kappa>0$ (Theorem \ref{Thm:SDEstrong}). Both existence and pathwise uniqueness are local properties, hence hold in a variety of set-ups.

The article is organized as follows. In Section 1, we discuss discrete couplings. Section 2 contains results on the Brownian loop measure, in particular in relation with functional ($\zeta$ and Fredholm) determinants. Schramm-Loewner Evolutions are discussed in Section 3, with emphasis on partition functions. Section 4 gathers material on the free field. Relevant boundary conditions are introduced in Section 5, before establishing identities of partition function. Section 6 is concerned with local and global couplings. Stochastic equations driven by the free field are discussed in Section 7, where uniqueness results are proved. Some consequences (Temperley's bijection in the continuum and strong duality identities) are discussed in Section 8.


\section{Discrete couplings}

In this section, we discuss some examples of discrete couplings between a path converging to some $\SLE$ and a field converging to the free field, for motivation and intuition.  

\subsection{The Temperley coupling}

A complete discussion would involve introducing a lot of material that will not be used later in the article, so we shall only sketch the construction, in the case of the square lattice. For a detailed treatment, see eg \cite{KPW} and references therein.

\begin{figure}[htb!]
\begin{center}
\centerline{\scalebox{.6}{\rotatebox{270}{\psfig{file=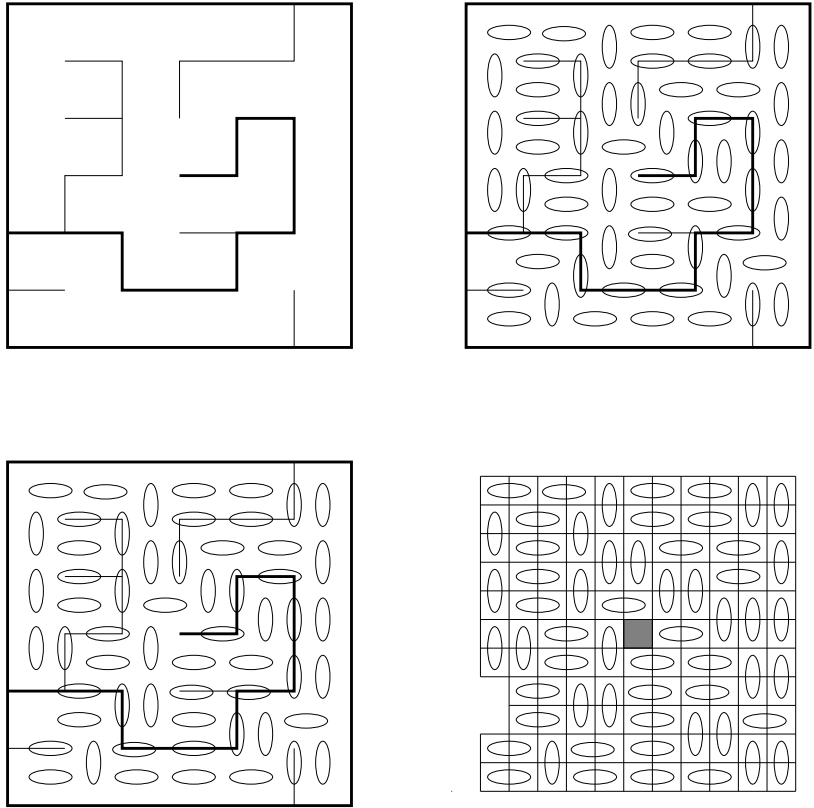}}}}
\end{center}
\caption{1. a LERW in a spanning tree rooted at the boundary. 2. Associated dime
r configuration. 3,4. Dimer configuration after slide.}\label{fig:radial}
\end{figure}

Consider a portion of the square lattice approximating a simply connected domain: this gives a finite graph $\Gamma$. The outer boundary is seen as a single extended vertex. A spanning tree on the graph rooted at the extended vertex determines by planar duality a spanning tree on the dual graph $\Gamma^\dagger$.
The two graphs may be oriented towards their root (this involves picking a root on the dual graph). From each vertex of the graph starts an outgoing edge in the tree. A square lattice with a twice smaller mesh can be constructed by superimposing the original graph $\Gamma$ and its dual $\Gamma^{\dagger}$: this gives a graph $D\Gamma$, which is bipartite (black vertices of $D\Gamma$ correspond to vertices and faces of $\Gamma$, white vertices to edges of $\Gamma$). To an oriented edge in the original graph $\Gamma$ or the dual $\Gamma^\dagger$, one can associate an edge of the new graph $D\Gamma$: its initial half. Thus a spanning tree on $\Gamma$ determines first a dual tree on $\Gamma^\dagger$ and this two trees yield a collection of edges in $D\Gamma$. Being careful with the treatment of the boundary, this collection of edges is a perfect matching on the bipartite graph $D\Gamma$. This describes Temperley's bijection between uniform spanning trees on $\Gamma$ and perfect matching of a related graph $D\Gamma$; see Figure \ref{fig:radial}, panel 2.


To the dimer configuration is associated an integer valued height function on vertices of the medial lattice $(D\Gamma)^\dagger$, as defined by Thurston. The variation of the height along an edge of $(D\Gamma)^\dagger$ with a black vertex of $D\Gamma$ to its left is $(-3)$ if it crosses a matched edge of $D\Gamma$ and $1$ otherwise. 
The data of the original tree, the dimer configuration, and the (admissible) height function are equivalent. Kenyon proved (\cite{Ken_domino_conformal,Ken_domino_GFF}) that, as the mesh of the lattice goes to zero, the height converges in distribution (in a weak topology) to the free field, with boundary value given by a multiple of the winding of the boundary, measured from the root. 

In order to make the connection with $\SLE$, it seems convenient to modify the situation as follows. Picking a point $y$ inside the domain, one can consider the branch of the tree from $y$ to the boundary. The branch hits the boundary at $x$, it is also convenient to condition on $x$. 
It is known that the branch is distributed as a Loop Erased Random Walk and that its scaling limit is radial $\SLE_2$ from $x$ to $y$ (\cite{LSW_LERW}). Given the tree, one obtains a dimer configuration by Temperley's bijection; one can create a hole by sliding the dimers along the LERW, a construction introduced by Kenyon; see Figure \ref{fig:radial}. The height function has now an additive monodromy around the puncture $y$ (ie is additively multiply valued, picking an additive constant 4 when traced counterclockwise along a simple loop around $y$). We stress however that the dimer configuration is not uniform on tilings of the punctured graph (it is uniform on a set of admissible matchings).

In this context, one also gets a natural definition of the partition function of the LERW from $y$ to $x$: the number of spanning trees of  $\Gamma$ such that the branch starting from $y$ exits at $x$. The total number of spanning trees in $\Gamma$ is, by the well known Matrix Tree Theorem, $\det(\Lap_\Gamma)$ where $\Lap_\Gamma$ is the combinatorial Laplacian on $\Gamma$ with Dirichlet conditions on the boundary (so that it is invertible). The probability that the LERW exits at $x$ is the probability that the underlying random walk exits at $x$, ie $\Harm_\Gamma(y,\{x\})$. This gives a partition function:
$${\mc Z}^{\rm LERW}=\det(\Lap_\Gamma)\Harm_\Gamma(y,\{x\})$$ 
This is coherent with the partition functions of $\SLE$ in the continuum that we will consider later on; in the case of radial $\SLE_2$, it is written:
$${\mc Z}^{\SLE_2}={\det}_\zeta(\Lap)\Harm(y,dx)$$
An important point is that the partition function accounts for an ambient ``environment", through $\det(\Lap)$. Justification for including this contribution will be given later, in particular in terms of absolute continuity properties.

\subsection{Discrete free field and double dominos}

\begin{figure}[htb!]
\begin{center}
\centerline{\scalebox{.6}{\rotatebox{270}{\psfig{file=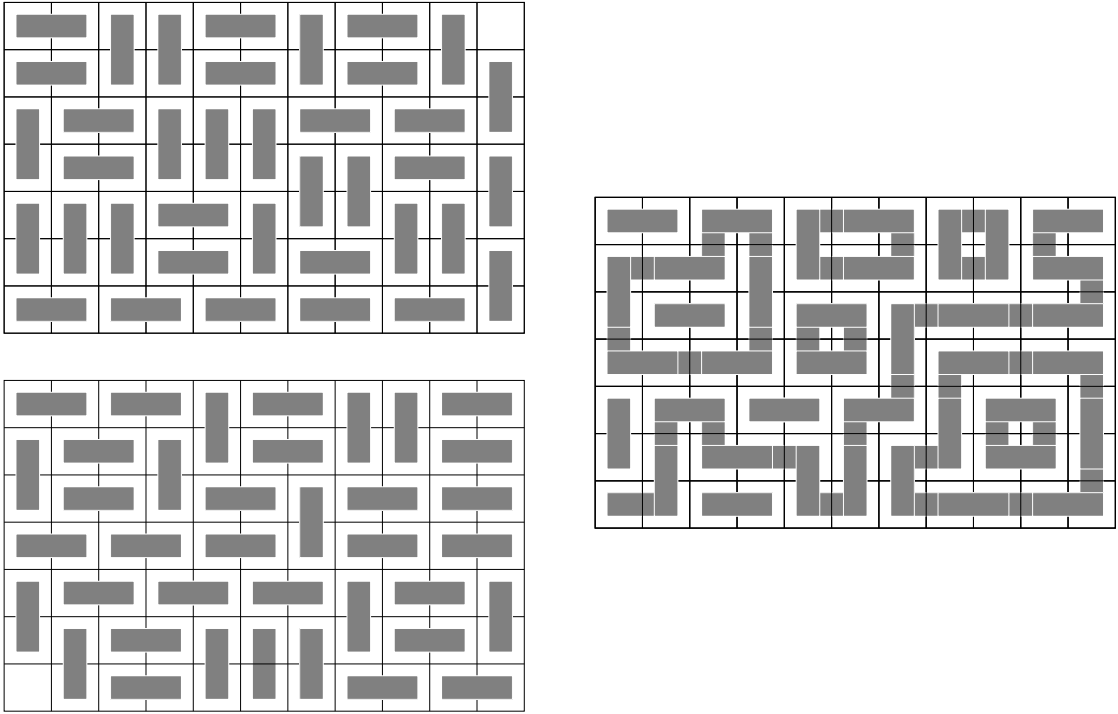}}}}
\end{center}
\caption{Up, down: two dimer tilings (rectangle minus corner). Right: superposition of the tilings, consisting of closed cycles, doubled dimers, and one open path.}
\label{fig:doubled}
\end{figure}

For continuity of the discussion, we start with the case of double domino dimers (Figure \ref{fig:doubled}). 
One possible setup is as follows. Consider a portion of the square lattice $\Gamma$, from which a boundary square is deleted (say a rectangle with odd sides minus a corner). One can sample uniformly a domino tiling of this domain. It is associated to a height function, that converges to the free field, with boundary condition given by a multiple of the winding of the boundary, jumping at the excised corner. One can proceed likewise by deleting another corner, sampling independently the domino tiling. The superposition of the two tilings consists of doubled dimers, closed loops (nested), and one open path connecting the two excised corners. Kenyon conjectured that this converges to chordal $\SLE_4$.

Each tiling is associated to an height function. The difference $h$ of the two height functions is such that $h$ jumps by $\pm 4$ when crossing a loop or the open path (and is constant along these paths). Hence one can think of the superposition of tilings in terms of contour lines of $h$. A double application of the invariance principle in \cite{Ken_domino_conformal,Ken_domino_GFF} shows that $h$ converges to a free field with piecewise constant boundary conditions, jumping at the excised corners (the winding contributions cancelling out when taking the difference of the two height functions).

A closely related situation is completely analyzed in \cite{SS_freefield}. 
Consider a portion $\Gamma$ of the triangular lattice, approximating a simply connected domain. One can project orthogonally the continuous free field on the space of functions that are piecewise constant on the triangulation. As an alternative description of the discrete field, one can take the Gaussian measure (on functions on the triangulation) with covariance given by the discrete Green kernel with Dirichlet boundary conditions; the mean of the field is the harmonic extension of a piecewise constant function on the boundary, with jumps at two marked points. From the field, one gets a coloring of vertices (black for positive values of the field, white for negative values). There is an interface running in the domain between black and white vertices, connecting the two marked boundary points. As the mesh of the lattice goes to zero, this interface converges in distribution to chordal $\SLE_4$, for a precisely tuned boundary condition. The discrete field trivially converges to the continuous free field with piecewise constant boundary conditions.

In the double domino model, the field jumps by $\pm 4$ at the marked point and, according to \cite{Ken_domino_GFF}, converges to $\frac{4\sqrt 2}{\sqrt\pi}$ times a standard free field. Normalizing the field, we find a jump of $\sqrt{\frac\pi 2}$. In the case of the discrete Gaussian free field, it is mentioned in \cite{SS_freefield} that the jump is $2\sqrt{\frac\pi 8}$ for a normalized field. In the continuum, we will consider a jump of $\pi\sqrt{\frac 2{\pi\kappa}}$, $\kappa=4$. All these expressions are thus coherent.

In these situations, it is quite natural to consider multiple paths created by, say, excising an even number of boundary squares in the double domino model or flipping the boundary condition an even number of times for the discrete free field. This is analyzed combinatorially in \cite{KenWil_grove}; see also \cite{Dub_Comm,Dub_Euler}.

In the context of the discrete free field, the Gaussian structure gives a natural definition for the partition function:
$${\mc Z}^{\FF}=\det(\Lap_\Gamma)^{-\frac 12}\exp(-\frac 12\langle m,m\rangle_\nabla)$$
where $\langle m,m\rangle_\nabla$ is the (discrete) Dirichlet energy of the mean of the field (which is also the state of minimal energy, under the boundary condition constraints). We will consider regularized versions of this in the continuum.

\section{Loop measure, determinants}

In order to relate quantities arising from $\SLE$ and free field densities, we need to introduce the loop measure \cite{LSW_restr,LW} and relate the masses of some sets under this measure to functional determinants of two types: Fredholm and $\zeta$-regularized. Some relations between loop measures, free fields and functional determinants are discussed in \cite{LJ_loops}. Most of the discussion here can be carried out at the level of Markov chains (\cite{Law_Fer_RWLS,LJ_loops}) or diffusions on manifolds; only the conformal anomaly formula is specific to the two-dimensional case.

Consider the (positive) Laplacian $\Lap$ on a compact manifold $M$ with boundary, with Dirichlet condition on the boundary (more generally, the negative generator of a diffusion). Following Ray-Singer, one attaches to $\Lap$ a $\zeta$-function:
$$\zeta_\Lap(s)=\frac{1}{\Gamma(s)}\int_0^\infty \Tr(e^{-t\Lap})t^{s-1}dt$$
where $P_t=e^{-t\Lap}$ is the transition kernel for Brownian motion (running at speed 2), trace class in $L^2(M)$ for $t>0$, and $\Tr(P_t)=\int_M P_t(x,x)dvol(x)$. This is absolutely convergent for $\Re s>(\dim M)/2$. If $\lambda_1\leq \lambda_2\leq\cdots$ is the spectrum of $\Lap$,  $\zeta_\Lap(s)=\sum_{n\geq 1}\lambda_n^{-s}$. Under regularity assumptions on the boundary, $\zeta_\Lap$ has a meromorphic extension to $\C$, in particular regular at $0$, so that one can define the spectral invariant:
$${\det}_\zeta(\Lap)=e^{-\zeta'_\Lap(0)}$$
Formally ($\Gamma(s)\sim s^{-1}$ as $s\rightarrow 0$), 
$$``-\log{\det}_\zeta(\Lap)=\int_0^\infty \Tr(e^{-t\Lap})t^{-1}dt"$$

Note that $P_t(x,y)dvol(y)$ is the disintegration of the measure on paths starting from $x'$, stopped at $t$, w.r.t. the endpoint $y$. Let us denote ${\mb W}^t_{x\rightarrow y}$ this subprobability measure (killing on the boundary). 
This gives another justification to introduce the (rooted) loop measure (\cite{LW}) (and also to normalize it this way):
$$\mu^{loop}_r=\int_0^\infty\frac{dt}{t}\int_M{\mb W}^t_{x\rightarrow x}dvol(x)$$
It is a measure on rooted loops, ie on functions $\delta:[0,\tau]\rightarrow M$ ($\tau$ the lifetime of the loop); since the endpoints are identical, $\theta_u\delta:t\mapsto \delta_{u+s\bmod\tau}$ is again a loop. This defines an equivalence relation on loops. The quotient of $\mu^{loop}_r$ under this relation is the loop measure $\mu^{loop}$. Two important properties are:
\begin{itemize}
\item (Restriction property) If $D'\subset D$, $L$ the generator of a diffusion on $D$, 
$$d\mu^{loop}_{L,D'}(\delta)=\ind_{\delta\subset D'}d\mu^{loop}_{L,D}(\delta)$$
\item (Conformal invariance) If $L$ a generator on $D$, $\sigma$ a function on $D$, then up to time reparameterization of loops,
$$\mu^{loop}_{e^{2\sigma}L,D}=\mu^{loop}_{L,D}$$
\end{itemize}
Moreover,
$$\zeta_L(s)=\int \frac{\tau(\delta)^s}{\Gamma(s)}d\mu^{loop}(\delta)$$
for $\Re(s)>1$. Again formally, we have:
$$``-\log(\det(L))=\int d\mu_L^{loop}{\rm\ }"$$
The divergence of the RHS comes from small loops. We will be able to phrase identities by various inclusion-exclusion arguments that cancel the small loops.

Of particular interest is the following quantity: if $D$ is a domain, $K_1,K_2$ are disjoint subsets of $D$ (say connected closed sets), then define:
$$m^l(D;K_1,K_2)=\mu^{loop}_D\{\delta\cap K_1\neq\varnothing, \delta\cap K_2\neq\varnothing\}$$  
It is a conformal invariant (if $\psi:D\rightarrow D'$ is a conformal equivalence, $m^l(\psi(D);\psi(K_1),\psi(K_2))=m^l(D;K_1,K_2)$). Typically, $K_i$ is a crosscut or a hull attached to $\partial D$. We have:

\begin{Prop}\label{loop_zeta}
If $D,K_1,K_2$ are bounded with smooth boundary, then:
$$\exp(-m^l(D;K_1,K_2))=\frac{{\det}_\zeta(\Lap_D){\det}_\zeta(\Lap_{D\setminus (K_1\cup K_2)})}{{\det}_\zeta(\Lap_{D\setminus K_1}){\det}_\zeta(\Lap_{D\setminus K_2})}$$
\end{Prop}
\begin{proof}
Under these assumptions, the $\zeta$ functions have a meromorphic extension to $\C$ and are regular at $0$. Then:
\begin{align*}
\Gamma(s)\left(\zeta_D+\zeta_{D\setminus (K_1\cup K_2)}-\zeta_{D\setminus K_1}-\zeta_{D\setminus K_2}\right)(s)
&=\int \tau(\delta)^s d\mu^{loop}_{D}(\delta)+\int \tau(\delta)^s d\mu^{loop}_{D\setminus (K_1\cup K_2)}(\delta)-\cdots\\
&=\int \tau(\delta)^s(1+\ind_{\delta\cap(K_1\cup K_2)=\varnothing}-\ind_{\delta\cap K_1=\varnothing}-\ind_{\delta\cap K_2=\varnothing})d\mu^{loop}_D(\delta)\\
&=\int \tau(\delta)^s \ind_{\delta\cap K_1\neq\varnothing,\delta\cap K_2\neq\varnothing}d\mu^{loop}_D(\delta)
\end{align*}\\
using the restriction property. It is easy to see that the measure $\mu^{loop}_D$ restricted to loops intersecting both $K_1$ and $K_2$ is finite (with mass $m^l(D;K_1,K_2)$), and that the RHS is an entire function in $s$ (the mass of loops connecting $K_1$ to $K_2$ in a short time $t$ is of order $\exp(-dist(K_1,K_2)^2/t)$, from the Varadhan large deviation estimate for the heat kernel). So taking the derivative at 0 gives the result.
\end{proof}
Note that the LHS is defined under more general assumptions ($K_1,K_2$ have positive capacity) than the RHS.

Another expression in terms of Fredholm determinants $\det_F$ is also useful. Again $D$ is a domain; $K_1$, $K_2$ are smooth curves (typically, crosscuts). The metric on $D$ induces a length on $K_1,K_2$. Let 
us define a map $T_{12}: L^2(K_1)\rightarrow L^2(K_2)$ by:
$$(T_{12}f)(x)=\int_{K_1}f(y)\Harm_{D\setminus {K_1}}(x,dy)$$
and $T_{21}: L^2(K_2)\rightarrow L^2(K_1)$ is defined similarly. These operators have smooth kernels and $T=T_{12}T_{21}$ is trace class on $L^2(K_1)$. 

\begin{Prop}\label{loopFred}
Under the above assumptions:
$$\exp(-m^l(D;K_1,K_2))={\det}_F(1-T_{12}T_{21})={\det}_F(1-T_{21}T_{12})$$
\end{Prop}
\begin{proof}
We have the expansion (see eg \cite{Simon_Trace}, Chapters 2-3):
$$-\log{\det}_F(1-T_{12}T_{21})=\sum_{n\geq 1}\frac 1n\Tr((T_{12}T_{21})^n)$$
More precisely, $T_{12}, T_{21}$ are strictly uniformly subMarkov kernel, ie: 
$$\max\left(\sup_{x\in K_1}\int_{K_2}T_{21}(x,y)dl(y),\sup_{y\in K_2}\int_{K_1}T_{12}(y,x)dl(x)\right)\leq 1-p$$ 
where $p>0$ is the infimum of probabilities that a particle starting from $K_i$ hits $\partial D$ before $K_{2-i}$, $i=1,2$. It follows that:
$$\sup_{x,x'\in K_1}(T_{12}T_{21})^n(x,x')\leq (1-p)^{2n-1}\sup_{x\in K_1,y\in K_2} T_{12}(x,y)$$ 
so that $\Tr((T_{12}T_{21})^n)=O((1-p)^{2n})$, and the determinant expansion is legitimate.

This expansion counts loops intersecting both $K_1$ and $K_2$ ($T_{12}$ corresponds to paths starting from $x\in K_2$ and stopped when they hit $K_1$). We just have to check that the count is correct.

Let us go back to the rooted loop measure. For a rooted loop $\delta$, consider the sequence of successive hits of $K_1$ and $K_2$ by $\delta$; this can be represented by an alternating sequence $\sigma$ of $1$'s and $2$'s.
We consider the parabolic harmonic measure $PH_{D\setminus K_1}(x,t,y)dtdl(y)$ seen for $x\in D\setminus K_1$, for $t>0$, $y\in K_1$. The mass under ${\mb W}^t_x$ of loops that start from $x$, hit $K_1$, then $K_2$, and return to $dA(x)$ at time $t$ without returning to $K_1$ is:
\begin{align*}
\mu^{loop}_r\{\sigma=(12),(212)\}=\\
\int_0^\infty\frac{dt}t\int_DdA(x)\int_{K_1}\int_0^tPH_{D\setminus K_1}(x,t_1,y)dt_1dl(y)
\int_{K_2}\int_{t_1}^tPH_{D\setminus K_2}(y,t_2,z)dt_2dl(z)P_{t-t_1-t_2,D\setminus K_1}(z,x)
\end{align*}
By the semigroup property:
$$\int_D P_{t-t_1-t_2,D\setminus K_1}(z,x)PH_{D\setminus K_1}(x,t_1,y)dA(x)=PH_{D\setminus K_1}(z,t-t_2,y)$$
So integrating out $x$, one gets:
$$\mu^{loop}_r\{\sigma=(12),(212)\}=\int_0^\infty\frac {dt}t\int_0^tdt_1\int_{t_1}^tdt_2\int_{K_1}dl(y)\int_{K_2}dl(z)PH_{D\setminus K_1}(z,t-t_2,y)PH_{D\setminus K_2}(y,t_2,z)$$
which can be rewritten as (integrating $t_1$):
$$\int_0^\infty\int_0^\infty\frac {s_2ds_1ds_2}{s_1+s_2}\int_{K_1}dl(y)\int_{K_2}dl(z)PH_{D\setminus K_1}(z,s_1,y)PH_{D\setminus K_2}(y,s_2,z)$$
If we add the symmetrized term (obtained by interchanging $K_1,K_2$), we get:
\begin{align*}
\int_0^\infty\int_0^\infty ds_1ds_2\int_{K_1\times K_2}dl(y)dl(z)PH_{D\setminus K_1}(z,s_1,y)PH_{D\setminus K_2}(y,s_2,z)\\
=\int_{K_1\times K_2}dl(y)dl(z) \Harm_{D\setminus K_1}(z,y)\Harm_{D\setminus K_2}(y,z)\\
=\Tr(T_{12}T_{21})=\Tr(T_{21}T_{12})
\end{align*}
 We have shown that:
$$\Tr(T_{12}T_{21})=\mu^{loop}_r\{\sigma=(12),(121),(21),(212)\}$$
Proceeding as above, one gets the following expression for $\mu^{loop}_r\{length(\sigma)=2n,2n+1\}$, after having integrated out the root $x$:
$$\int_{[0,\infty)^{2n}}\frac{(s_1+s_2)ds_1\dots ds_{2n}}{s_1+\cdots+s_{2n}}\int_{K_1^n\times K_2^n}\prod_{i=1}^n dl(y_i)dl(z_i)PH_{D\setminus K_2}(y_i,s_{2i-1},z_i)PH_{D\setminus K_1}(z_i,s_{2i},y_{i+1})$$
with cyclical indexing ($y_{n+1}=y_1$). Plainly, taking a cyclic permutations of indices does not change the value of that expression. So averaging over the $n$ cyclic permutations, one gets:
\begin{align*}
\mu^{loop}_r\{length(\sigma)=2n,2n+1\}&=
\frac 1n \int_{[0,\infty)^{2n}}d{\bf s}\int_{K_1^n\times K_2^n}\prod_{i=1}^n dl(y_i)dl(z_i)PH_{D\setminus K_2}(y_i,s_{2i-1},z_i)PH_{D\setminus K_1}(z_i,s_{2i},y_{i+1})\\
&=\int_{K_1^n\times K_2^n}\prod_{i=1}^n dl(y_i)dl(z_i)\Harm_{D\setminus K_2}(y_i,z_i)\Harm_{D\setminus K_1}(z_i,y_{i+1})\\
&=\frac 1n \Tr((T_{12}T_{21})^n)
\end{align*}
which concludes the proof.
\end{proof}

To illustrate the result and to fix normalization, we embark on a sample computation, similar to the one in \cite{LW}. Let $D=\H$, $K_1$ a small hull near $0$ (with half-plane capacity $2t$, see Section 3.1), $K_2$ the unit semicircle. The harmonic measure in the semidisk $\D^+=\D\cap\H$ (on the semicircle) can be obtained from the harmonic measure in the disk by a reflection principle argument:
$$\Harm_{\D^+}(z,y)dl(y)=(\Harm_\D(z,y)-\Harm_\D(z,\overline y))dl(y)$$
and classically $\Harm_\D(z,y)=\frac 1{2\pi}\Re\frac{y+z}{y-z}$. So for $z$ close to 0, 
$$\Harm{\D^+}(z,y)\simeq -2\Im(z)\frac{1}{2\pi}\Im(\frac 2y)\simeq -\frac 2\pi\Im(z)\Im(y^{-1})$$
On the other hand, starting from $y\in\U$, if $X_\tau$ is Brownian motion stopped on exiting $\H\setminus K_1$, $\E(\Im(X_\tau))\simeq -2t\Im(y^{-1})$. It follows that on $L^2(K_2)$, the operator $T_{21}T_{12}$ has kernel:
$$(T_{21}T_{12})(y_1,y_2)\simeq \frac{4t}\pi\Im(y_1^{-1})\Im(y_2^{-1})$$
It follows that:
$$m^l(\H;K_1,K_2)\simeq\Tr(T_{21}T_{12})\simeq \frac{4t}\pi\int_{-\pi}^\pi(\sin\theta)^2d\theta\simeq 4t$$
Let $x\in(0,1)$. Consider the homography: $\varphi(z)=\frac{1-xz}{x-z}$. It permutes $-1,1$, hence preserves $K_2$ (an hyperbolic geodesic); moreover $\varphi(0)=x^{-1}$, $\varphi'(0)=x^{-2}-1$. Besides, it is easy to see that $\psi(z)=z+z^{-1}$ is the conformal equivalence $\H\setminus\D\rightarrow\H$ with hydrodynamic normalization at $\infty$. Then:
\begin{align*}
(S\psi)(z)&=\frac{\psi'''}{\psi'}(z)-\frac 32\cdot\left(\frac{\psi''}{\psi'}\right)^2(z)\\
&=-\frac 6{z^4}\cdot\frac{1}{1-z^{-2}}-\frac 32\left(\frac 2{z^3}\cdot\frac 1{1-z^{-2}}\right)^2=-\frac{6}{z^2(z^2-1)}\left(1+\frac{1}{z^2-1}\right)=-\frac 6{(z^2-1)^2}
\end{align*}
Hence: 
$$m^l(\H;\varphi(K_1),\varphi(K_2))=m^l(\H,K_1,K_2)\simeq 4t= 2.2t(1-x^{-2})^2.\frac{1}{(1-x^{-2})^2}=2.hcap(\varphi(K_1)).\left(-\frac{S\psi}6(x^{-1})\right)$$

We now discuss the Polyakov-Alvarez (\cite{Pol_bosonic,Alv_boundary, OPS_extr}) conformal anomaly formula, that describes the transformation of ${\det}_\zeta(\Lap)$ under a conformal change of metric. The key point is that under a change of metric $g_0\rightarrow g=e^{2\sigma}g_0$, the Laplacian transforms as $\Lap\rightarrow e^{-2\sigma}\Lap_0$ (this is particular to dimension 2). It is then a matter of short time heat kernel asymptotics (Pleijel-Minakshisundaram in the bulk, McKean-Singer near the boundary). We also give a simplified version in the case of planar domains; this expression is used in \cite{OPS_extr} to prove that among simply connected Riemannian surfaces with given boundary length, flat disks have extremal Laplacian determinants.

\begin{Prop}[Polyakov-Alvarez conformal anomaly formula]\label{PA}\begin{enumerate}
\item Let $(M,g_0)$ be a compact surface with boundary, $g=e^{2\sigma} g_0$. Then:
\begin{align*}
\log{\det}_\zeta\Lap_\sigma-\log{\det}_\zeta\Lap_0&=-\frac 1{6\pi}\left(\frac 12\int_M|\nabla_0\sigma|^2dA_0+\int_M K_0\sigma dA_0+\int_{\partial M}k_0\sigma dl\right)-\frac 1{4\pi}\int_{\partial M}\partial_n\sigma dl
\end{align*}
where $K_0$ is the Gauss curvature and $k_0$ is the geodesic curvature of the boundary for $g_0$.
\item Let $D$ be a planar simply connected domain with smooth boundary and Euclidean metric, $\D$ the unit disk. Let $\varphi: \D\rightarrow D$ be a conformal equivalence, $\sigma=\log|\varphi'|$. Then:
$$\log{\det}_\zeta\Lap_D-\log{\det}_\zeta\Lap_\D=-\frac 1{6\pi}\left(\frac 12\int_\D|\nabla\sigma|^2dA+\int_{\partial \D}\sigma dl\right)$$
\end{enumerate}
\end{Prop}
\begin{proof}
One deduces 2 from 1 as follows. It is equivalent to consider the Euclidean Laplacian in $D$ and to consider the Laplacian in $\D$ with pulled back metric $g=|\varphi'|^2g_0$, so that the conformal factor is $\sigma=\log|\varphi'|$. In the Euclidean metric of $\D$, $K_0=0$ and $k_0\equiv 1$. Also, $\sigma$ is harmonic so that $\int_{\partial \D}\partial_n\sigma dl=0$. 
\end{proof}

\section{Schramm-Loewner Evolutions}

\subsection{Chordal SLE}

First we recall some definitions and fix notations. We briefly discuss here chordal SLE in the upper half-plane
$\H$, from a real point to $\infty$. 
Chordal SLE in other (simply connected) domains are obtained by conformal equivalence. We will use chordal $\SLE$ both in itself and as a reference measure. For general background on SLE, see \cite{RS01,W1,Law}. 

Consider the family of ODEs, indexed by $z$ in $\H$:
$$\partial_tg_t(z)=\frac 2{g_t(z)-W_t}$$ 
with initial conditions $g_0(z)=z$, where $W_t$ is some real-valued
(continuous) function. These chordal Loewner equations are defined up
to explosion time $\tau_z$ (possibly infinite). Define:
$$K_t=\overline{\{z\in\H:\tau_z<t\}}.$$
Then $(K_t)_{t\geq 0}$ is an increasing family of compact subsets of
$\overline\H$; moreover, $g_t$ is the unique conformal equivalence
$\H\setminus K_t\rightarrow \H$ such that (hydrodynamic normalization
at $\infty$):
$$g_t(z)=z+o(1).$$ 
The coefficient of $1/z$ in the Laurent expansion of $g_t$ at $\infty$
is by definition the half-plane capacity of $K_t$ at infinity; this
capacity equals $(2t)$.

If $W_t=x+\sqrt\kappa B_t$ where $(B_t)$ is a standard Brownian
motion, then the Loewner chain $(K_t)$ (or the family $(g_t)$) defines
the chordal Schramm-Loewner Evolution with parameter $\kappa$ in
$(\H,x,\infty)$. The chain $K_t$ is generated by the trace $\gamma$, a
continuous process taking values in $\overline\H$, in the following
sense: $\H\setminus K_t$ is the unbounded connected component of
$\H\setminus\gamma_{[0,t]}$.

The trace is a continuous non self-traversing curve. It is a.s. simple
if $\kappa\leq 4$ and a.s. space-filling if $\kappa\geq 8$ (\cite{RS01}). The boundary of a nonsimple $\SLE_\kappa$ ($\kappa>4$) is locally absolutely continuous w.r.t. $\SLE_{\hat\kappa}$, $\hat\kappa=16/\kappa$ (SLE duality, \cite{Dub_dual,Zhan_dual}).


Note that chordal SLE depends only on two boundary points, and radial
SLE depends on one boundary and one bulk point. In several natural
instances, one needs to track additional points on the boundary. This
has prompted the introduction of $\SLE_\kappa(\rho)$ processes in
\cite{LSW_restr}, generalized in \cite{Dub_kapparho}. The driving Brownian motion is
replaced by a semimartingale which has local Girsanov density w.r.t.
the original Brownian motion. These turn out to be technically useful processes (eg \cite{Dub_dual}).

In the chordal case, let $\underline\rho$ be a multi-index, i.e. :
$$\underline\rho\in \bigcup_{i \ge 0} \R^i$$
Let $k$ be the length of $\underline\rho$; if $k=0$, one simply
defines $\SLE_\kappa(\varnothing)$ as a standard $\SLE_\kappa$. If
$k>0$, assume the existence of processes $(W_t)_{t\geq 0}$ and
$(Z^{(i)}_t)_{t\geq 0}$, $i\in\{1\dots k\}$ satisfying the SDEs:
\begin{equation}\label{E1}
\left\{\begin{array}{l}dW_t=\sqrt\kappa dB_t+\sum_{i=1}^k\frac{\rho_i}{W_t-Z^{(i)}_t}dt\\
dZ^{(i)}_t=\frac 2{Z^{(i)}_t-W_t}dt\end{array}\right.
\end{equation}
and such that the processes $(W_t-Z^{(i)}_t)$ do not change sign.
Then we define the chordal $\SLE_\kappa(\underline\rho)$ process starting from
$(w,z_1,\dots z_k)$ as a chordal Schramm-Loewner evolution the driving process
of which has the same law as $(W_t)$ as defined above, with
$W_0=w,Z^{(i)}_0=z_i$.

\subsection{Partition functions}

In this subsection we introduce partition functions of $\SLE$ (a predefinition is in \cite{Dub_dual}), and give some basic properties. These partition functions are null vectors of some canonical Virasoro representations (\cite{Dub_Vir}). They correspond to some correlators in Conformal Field Theory.
 
We begin with an informal discussion to motivate the definition (see eg \cite{BB_review} and references therein for related topics). Consider the Ising model on, say, the triangular lattice. Let $D$ be a (simply connected) portion of the triangular lattice with boundary vertices partitioned in two arcs $\partial^-$, $\partial^+$. A spin configuration $\eps$ consists of an assignment of $\pm$ spins to vertices of $D$, the spins being fixed on the boundary ($\pm$ on $\partial^\pm$). The energy of a configuration is $H(\eps)=-\beta\sum_{i\sim j}\delta_{\eps_i,\eps_j}$. The partition function ${\mc Z}$ is defined as:
$${\mc Z}(D)={\mc Z}(D,\partial^-,\partial^+)=\sum_\eps\exp(-H(\eps)).$$
Except for exceptional cases (torus), there is no explicit asymptotic expansion (as the mesh of the lattice goes to 0) of this.

In this situation, one can define an interface $\gamma$ running between the connected clusters of negative spins attached to $\partial^-$ and positive spins attached to $\partial^+$; it connects the two marked boundary points $x,y$ separating $\partial^-$, $\partial^+$. Consider the following relative situation: $D'$ is another configuration which is identical to $D$ in a neighbourhood $U$ of $x$. The two models induce measures $\mu,\mu'$ on $\gamma^U$, the interface $\gamma$ started from $x$ stopped upon exiting $U$. This defines two new configurations, denoted simply by $D\setminus\gamma^U$, $D'\setminus\gamma^U$, in which the spins neighbouring $\gamma^U$ (which are fixed by construction) are taken as part of the boundary, and the marked point $x$ is moved to the tip of $\gamma^U$. Then it is easy to see that:
$$\frac{d\mu'}{d\mu}(\gamma^U)=\frac{{\mc Z}(D'\setminus\gamma^U){\mc Z}(D)}{{\mc Z}(D'){\mc Z}(D\setminus\gamma^U)}$$
This is only using the local form of the interaction (and the existence of a ``Markovian" set of boundary conditions). If this converges to $\SLE$ (for critical $\beta$), the LHS is well-defined; this suggests looking for continuous analogues of ${\mc Z}$ compatible with Radon-Nikod\`ym derivatives. This is achieved by the

\begin{Def}
Let $c$ be a configuration $c=(D,x,y)$ consisting of a simply connected Riemannian surface $D$ with metric $g$ smooth up to the boundary, two marked boundary points $x,y$ with analytic local coordinates. The partition function ${\mc Z}^{\SLE}_{c,\kappa}$ of chordal $\SLE_\kappa$  is:
$${\mc Z}^{\SLE}_{c,\kappa}={\det}_\zeta(\Lap_D)^{-\frac{\cc}2}H_D(x,y)^{h_{1;2}}$$
where $\cc=1-\frac 32\cdot\frac{(\kappa-4)^2}{\kappa}$ is the central charge, $h_{1;2}=h_{1;2}(\kappa)=\frac{6-\kappa}{2\kappa}$.
\end{Def}

Here 
$$H_D(x,y)=\lim_{x'\rightarrow x,y'\rightarrow y}\frac{G_D(x',y')}{\Im(z_x(x'))\Im(z_y(y'))}$$
 (Poisson excursion kernel, relative to the local coordinates $z_x,z_y$). This can also be seen as a tensor. The local coordinate $z_x$ maps a neighbourhood of $x$ in $D$ conformally to a neighbourhood of $0$ in $\H$.

The following situation will be typical. Let $c_{1}=(D_{1},x_1,y_1)$ be a configuration; $\delta$ a crosscut separating $x_1,y_1$, $C$ a collar neighbourhood of $\delta$ at positive distance of $x_1,y_1$. Let $c_{2}=(D_{2},x_2,y_2)$ be another configuration that agrees with $D_{1}$ in the collar $C$. One can generate hybrid configurations $c_{ij}=(D_{ij},x_i,y_j)$, such that $c_{ij}$ agrees with $D_i$ left of $\delta$ and with $D_j$ right of $\delta$, $i,j\in\{1,2\}$. The local coordinates at $x_i$ of $D_{i1}$, $D_{i2}$ are the same, symmetrically at $y_i$. The metrics of $D_{i1}$, $D_{i2}$ agree to the left of $\delta$ (and a bit further), symmetrically for $D_{1j}$, $D_{2j}$. Then one can form the ratio:
$$\frac{{\mc Z}(c_{11}){\mc Z}(c_{22})}{{\mc Z}(c_{21}){\mc Z}(c_{12})}$$
The point is that this is independent of choices of local coordinates (as tensor dependences cancel out) and of metrics (due to the local form of the Polyakov-Alvarez formula, Proposition \ref{PA}). This is an analogue with boundary of the ``train track" argument of \cite{KontSuh}.

The definition of the partition function is now justified by the following result.

\begin{Prop}
Let $c=(D,x,y)$, $c'=(D',x,y')$ be two chordal configurations that agree in a neighbourhood $U$ of $x$, $U$ at positive distance of $y,y'$. Let $\gamma^\tau$ be the $\SLE$ trace started from $x$ stopped upon exiting $U$ at time $\tau$; $c_\tau=(D\setminus\gamma^\tau,\gamma_\tau,y)$, $c'_\tau=(D'\setminus\gamma^\tau,\gamma_\tau,y')$. Then:
$$\frac{d\mu^{\SLE}_{c'}}{d\mu^{\SLE}_{c}}(\gamma^\tau)=
\frac{{\mc Z}^{\SLE}_{c'_\tau}{\mc Z}^{\SLE}_{c}}
{{\mc Z}^{\SLE}_{c'}{\mc Z}^{\SLE}_{c_\tau}}.
$$
\end{Prop}
\begin{proof}
This is proved in Proposition 3 in \cite{Dub_dual}, based on results in \cite{LSW_restr}. The loop measure term is identified via Proposition \ref{loop_zeta}.
\end{proof}

There is a number of variants of $\SLE$. A configuration $c$ can consist of a Riemannian bordered surface (oriented, otherwise general topology) $D$ with marked points $x_1,\dots,x_n$ on the boundary and $y_1,\dots,y_m$ in the bulk. Analytic coordinates (or merely 1-jets of local coordinates) at the marked points are given. A partition function ${\mc Z}$ is a positive function of such configurations. It has a tensor dependence on analytic coordinates (ie it transforms as $\prod_i(dz_i)^{h_i}\prod_j|dw_j|^{2h_j}$, $z_i$ local coordinate at $x_i$, $w_j$ local coordinate at $y_j$), and depends on the metric as ${\det}_\zeta(\Lap_D)^{-\cc/2}$. The partition function can be seen as a section of a line bundle over a moduli space, as exposed in \cite{FriKal,Kont_arbeit}.

\begin{Def}
Let ${\mc Z}$ be such a partition function. Assume that $h_1=h_{1;2}(\kappa)$. An $\SLE_\kappa({\mc Z})$ is a random non self traversing curve on $D$ started at $x_1$ such that for any simply connected neighbourhood $U$ of $x_1$ in $D$ at positive distance of all other marked points, $\varphi:U\rightarrow V$ a conformal equivalence between $U$ and a bounded neighourhood of $x$ in a simply connected configuration $c=(D',x,y)$, $y$ at positive distance of $V$, one has:
$$\frac{d\varphi_*\mu_D^{\SLE_\kappa({\mc Z})}}{d\mu^{\SLE_\kappa}_{D'}}(\gamma^\tau)=\frac{{\mc Z}(D\setminus\varphi^{-1}(\gamma^\tau)){\mc Z}^{\SLE}(D')}{{\mc Z}(D){\mc Z}^{\SLE}(D'\setminus\gamma^\tau)}$$
where $\gamma^\tau$ is the $\SLE$ trace stopped upon exiting $V$. 
\end{Def}
This does not depend on the choice of $\varphi$, from the previous result. We proceed to show how some variants of $\SLE$ fit in this construction. An important situation is when the same partition function ${\mc Z}$ generates $\SLE$'s starting from different seeds : the two $\SLE$'s then satisfy {\em local commutation} (\cite{Dub_Comm}). This imposes precise conditions on ${\mc Z}$. Definitions of $\SLE$'s in general configurations from CFT correlators are considered in \cite{FriKal}.

In order to express partition functions invariantly, we need to introduce some harmonic constructions:
\begin{itemize}
\item If $D$ is a domain, $y\in D$, $x\in \partial D$, the Poisson kernel $P_D(y,x)$ is a 1-form in $x$ given by $P_D(y,x)dx=\Harm_D(y,dx)$.
\item If $D$ is a domain, $x,y\in\partial D$, the Poisson excursion kernel $H_D(x,y)$ is a 1-form in $x,y$ given by:
$$H_D(x,y)=\partial_{n_y}P_D(y,x)=\partial_{n_x}P_D(x,y)=\partial_{n_x}\partial_{n_y}G_D(x,y)$$
\item If $D$ is a simply connected domain, $y\in D$, $\varphi:D\rightarrow\D$, $y\mapsto 0$, a conformal equivalence, let $H_D(y)=|d\varphi|_{|y}$ ($\varphi$ is unique up to a phase); this is a version of the conformal radius.
\end{itemize}

Let us use chordal $\SLE_\kappa$ in $\H$ as reference measure (normalized at infinity as usual). Let $z_i$'s be marked points (initially) on the real line), and $Z^i_t=g_t(z_i)-W_t$.
Then a simple computation (eg Section 6.1 in \cite{Dub_dual}) shows that:
$$M_t=\prod_i g'_t(z_i)^{\alpha_i}(Z^i_t)^{\beta_i}\prod_{i<j}|Z^j_t-Z^i_t|^{\eta_{ij}}$$
is a local martingale (under the reference chordal measure) if $2\alpha_i=\frac\kappa 2\beta_i(\beta_i-1)+2\beta_i$, $2\eta_{ij}=\kappa\beta_i\beta_j$. Using this as a density produces an $\SLE_\kappa(\underline\rho)$ process, $\underline\rho=\kappa\beta_1,\dots,\kappa\beta_n$. This process is invariant in distribution under homographies if $\rho_1+\cdots+\rho_n=\kappa-6$ (Lemma 3.2 in \cite{Dub_Comm}). In terms of partition functions, this can expressed by:
$${\mc Z}^{\SLE_\kappa(\underline\rho)}_c={\det}_\zeta(\Lap_D)^{-\frac\cc 2}\prod_{0\leq i<j}H_D(x_i,x_j)^{-\frac{\rho_i\rho_j}{4\kappa}}$$
in a configuration $c=(D,x_0,x_1,\dots,x_n)$ with the seed at $x_0$ and the convention $\rho_0=2$.

To treat the radial case, we use chordal results, together with a reflection argument.

The map $g_t:\H\setminus K_t\rightarrow\H$ can be extended by Schwarz reflection: $g_t(\overline z)=\overline{g_t(z)}$. This is compatible with Loewner evolution. In the result above, one can take $z_i$ on the real line or a pair of conjugate $z_i,z'_i=\overline{z_i}$, pairing terms so as to get a real process. For instance, take two marked points $y,y'=\overline y$; set $\rho=\rho'=(\kappa-6)/2$ (to get invariance under homographies). Then the resulting $\SLE$ is simply radial $\SLE_\kappa$ aiming at $y\in\H$.

More generally, mark $x=z_0,z_1,\dots,z_n\in\R$ and $y\in\H$; take $\rho=\rho'=(\kappa-6-\overline\rho)/2$, where $\overline\rho=\rho_1+\cdots+\rho_n$. Also set $\rho_0=2$ (at the seed of the $\SLE$). One gets a local martingale:
$$M_t=|g'_t(y)|^{2\alpha}\Im(Y_t)^{\frac{\rho^2}{2\kappa}}\prod_{i>0} g'_t(z_i)^{\alpha_i}\prod_{i\geq 0}|Y_t-Z^i_t|^{\frac{\rho\rho_i}\kappa}\prod_{0\leq i<j}|Z^j_t-Z^i_t|^{\frac{\rho_i\rho_j}{2\kappa}}$$
where $\alpha_i=\frac{\rho_i}{4\kappa}(\rho_i-\kappa+4)$.
This density generates radial $\SLE_\kappa(\underline\rho)$. A more invariant phrasing can be obtained as follows. In the upper half-plane $\H$, $H_\H(x,y)=\frac{dxdy}{(x-y)^2}$, $P_\H(y,x)=\Im(\frac 1{x-y})dx=\frac{\Im(y)}{|x-y|^2}dx$, $H_\H(y)=\frac{1}{\Im y}|dy|$, up to multiplicative constants. This identifies the partition function as:
$${\mc Z}^{\SLE_\kappa(\underline\rho)}_c={\det}_\zeta(\Lap_D)^{-\frac\cc 2}H_D(y)^{2\alpha}\prod_{i\geq 0}P_D(y,z_i)^{-\frac{\rho\rho_i}{2\kappa}}\prod_{0\leq i<j}H_D(z_i,z_j)^{-\frac{\rho_i\rho_j}{4\kappa}}$$
in the configuration $c=(D,x_1,\dots,x_n,y)$, where $\rho=(\kappa-6-(\rho_1+\cdots+\rho_n))/2$, $\alpha=\frac\rho {4\kappa}(\rho-\kappa+4)$.

As pointed out in \cite{Dub_Comm}, a particularly interesting situation is when $\rho_1=\cdots=\rho_n=2$. In this case, one obtains a symmetric partition function, so that one can grow simultaneously $n$ commuting $\SLE$'s starting from the boundary marked points, aiming at the bulk point. In this case, the partition function is written
$${\mc Z}^{\SLE_\kappa(\underline 2)}_{c}={\det}_{\zeta}(\Lap)^{-\frac\cc 2}H_D(y)^{2h_{0;n/2}}\prod_iP_D(y,x_i)^{-\frac\rho\kappa}\prod_{i<j}H_D(x_i,x_j)^{-\frac 1\kappa},$$
using the ``highest weight" notation:
$$h_{p;q}(\kappa)=\frac{(p\kappa-4q)^2-(\kappa-4)^2}{16\kappa}=h_{q;p}(16/\kappa).$$
We refer to these as multiple radial $\SLE$'s.

\section{Massless Euclidean free field}

In this section, we gather a few facts on the massless (Euclidean) free field that will be needed later. We consider here only fields with Dirichlet boundary conditions. See \cite{Simon_Pphi,Glimm_Jaffe} for background on the free field, \cite{Janson} for Gaussian Hilbert spaces; also the survey \cite{Sheff_GFF}.

\subsection{Discrete free field}

To illustrate some of the notions while avoiding technicalities, we consider first a discrete analogue of the situation (leading to finite dimensional Gaussian vectors).

Let $\Gamma$ be a connected graph with some vertices marked as the boundary $\partial\Gamma$. Fields $\phi$ with Dirichlet boundary conditions are elements of $\sH^1_0(\Gamma)$, that is functions on vertices that vanish on the boundary, with centered Gaussian distribution relative to the Dirichlet inner product:
$$\langle f,g\rangle=\sum_{x\sim y}(f(y)-f(x))(g(y)-g(x))=\sum_{x\in\mathring{\Gamma}}f(\Lap g)(x)$$
where $\Lap$ is the (positive) combinatorial Laplacian:
$$(\Lap g)(x)=\sum_{y\sim x} (g(x)-g(y)).$$
Hence $(\phi(x))_{x\in\mathring{\Gamma}}$ is a Gaussian vector with distribution:
$${\mc Z}_\Gamma^{-1}\exp\left(-\frac 12 \langle \phi,\phi\rangle\right)\prod_{x\in\mathring{\Gamma}}\frac{d\phi(x)}{\sqrt{2\pi}}$$
where the normalization constant ${\mc Z}_\Gamma$ is given by
$${\mc Z}_\Gamma=\det(\Lap)^{-\frac 12}.$$

The free field with boundary conditions $\phi_\partial\in\R^{\partial\Gamma}$ is the Gaussian variable on the affine space $\{\phi\in\R^\Gamma, \phi_{|\partial\Gamma}=\phi_\partial\}$, with covariance operator $\Lap^{-1}$. It is easy to see that the mean $m$ of the field is the harmonic extension of $\phi_\partial$ to $\Gamma$. Furthermore, $\phi$ is distributed as $\phi=m+\phi_0$, where $\phi_0$ is a free field with (zero) Dirichlet boundary conditions. Finally, the partition function can be expressed as:
$${\mc Z}_{\Gamma,\phi_\partial}=\int_{\phi:\phi_{|\partial\Gamma}=\phi_\partial}\exp(-\frac 12\langle \phi,\phi\rangle)\prod_{x\in\mathring{\Gamma}}\frac{d\phi(x)}{\sqrt{2\pi}}=\det(\Lap)^{-\frac 12}\exp(-\frac 12\langle m,m\rangle) $$

One can also put weights on (unoriented) edges and vertices. 
Assume that the vertex set is partitioned in connected subsets $V_l,\delta,V_r$, in such a way that no vertex of $V_l$ is adjacent to $V_r$ (one may also require: no vertex of $\delta$ has all its neighbours in $\delta$). It is easy to see that:
$$\sH^1_0(\Gamma)\simeq \sH^1_0(\Gamma_l)\oplus^\perp W \oplus^\perp \sH^1_0(\Gamma_r)$$
 where $W$ is the space of functions in $\sH^1_0(\Gamma)$ that are harmonic except on $\delta$, $\Gamma_l$ is the graph with inner vertices $V_l$, so that $\delta$ is part of its boundary. Functions in $W$ are in bijection with functions on $\delta$ vanishing on $\delta\cap\partial\Gamma$ (unique harmonic extension to $\Gamma_l,\Gamma_d$). Thus there is an inner product on functions on $\delta$ induced by the inclusion $W\hookrightarrow \sH^1_0(\Gamma)$. Let $P_r$ (resp. $P_l$) denotes the operator from $W$ to functions on $\Gamma_r$ (resp. $\Gamma_l$) that associates to $w\in W$ its unique harmonic extension to $\Gamma_r$ (with Dirichlet boundary condition on $\partial\Gamma$). Define a Neumann jump operator in $\End(W)$ as follows:
$$(Nw)(x)=\sum_{y\in\Gamma_r, y\sim x}(w(x)-(P_rw)(y))+\sum_{y\in\Gamma_l, y\sim x}(w(x)-(P_lw)(y)).$$
If $(Pw)\in H_0^1(\Gamma)$ is the function equal to $P_lw$ (resp. $P_rw$, $w$) on $\Gamma_l$ (resp. $\Gamma_r$, $\delta$), then
$$\langle Pw,Pw\rangle=\sum_{x\in\mathring{\Gamma}}(Pw)(\Lap Pw)(x)=\sum_{x\in\delta}(Pw)(\Lap Pw)(x)=\sum_{x\in\delta} w(Nw)(x)$$
This shows that a free field $\phi$ on $\Gamma$ is the sum of three independent components:
$$\phi=\phi_l+Pw+\phi_r$$
where $\phi_l$ (resp. $\phi_r$) is a free field on $\Gamma_l$ (resp. $\Gamma_r$) with Dirichlet boundary conditions on $\partial\Gamma\cup\delta$ and $w=\phi_{|\delta}$ is a Gaussian variable taking values in $W$ with covariance operator $N^{-1}$ (which is the restriction of $\Lap^{-1}$ to $\delta$). The restrictions of $\phi$ to $\Gamma_l$, $\Gamma_r$ (not to confuse with $\phi_l$, $\phi_r$) are independent conditionally on $w$:
$$\phi_{|\Gamma_l}=\phi_l+P_lw{\textrm\ ,\ \ \ }\phi_{|\Gamma_r}=\phi_r+P_rw.$$ 

Chasing normalizing constants in Gaussian integrals, one also get the identity:
$$\det(\Lap_D)=\det(\Lap_{\Gamma_l})\det(N)\det(\Lap_{\Gamma_r}).$$

Consider now the following situation: $\Gamma_1$, $\Gamma_2$ are graphs as above that agree in a neighbourhood of $\Gamma_r$. 
Let $\mu_i$ be the discrete free field measure on $\sH^1_0(\Gamma_i)$, $i=1,2$; let $R$ be the restriction $\phi\mapsto\phi_{|\Gamma_r}$. We are interested in the Radon-Nikod\`ym derivative:
$$\left(\frac{dR_*\mu_2}{dR_*\mu_1}\right)(\phi_{|\Gamma_r}).$$
From the decomposition $\phi_{|\Gamma_r}=\phi_r+P_rw$ where $\phi_r$ is independent of $w=T\phi$ and its distribution is the same for $\Gamma_1$, $\Gamma_2$, it is clear that:
$$\left(\frac{dR_*\mu_2}{dR_*\mu_1}\right)(R\phi)=
\left(\frac{dT_*\mu_2}{dT_*\mu_1}\right)(T\phi)$$
Looking at the marginal distribution $T\phi$, we may as well assume that $\Gamma_1,\Gamma_2$ only agree in a collar neighbourhood of $\delta$. Since these distributions are Gaussian, we have:
$$\left(\frac{dT_*\mu_2}{dT_*\mu_1}\right)(w)=\frac{\det(N_1)^{1/2}}{\det(N_2)^{1/2}}\exp\left(\frac 12\langle w,(N_1-N_2)w\rangle_{L^2(\delta)}\right)$$
where $N_i$, $i=1,2$, is the jump operator in each situation. We note that:
$$\frac{\det(N_1)}{\det(N_2)}=\frac{\det(\Lap_{\Gamma_1})}{\det(\Lap_{\Gamma_{l,1}})\det(\Lap_{\Gamma_{r,1}})}
\cdot\frac{\det(\Lap_{\Gamma_{l,2}})\det(\Lap_{\Gamma_{r,2}})}{\det(\Lap_{\Gamma_2})}
=\frac{\det(\Lap_{\Gamma_1})\det(\Lap_{\Gamma_{l,2}})\det(\Lap_{\Gamma_{r,2}})}{\det(\Lap_{\Gamma_2})\det(\Lap_{\Gamma_{l,1}})\det(\Lap_{\Gamma_{r,1}})}$$
which is better suited to scaling limits. Also, while $N_1,N_2$ will converge to first order pseudodifferential operators, $N_1-N_2$ will converge to a smoothing kernel operator.

We conclude with a computation of partition functions, that is an elementary discrete analogue of Lemma \ref{loccoupl}. Let $\Gamma_1$, $\Gamma_2$ be graphs that agree in a neighbourhood of a cut $\delta$, with boundary conditions $\phi_{\partial_1}$, $\phi_{\partial_2}$. Let $\Gamma_{ij}$, $i,j\in\{1,2\}$, be the graph that agrees with $\Gamma_i$ (resp. $\Gamma_j$) left (resp. right) of $\delta$, with induced boundary conditions $\phi_{\partial_{ij}}$. Consider the measure on $\phi_{|\delta}$:
$$d\mu_\delta(\phi_{|\delta})=\exp(-\frac 12\sum_{x,y\in\delta, x\sim y}(\phi(x)-\phi(y))^2)\prod_{x\in\mathring{\delta}}\frac{d\phi(x)}{\sqrt{2\pi}}$$
and the one-sided partition function ${\mc Z}_{\Gamma^l_i,\phi_{\partial_i},\phi_{\delta}}$ for the field left of $\delta$ with boundary conditions $\phi_{\partial_i}$ on the part of $\partial_i$ left of $\delta$ and $\phi_\delta$ on the cut $\delta$; ${\mc Z}_{\Gamma^r_i,\phi_{\partial_i},\phi_{\delta}}$ is defined similarly. Chasing definitions, one gets the decomposition
$${\mc Z}_{\Gamma_{ij},\phi_{\partial_{ij}}}=\int {\mc Z}_{\Gamma^l_i,\phi_{\partial_i},\phi_{\delta}}{\mc Z}_{\Gamma^r_j,\phi_{\partial_j},\phi_{\delta}}d\mu_\delta(\phi_{\delta})$$
It is then immediate that ($T$ denotes the restriction to $\delta$):
$$dT_*\mu_{ij}(\phi_{\delta})=\frac{{\mc Z}_{\Gamma^l_i,\phi_{\partial_i},\phi_{\delta}}{\mc Z}_{\Gamma^r_j,\phi_{\partial_j},\phi_{\delta}}}{{\mc Z}_{\Gamma_{ij},\phi_{\partial_{ij}}}}d\mu_\delta(\phi_{\delta})$$
where $\mu_{ij}$ is the measure of the discrete free field in $\Gamma_{ij}$ with boundary conditions $\partial_{ij}$. One deduces the identity:
$$\int \frac{dT_*\mu_{21}}{dT_*\mu_{11}}\cdot\frac{dT_*\mu_{12}}{dT_*\mu_{11}}(\phi_{|\delta})dT_*\mu_{11}(\phi_{|\delta})=\frac{{\mc Z}_{\Gamma_{11},\phi_{\partial_{11}}}{\mc Z}_{\Gamma_{22},\phi_{\partial_{22}}}}{{\mc Z}_{\Gamma_{12},\phi_{\partial_{12}}}{\mc Z}_{\Gamma_{21},\phi_{\partial_{21}}}}$$
which is better suited to scaling limits (see Lemma \ref{loccoupl}). Without additional difficulty, one gets a similar identity when $\delta$ is a region with $m$ connected components in its complement.

\subsection{Continuous free field}

Let $D$ be a bounded planar domain with Jordan boundary (allowing the bounding Jordan arc to have double points). The massless (Euclidean) free field is a random distribution $\phi$, ie a random element of $C^\infty_0(D)'$. It has a Gaussian distribution, with mean 0 and covariance operator $G_D$ (the Green kernel with Dirichlet boundary conditions $G_D$). As in the case of Brownian motion, it is sometimes convenient (if only for psychological reasons) to take as model of the underlying probability space a ``path space". The space $C^\infty_0(D)'$ is usually taken as Wiener space; we will also use a tighter $\sH^{-s}(D)$ for some $s>0$.

Let $\sH^1_0(D)$ be the closure of $C^\infty_0(D)$ for the norm:
$$||f||^2_{\sH^1}=\int_D||\nabla f||^2dA$$
where $dA$ is the Lebesgue measure. We have a Poincar\'e inequality $||f||^2_{L^2}\leq(\lambda_1)^{-1}||f||^2_{\sH^1}$, $\lambda_1$ the lowest eigenvalue of $\Lap$ in $D$ (Dirichlet boundary conditions), so that $||.||_{\sH^1}$ is equivalent to the usual Sobolev norm. One feature of this norm is the conformal invariance: 
$$||f||_{\sH^1_0(D)}=||f\circ\psi||_{\sH^1_0(D')}$$
for $\psi:D'\rightarrow D$ a conformal equivalence.

By duality, one defines $\sH^{-1}(D)$, a space of distribution with norm:
$$||f||_{\sH^{-1}}=\sup_{g\in C^\infty_0(D),||g||_{\sH^1}\leq 1} \langle f,g\rangle $$
where $\langle,\rangle$ is the evaluation of the distribution $f$ against the test function $g$. Since $\langle f,g\rangle_{\sH^1}=\langle f,\Lap g\rangle_{L^2}$ (positive Laplacian) for $f,g\in C^\infty_0(D)$, it follows that:
$$\langle f,g\rangle_{\sH^{-1}}=\langle f,\Lap^{-1} g\rangle_{L^2}$$
where $\Lap^{-1}$ is given by convolution with the Green kernel $G_D$.

One can give a first definition of the free field. Let $(\Omega,{\mc F},P)$ be a probability space carrying a sequence $(\eps_n)_n$ of iid centered, unit variance Gaussian variables, ${\mc F}$ the Borel algebra generated by cylinder events. Let $e_n$ be a Hilbert basis of $\sH^{-1}(D)$. Denote $\phi(e_n)=\eps_n\in L^2(\Omega,P)$. This maps the $e_n$'s isometrically from $\sH^{-1}(D)$ to $L^2(\Omega,P)$, so this can be extended to an isometric embedding $\sH^{-1}(D)\hookrightarrow L^2(\Omega,P)$, denoted by $\phi(.)$.
In the language of Gaussian processes, $L^2(\Omega,P)$ is a Gaussian Hilbert space; it is indexed by the Hilbert space $\sH^{-1}$; it is also a special case of a Gaussian stochastic process, with index set 
$\sH^{-1}$ and covariance function $\rho(f,g)=\langle f,g\rangle_{\sH^{-1}}$. Plainly, $\phi(f)$ is a Gaussian variable for any $f\in\sH^{-1}$, and $\E(\phi(f)\phi(g))=\rho(f,g)=\langle f,g\rangle_{\sH^{-1}}$.

By duality, one can also think of $\sH^1_0(D)$ as the index set, via
$\langle\phi,f\rangle_{L^2}=\langle \phi,\Lap^{-1}f\rangle_{\sH^1}$, where $(\Lap^{-1}f)\in \sH^1_0(D)$ for $f\in \sH^{-1}(D)$.

A more explicit construction goes as follows. Let $(e_n)$ be an orthonormal basis of $\sH^1_0(D)$ consisting of smooth functions. Formally, $\phi=\sum_n\eps_n e_n$, where $(\eps_n)_n$ is a sequence of iid random variables; so that for $f$ a test function, $\phi(f)=\sum\eps_n\langle e_n,f\rangle_{L^2}=\sum\eps_n\langle e_n,\Lap^{-1}f\rangle_{\sH^1}$ has variance $||\Lap^{-1}f||^2_{\sH^1}=||f||^2_{\sH^{-1}}$. We have to determine a space in which this is a.s. convergent. 

For simplicity (and by virtue of the conformal invariance of the $\sH^1$ norm), consider the square $D=[0,1]^2$. For a smooth function $f$ on $D$ vanishing on the boundary, consider the Fourier decomposition:
$$f(x,y)=\sum_{j,k>0}2a_{jk}\sin(\pi jx)\sin(\pi ky).$$ 
Then $||f||^2_{L^2}=\sum_{j,k>0}|a_{jk}|^2$ and $||f||^2_{\sH^1}=\sum_{j,k>0}(j^2+k^2)|a_{jk}|^2$. More generally, for $s\geq 0$, $\sH^s_0(D)$ can be defined via: $||f||^2_{\sH^s}=\sum_{j,k>0}(j^2+k^2)^s|a_{jk}|^2$. 
One can define $\sH^{-s}(D)$ by duality as follows: if $f\in \sH^s_0(D)$, $f$ induces a bounded linear form on $\sH^s_0(D)$ by $g\mapsto \langle f,g\rangle_{L^2}$; let $||f||_{\sH^{-s}}$ be the norm of this bounded operator. The completion of $\sH^s_0(D)$ for this norm is $\sH^{-s}(D)$. In terms of Fourier coefficients, $||f||_{\sH^{-s}}^2=\sum_{j,k>0}(j^2+k^2)^{-s}|a_{jk}|^2$.

Let $(\eps_{jk})$ be iid centered, unit variance Gaussian variables. Define $\phi=\sum_{jk}\eps_{jk}e_{jk}$ where $e_{jk}(x,y)=\frac{2\sin(pi jx)\sin(\pi ky)}{\sqrt{j^2+k^2}}$, in such a way that $e_{jk}$ is an orthonormal basis of $\sH^1_0(D)$. It is easy to see that $h$ converges a.s. in any $\sH^{-s}(D)$, $s>0$ (since $\E(||f||_{\sH^{-s}}^2)=\sum_{j,k}(j^2+k^2)^{s-1}<\infty$: Kolmogorov's One Series Theorem). 

Thus we can use $\sH^{-s}(D)$, $s>0$, as a Wiener space for the free field. Classically (Rellich theorem, also clear here from the Fourier representation), $\sH^{s_2}(D)$ is compactly embedded in $\sH^{s_1}(D)$ for $s_2>s_1$. It follows that the measure on $\sH^{-s}(D)$ is a Radon measure : with probability at least $1-\eps$, $||h||_{\sH^{-s/2}}\leq M(\eps)$, and this ball maps to a compact set of $\sH^{-s}(D)$.

Hence we can take $\Omega=\sH^{-s}(D)$, ${\mc F}$ its (countably generated) Borel algebra, $P$ the measure described above. 
As before, for any $f\in\sH^{-1}$, $\phi(f)$ is defined as an element of $L^2(\Omega,P)$. It is easy to see that ${\mc F}$ is generated by these random variables. This reconstructs the Gaussian Hilbert space indexed by $\sH^{-1}(D)$ from a Radon measure on $\sH^{-s}(D)$.
Moreover, for a fixed $\phi\in\Omega$, $f\mapsto \langle \phi,f\rangle_{L^2}$ defines a bounded linear map on $\sH^s_0(D)$.
 
For a general domain $D$, one can map conformally the square $D_0=[0,1]^2$ to $D$. This maps $\sH^{-s}(D_0)$ bicontinuously to $\sH^{-s}_{loc}(D)$, by standard change of coordinates results for Sobolev spaces (this takes care of both unbounded domains and domains with rough boundaries); and as noted earlier, this preserves the $\sH^1$ norm. So one can take here $\Omega=\sH^{-s}_{loc}(D)$, a Fr\'echet space.

\subsection{Field decompositions, trace}

To prepare the description of the spatial Markov property, we describe
decompositions of a free field in different areas of a domain $D$. More specifically, $\delta$ is a smooth crosscut separating $D$ into two open subdomains $D_l$, $D_r$. Most of what is discussed there works for more general topologies (and also in higher dimension). We have already seen a discrete analogue of the situation. We describe here first a Gaussian space approach, and then a pathwise construction. The main result is that one can define a trace of the free field on $\delta$, a Gaussian variable in $\sH^{-s}(\delta)$.

\subsubsection{Decomposition in Gaussian spaces}

We consider decompositions of free fields, from the point of view of Gaussian spaces. The index set $\sH^1_0(D)$ splits as:
$$\sH^1_0(D)\simeq \sH^1_0(D_l)\oplus^\perp W \oplus^\perp \sH^1_0(D_r)$$
where $W$ is the closure of functions that are harmonic on $D_l\sqcup D_r$, and as above can be identified as a function space on $\delta$. 

The Neumann jump operator is defined as:
$$Nw=\partial_n(P^lw)-\partial_n(P^rw)$$
where: $P^l$ is the Poisson operator extending $w$ to a harmonic function on $D^l$ with Dirichlet boundary condition on $\partial D_l\cap\partial D$; similarly $P^r$ is the Poisson operator on $D^r$; the crosscut $\partial$ is oriented (with $D^l$ to its left) and $\partial_n$ is the normal derivative pointing to $D^l$. We also denote by $P^{lr}$ the harmonic extension of $w$ to $D_l\sqcup D_r$. The Green's formula readily shows that $\int_\delta wNw dl=\int_D |\nabla P^{lr}w|^2 dA$ ($dl$ is the length element on $\delta$ induced by the metric on $D$). 

The Neumann operator $N$ is a first order pseudodifferential operator. Let $H^l(x,y)=\partial_{n_x} P^l(x,y)\propto \partial_{n_x}\partial_{n_y}G_{D_l}(x,y)$ be the Poisson excursion kernel. Up to a multiplicative constant, $H^l-H^r$ is the kernel of the operator $N$. Let $\psi$ be a conformal equivalence from $D^l$ to the upper half-plane $\H$. If $\delta$ is smooth enough ($C^{2+\eps}$), $\psi'$ extends to the boundary. Then $H_l(x,y)=\frac{\psi'(x)\psi'(y)}{(\psi(y)-\psi(x))^2}$ by conformal invariance of the Green's function and explicit computations in $\H$ (this expression is Moebius invariant). This shows that $N(x,y)\asymp (x-y)^{-2}$ at short distance. 

We note that in the unit disk $\D=D(0,1)$, say, the Dirichlet energy of $P_\D w$ (harmonic extension of $w$, $w$ a continuous function on the circle) can be expressed in terms of $w$ as follows:
\begin{align*}
\int_{\D}|\nabla P_\D w|^2dA&=\lim_{r\nearrow 1}\int_{D(0,r)}|\nabla P_\D w|^2dA=\lim_{r\nearrow 1}
\int_{C(0,1)^2}w(x)w(y)dl(x)dl(y)\int_{D(0,r)}\nabla P_\D(x,.).\nabla P_\D(y,.)dA\\
&=\lim_{r\nearrow 1}
\int_{C(0,1)^2}w(x)(w(y)-w(x))dl(x)dl(y)\int_{D(0,r)}\nabla P_\D(x,.).\nabla P_\D(y,.)dA\\ 
&=\lim_{r\nearrow 1}
\int_{C(0,1)^2}\frac 12(w(y)-w(x))^2dl(x)dl(y)\int_{D(0,r)}\nabla P_\D(x,.).\nabla P_\D(y,.)dA\\
&=\int_{C(0,1)^2}\frac 12(w(y)-w(x))^2H_\D(x,y)dl(x)dl(y)
\end{align*}
and this is a conformally invariant expression.
Hence we have:
$$\int_\delta wNw dl=\int (w(x)-w(y))^2(H_l-H_r)(x,y)dl(x)dl(y)\asymp \int \frac{(w(x)-w(y))^2}{(x-y)^2}dl(x)dl(y)\asymp||w||^2_{\sH^{1/2}(\delta)}$$
The last asymptotic follows from the local characterization of Sobolev spaces: for $0<s<1$, in dimension $n$, an element $f$ of $\sH^s$ is an element of $L^2$ such that:
$$\int \frac{(f(x)-f(y))^2}{(x-y)^{n+2s}}dvol(x)dvol(y)<\infty$$
and this quantity gives an equivalent norm (modulo constant functions).

We can construct $w$ from $\phi$ as follows. As noted earlier, the $\langle \phi,f\rangle_{\sH^1}$ are elements of $L^2(\Omega,P)$. Take $(e_n)$ a Hilbert basis of $W$, which is isometrically embedded in $\sH^1_0(D)$; one can choose the $e_n$'s with smooth restriction on $
\delta$. Then consider:
$$T_\delta \phi\stackrel{def}{=}\sum_n\langle e_n,\phi\rangle_{\sH^1}e_n$$
This converges a.s. in any $\sH^{-s}(\delta)=\sH^{1/2-n/2-s}(\delta)$ where $s>0$, $n=1$ (dimension of $\delta$); it follows from the equivalence of norms $\langle .,N.\rangle_{L^2}$ and $||.||_{\sH^{1/2}}^2$. We will give a more explicit construction later on.

This can also be seen from the point of view of Wiener chaos decomposition (isometrically, the Fock space, \cite{Simon_Pphi}, I.4). We have:
$$L^2(\Omega,P)\isom \Gamma(\sH^1_0(D))=\bigoplus_{n\geq 0} \sH^1_0(D)^{\odot n }$$
(complexified, symmetrized tensor algebra). For the probability space associated with $T_\delta \psi$:
$$L^2(\Omega_\delta,P_\delta)\isom\Gamma(W)=\bigoplus_{n\geq 0} W^{\odot n }$$
and the isometric embedding $L^2(\Omega_\delta,P_\delta)\hookrightarrow L^2(\Omega,P)$ is induced by the isometric embedding $W\hookrightarrow \sH^1_0(D)$ (second quantization). This embedding is also positive and preserves 1. It follows that it is induced by a measurable map $T_\delta$, with $P_\delta=(T_\delta)_*P$. (If $\ind_A\in L^2(\Omega_\delta,P_\delta)$ is mapped to $f\in L^2(\Omega,P)$, then $f\geq 0$, $1-f\geq 0$, and $\int f(1-f)dP_\delta=0$ by isometry, hence $f=\ind_{A'}$ for some $A'$).

Note also that $\phi(f)$ is defined for any $f\in \sH^{-1}(D)$. There is a bounded operator, the Sobolev trace, from $\sH^{s}_0(D)$ to $\sH^{s-\frac 12}(\delta)$ for $s>\frac 12$. Applying the transpose of the Sobolev trace operator to an element of $\sH^{-\frac 12}(\delta)$ yields a distribution with support on $\delta$ that belongs to $\sH^{-1}(D)$, hence can be evaluated against $\phi$.

\subsubsection{A pathwise construction}

We discuss here a pathwise construction of the trace of the field on the crosscut $\delta$. The issue is that an instance $\phi$ of the free field  lies in $\sH^{-s}(D)$; the Sobolev trace theorem defines a bounded linear map from $\sH^s(D)$ to $\sH^{s-\frac 12}(\delta)$ for $s>\frac 12$, thus cannot be applied here. But we are only concerned with defining a trace almost everywhere on $\sH^{-s}(D)$. For lightness of notation, we take $D$ bounded, with smooth boundary, and so drop the $loc$ subscripts.

Starting from $\phi\in \sH^{-s}(D)$, we want to define a trace $T\phi$ in some function (or distribution) space on $\delta$, almost everywhere in $\phi$. For simplicity, we will define $T\phi$ in $\sH^{-1}(\delta)$. The topological condition that $\delta$ is a crosscut plays no role here, we can simply assume that $\delta$ is a smooth curve in $\overline D$, possibly intersecting $\partial D$ only at its endpoints. It is parameterized by $t\in [0,1]$: $\delta=(\delta_t)_{t\in [0,1]}$, thus identifying $\sH^{-1}(\delta)$ with $\sH^{-1}([0,1])$.

Let us consider a kernel operator $K:\sH^{-s}(D)\rightarrow \sH^{-1}(\delta)$
$$(Kh)(t)=\int_D K(t,z)h(z)dA(z)$$
where the kernel $K$ is smooth, Markov ($K\geq 0$, $K1=1$), and with finite range $\eps>0$ ($K(t,z)=0$ when $|\delta_t-z|\geq\eps$). It follows that $Kh$ is a smooth function on $\delta$. We want to take $\eps\searrow 0$. For this we need to estimate $||Kh||_{\sH^{-1}(\delta)}$. 

Let $\psi$ be a smooth test function on $\delta$ with compact support (ie vanishing in a neighbourhood of the endpoints), and $f$ be a continuous function, $F=\int_0^t f(s)ds$. Then:
$$||f||_{\sH^{-1}}=\sup_\psi\frac{\int \psi f}{||\psi'||_{L^2}}$$
Since $\int \psi f=-\int \psi' F$, it is easily seen that $\psi'\propto F-\int_0^1F$ is optimal (under the constraint $\int_0^1\psi'=0$). This leads to $||f||_{\sH^{-1}}^2=\int F^2-(\int F)^2$, and after some manipulations:
$$||f||_{\sH^{-1}}^2=\int_0^1(\int_0^s f(u)du\int_s^1 f(v)dv)ds$$
which we now specialise to $f=K\phi$.

First we estimate:
$$\E(f(u)f(v))=\E\left(\int_D K(\delta_u,z)\phi(z)dA(z)\int_D(K(\delta_v,z)\phi(z)dA(z)\right)=\int_{D^2} K(\delta_u,z_1)G_D(z_1,z_2)K(\delta_u,z_2)dA(z_1)dA(z_2)$$
given that $G_D(z_1,z_2)=O(\log|z_1-z_2|)$. One can think of the RHS as drawing $z_1$ (resp. $z_2$) from the distribution $K(\delta_u,.)dA$ (resp. $K(\delta_v,.)dA$) and taking the expectation of $G_D(z_1,z_2)$. Given $z_1$, the probability that $|z_2-z_1|<\eta\eps$ is at worse of order $\eta^2$ (if $\delta_u,\delta_v$ are very close). This gives a contribution of order $\int_0^{\eta\eps}\eps^{-2}rdr=O(|\log\eps|)$ (fixing $\eta=1/10$, say). If $|z_1-z_2|\geq\eta\eps$, one also get a contribution of order $|\log\eps|$. This gives the (crude) uniform estimate: $\E(f(u)f(v))=O(|\log\eps|)$. 

If $|\delta_u-\delta_v|>3\eps$, one has $\E(f(u)f(v))=O(\log|\delta_u-\delta_v|)$. It follows that $\E(||K\phi||_{\sH^{-1}}^2)=O(1)$. 

If $K_1,K_2$ are two smooth Markov kernels as above with range $\eps$, we have the following estimate if $|\delta_u-\delta_v|>3\eps$:
$$\E\left((K_1-K_2)\phi(u)(K_1-K_2)\phi(v)\right)=O(\eps^2/|\delta_u-\delta_v|^2)$$
since $G_D(z_1,z_2)=G_D(\delta_u,\delta_v)+O(\eps^2/|\delta_u-\delta_v|^2)$ for $|z_1-\delta_u|<\eps$, $|z_2-\delta_v|<\eps$. (We use here the fact that $K_i1=1$). Combining with the uniform estimate above, we get:
$$\E(||(K_1-K_2)\phi||^2_{\sH^{-1}})=O(\eps^2|\log\eps|)$$

Consider now a sequence of Markov kernels $K_n$ with range $\eps_n=O(n^{-1-\eta})$ for some $\eta>0$, say. Then:
$$\E(||(K_n-K_{n+1})\phi||_{\sH^{-1}})=O(\sqrt{n^{-2-2\eta}\log(n)})=O(n^{-1-\eta/2})$$
which is summable. It follows that:
$$T_\delta \phi\stackrel{def}{=}\lim_{n\rightarrow\infty}K_n\phi$$
exists a.s. in $\sH^{-1}(\delta)$. For another choice of kernel sequence, one gets a.s. the same element. Since the $K_n$'s are bounded linear operators, $T_\delta$ is Borel measurable.

\begin{Rem} One can proceed similarly if $\delta$ is a boundary arc. Estimates of the Green kernel near the boundary show that the trace of the field on the boundary vanishes a.s., as it should (Dirichlet boundary conditions).
\end{Rem}

Since the $K_n$'s are linear, it appears readily that $T_\delta \phi$ is also Gaussian, with covariance operator given by the restriction of the covariance operator of the free field (that is, $G_D$). Note that while the inverse of $G$ is the Laplacian (a differential operator), the inverse of its restriction to $\delta$ is the Neumann jump operator (a first order pseudodifferential operator, which is nonlocal). This can be checked directly: if $\psi$ is a smooth function on $\delta$, say with compact support, then: $\int_\delta G_D(x,.)\psi(x)dl(x)$ is a continuous function on $D$ that vanishes on $\partial D$, is harmonic on $D\setminus\delta$, and its normal derivative across $\delta$ jumps by $f(x)$ at $x\in\delta$.

The Poisson kernel is smooth, so that $PT_\delta \phi$ defines a harmonic function away from $\delta$. One can recover $\phi_l$ (resp. $\phi_r$) by $\phi_l=\phi_{|D_l}-P^lT_\delta \phi$.

\subsection{Markov property}

Let us consider a free field in a domain $D$ with Dirichlet boundary condition (probability space $(\Omega,{\mc F},P)$, $\Omega=\sH^{-s}(D)$), and $\delta$ a crosscut that splits $D$ into two subdomains $D_l$, $D_r$. (This works for more general topology). We describe here a spatial Markov property of the free field first pointed out by Nelson and Symanzik; we essentially follow \cite{Simon_Pphi} here. We also describe some marginal and conditional distributions that will be needed later on.

Let $U$ be an open subset of $D$. Define ${\mc F}_U$ as the subalgebra of ${\mc F}$ generated by the variables $\langle \phi,f\rangle_{L^2}$, where $f\in\sH^{-1}$ is supported in $U$. For a closed set $K$, ${\mc F}_K=\bigcap_{U{\rm\ open,\ }U\supset K}{\mc F}_U$. In Wiener chaos decomposition, the conditional expectation operator:
$$\E(.|{\mc F}_U):L^2(\Omega,{\mc F},P)\rightarrow L^2(\Omega,{\mc F}_U,P)$$
is generated by the projection of $\sH^1_0(D)$ onto its closed subspace $\sH^1_0(U)$ (this is a contraction).
Similarly, $\E(.|{\mc F}_K)$ corresponds to the projection $\sH^1_0(D)\rightarrow \sH^1_0(D\setminus K)^\perp$.

Consider the trace $T_\delta \phi$ of $\phi$ on $\delta$ and the decomposition: $\phi=\phi_l+PT_\delta h+\phi_r$, $\phi_{|D_l}=\phi_l+P^lT_\delta \phi$.

We have the following description:

\begin{Prop}
\begin{enumerate}
\item ${\mc F}_{D_l}$ (resp. ${\mc F}_{D_r}$, ${\mc F}_{\delta}$) is generated by $\phi_{|D_l}$ (resp. $\phi_{|D_r}$, $T_\delta \phi$).
\item ${\mc F}_{D_l}$ and ${\mc F}_{D_r}$ are independent conditionally on ${\mc F}_{\delta}$.
\item $\phi_l$, $\phi_r$, $T_\delta h$ are independent Gaussian, centered, with covariance $G_{D_l}$, $G_{D_r}$, $(G_D)_{|\delta}$. 
\end{enumerate}
\end{Prop}
\begin{proof}
1. For $f\in C^\infty_0(D_l)$, $\langle\phi,f\rangle_{L^2}=\langle\phi_{|D_l},f\rangle_{L^2}$, so by density ${\mc F}_{D_l}$ is generated by $\phi_{D_l}$. For ${\mc F}_\delta$, notice that a function in $C^\infty_0(\delta)$ induces a distribution on $D$ with support in $\delta$, and this distribution is in ${\mc H}^{-1}$.\\
2. This follows from the representation of $\E(.|{\mc F}_U)$ as $\Gamma(p_U)$ (second quantization), where $p_U$ is the orthogonal projection of ${\mc H}^{-1}(D)$ on the (closure of) the space spanned by distributions with support in $U$. So the statement on independence boils down to:
$$p_{D_l}p_{D_r}=p_{\delta}p_{D_r}$$
In turn this follows from the locality of the inverse covariance (viz. the Laplacian). More precisely, if $f$ has support in $\overline{D_r}$, we have to prove that $p_{D_l}f=p_{\partial}f$. It is enough to see that for $g\in C^\infty_0(D_l)$, $\langle p_{D_l}f,g\rangle_{L^2}=0$. Now:
\begin{align*}
\langle p_{D_l}f,g\rangle_{L^2}&=\langle p_{D_l}f,\Lap g\rangle_{\sH^{-1}}=\langle f,p_{D_l}\Lap g\rangle_{\sH^{-1}}\\
&=\langle f,\Lap g\rangle_{\sH^{-1}}=\langle f,g\rangle_{L^2}=0
\end{align*}
since $\Lap g\in C^\infty_0(D_l)$.\\
3. This is an expression of the orthogonal decomposition :
$$\sH^1_0(D)\simeq \sH^1_0(D_l)\oplus^\perp W \oplus^\perp \sH^1_0(D_r).$$
\end{proof}

\subsection{Absolute continuity}

We begin by recalling some general results (following here \cite{DaPrato_LMS}, see also \cite{Simon_Pphi}, I.6). Let $H$ be a Hilbert space, $Q$ a trace class (symmetric, positive) covariance operator. Denote by $dN_Q$ the centered Gaussian measure on $H$ with covariance $Q$ (it exists since $Q$ is trace class) and by $dN_{m,Q}$ the Gaussian measure with mean $m$, covariance $Q$. Then:

\begin{Prop}
\begin{enumerate}
\item Let $Q$ be a positive, trace class operator; $M$ a symmetric operator such that $Q^{1/2}MQ^{1/2}<1$; and $m\in H$. Then:
\begin{equation}\label{Gint}
\int_H\exp\left(\frac 12 \langle Mh,h\rangle_H+\langle m,h\rangle_H\right)dN_Q(h)=\left[{\det}_F(1-Q^{1/2}MQ^{1/2})\right]^{1/2}\exp\left(\frac 12||(1-Q^{1/2}MQ^{1/2})^{-1/2}Q^{1/2}m||^2\right)
\end{equation}
\item (Cameron-Martin formula). The measures $dN_{m,Q}$, $dN_Q$ are mutually absolutely continuous iff $m\in Q^{1/2}(H)$ (the Cameron-Martin space), in which case:
\begin{equation}\label{CM}
\frac{dN_{m,Q}}{dN_Q}(h)=\exp\left(\langle Q^{-1/2}m,Q^{-1/2}h\rangle_H-\frac 12||Q^{-1/2}m||^2_H\right)
\end{equation}
\item If $Q,R$ are trace class covariance operators, then the measures $dN_Q$, $dN_R$ are mutually absolutely continuous iff $R=Q^{1/2}(1-S)Q^{1/2}$ for some symmetric Hilbert-Schmidt operator $S$. Moreover, if $S$ is trace class, $S<1$, one has the expression:
\begin{equation}\label{RNCov}
\frac{dN_R}{dN_Q}(h)=[{\det}_F(1-S)]^{-1/2}\exp\left(-\frac 12\langle S(1-S)^{-1}Q^{-1/2}h,Q^{-1/2}h\rangle_H\right)
\end{equation}
\end{enumerate}
\end{Prop}
Note that $\langle f,Q^{-1/2}h\rangle_H$ is well-defined for any $f\in \overline{Q^{1/2}(H)}$ (which is $H$ if $\Ker Q=\{0\}$). Indeed, if $f\in Q^{1/2}(H)$, $\langle Q^{-1/2}f,h\rangle_H$ is defined, and this mapping $Q^{1/2}(H)\rightarrow L^2(\Omega,N_Q)$ can be completed given the isometry property:
$$\int_H \langle Q^{-1/2}f,h\rangle_H\langle Q^{-1/2}g,h\rangle_H dN_Q(h)=\langle f,g\rangle_H.$$
Also the properties 2,3 can be combined to give a more general expression:
$$\frac{dN_{p,R}}{dN_{m,Q}}=\frac{dN_{p,R}}{dN_{R}}\cdot\frac{dN_{R}}{dN_{Q}}\cdot\frac{dN_{Q}}{dN_{m,Q}}$$

Let us now specialize this to the free field case. Let us take $H=\sH^{-s}(D)$, $s>0$, with inner product:
$\langle f,g\rangle_{\sH^{-s}}=\langle f,\Lap^{-s}g\rangle_{L^2}$. The covariance operator of the free field w.r.t. $\langle .,.\rangle_{L^2}$ is $G_D$. In terms of $\langle .,.\rangle_{\sH^{-s}}$:
$$\E(\langle f,\phi\rangle_{\sH^{-s}}\langle g,\phi\rangle_{\sH^{-s}})=\E(\langle \Lap^{-s}f,\phi\rangle_{L^2}\langle \Lap^{-s}g,\phi\rangle_{L^2})=\langle \Lap^{-s}f, G_D\Lap^{-s}g\rangle_{L^2}=\langle f,Qg\rangle_{\sH^{-s}}$$
where $Q=\Lap^{-s/2}G_D\Lap^{-s/2}$ is trace class. We can rewrite the expressions above in terms of $\langle,.\rangle_{L^2}$ given:
$$\langle Q^{-1/2}f,Q^{-1/2}g\rangle_{\sH^{-s}}=\langle G_D^{-1/2}f,G_D^{-1/2}g\rangle_{L^2}=\langle f,g\rangle_{\sH^1_0}$$
which is simply saying that while the choice of $\sH^{-s}$ is arbitrary, the Cameron-Martin space $\sH^1_0(D)$ is canonical.

Let us consider the following situation. Two domains $D_1,D_2$ agree in a subdomain containing a crosscut $\delta$. The crosscut splits $D_1$ in $D_{l,1}$, $D_r$ and $D_2$ in $D_{l,1}$, $D_r$ (so that the two domains agree in a neighbourhood of $\overline D_r$. We can define the massless free field in $D_1,D_2$, and then restrict it to $D_r$; in this way, we get absolutely continuous measures. We will need an expression for the Radon-Nikod\`ym derivative. As in the discrete case, the Markov property shows that the 
derivative factors through the trace of the field on $\delta$:
$$\left(\frac{dR_*\mu_2}{dR_*\mu_1}\right)(R\phi)=
\left(\frac{dT_*\mu_2}{dT_*\mu_1}\right)(T\phi)$$
where $\mu_i$ are the free field measures, $R$ is the restriction to $D_r$, $T$ the trace on $\delta$. So we have only to consider $\frac{dT_*\mu_2}{dT_*\mu_1}$; for this purpose it is enough to assume that $D_1,D_2$ agree in a collar neighbourhood $C$ of $\delta$. 

Let $N_i$ be the Neumann jump operator on $\delta\subset D_i$, $i=1,2$. While $N_1,N_2$ are first order pseudodifferential operators (on functions on $\delta$), the difference $N_2-N_1$ is a smoothing kernel. Indeed, $N_2-N_1=N_1(N_1^{-1}-N_2^{-1})N_2$, and $N_i^{-1}=(G_{D_i})_{|\delta}$; so it is enough to see that $G_{D_2}-G_{D_1}$ is smooth on $\delta^2$. This follows from the fact that $G_.(x,y)+\frac{1}{2\pi}\log|x-y|$ is smooth. Alternatively, $G_D(x,y)$ counts Brownian paths from $x$ to $y$ in $D$; decomposing w.r.t. the first exit of the collar $C$, one eliminates the contribution of paths staying in $C$ (viz. $G_C$) which accounts for the singularity. In this fashion, one can represent $(G_{D_2}-G_{D_1})_{|\delta}$ in terms of the Poisson kernel in $C$. It follows that $\langle w,(N_2-N_1)w\rangle_{L^2}$ is defined for all $w\in\sH^{-s}(\delta)$.

It will also be convenient to identify the normalization constant in \eqref{RNCov}; here $1-S=N_1N_2^{-1}$, so that $S=N_1(N_1^{-1}-N_2^{-1})$, a smooth kernel operator on $\delta$. Recall that $m^l(D;K_1,K_2)$ denotes the mass of loops in the loop measure in $D$ that intersect both $K_1$ and $K_2$.

\begin{Lem}\label{denscollar}
If $D_1,D_2$ agree in a collar neighbourhood $C$ of $\delta$, then:
$$\left(\frac{dT_*\mu_2}{dT_*\mu_1}\right)(w)=\exp\left(\frac 12\langle w,(N_1-N_2)w\rangle_{L^2}+\frac 12 (m^l(D_1;\delta,D_1\setminus C)-m^l(D_2;\delta,D_2\setminus C))\right)$$
\end{Lem}
The RHS does not depend on the choice of collar $C$, due to the restriction property of the loop measure.
\begin{proof}
We have to prove that:
$${\det}_F(N_1N_2^{-1})=\exp(m^l(D_2;\delta,D_2\setminus C)-m^l(D_1;\delta,D_1\setminus C)).$$
Given the multiplicative structure of the result, it is enough to prove it for $D_2=C\subset D_1$; one may even assume that $D=D_1$ and $C=D_2$ agree on one side of $\delta$. Taking $K_1=\delta$ and $K_2=\partial C$ a crosscut ``parallel" to $\delta$, we have to prove:
$${\det}_F(N_1N_2^{-1})=\exp(-m^l(D;K_1,K_2))$$
Let $T_{12}$, $T_{21}$ be as in Proposition \ref{loopFred}. Given the result there, we need only prove: $N_1N_2^{-1}=1-T_{21}T_{12}$.


In the reduced case, $D$ is a domain, $K_1,K_2$ are two disjoint crosscuts and $D_l\subset D$ (resp. $D_r\subset D$) is the connected component of $K_1$ in $D\setminus K_2$ (resp. of $K_2$ in $D\setminus K_1$); $N$ (resp. $N_l$) is the Neumann jump operator for $K_1$ in $D$ (resp. in $D_l$). We have $N^{-1}=(G_{D})_{|K_1}$, $N^{-1}_l=(G_{D_l})_{|K_1}$. Clearly, $G_{D}-G_{D_l}$ is positive, and we have the following path representation: let $x$ in $D_l$, $f$ a bounded positive Borel function with support in $\overline{D_l\setminus D_r}$, then 
$$\E_x(\int_{\sigma_1}^\tau f(X_t)dt)=\int_{\overline{D_r\setminus D_l}}(G_D-G_{D_l})(x,y)f(y)dA(y)$$
where $X$ is a Brownian motion (running at speed 2) started at $x$, killed when it hits $\partial D$ at time $\tau$; $\sigma_2$ is the first time it hits $K_2$; and $\sigma_1$ is the first time it hits $K_1$ after $\sigma_2$.  Disintegrating w.r.t. $X_{\sigma_2}$, $X_{\sigma_1}$, we get:
$$(G_D-G_{D_l})(x,y)=\int_{K_2}\Harm_{D_l}(x,z_2)dl(z_2)\int_{K_1}\Harm_{D_r}(z_2,z_1)dl(z_1)G_D(z_1,y)$$
When $x,y$ are on $K_1$, this can be phrased more tersely as:
$$N^{-1}-N_l^{-1}=T_{21}T_{12}N^{-1}$$
which is what we needed.

\end{proof}

\section{Boundary conditions and partition functions}

Quoting from \cite{Simon_Pphi}, ``While we will not use Gaussian variables of mean different from zero, they may well play a role in the future development of the theory". A free field in $D$ with boundary conditions $\phi_{|\partial D}=\phi_\partial$ is written as $\phi=m+\phi_0$ where $\phi_0$ is a free field with Dirichlet boundary conditions and $m$, the mean of the field, is the harmonic extension of $\phi_\partial$ to $D$. 

We will define here appropriate sets of boundary conditions that are continuous in Carath\'eodory-type topologies, and study partition functions of associated free fields. 

\subsection{Domain continuity}

In what follows, we will be primarily interested in the chordal case, in which a configuration $c=(D,x,y)$ consists of a simply connected domain $D$ with two points $x,y$ marked on the boundary. We begin with the case of smooth domains.
 
From the examples of the Temperley coupling and the discrete Gaussian Free Field, it is natural to consider free fields with the following boundary conditions in a configuration $c=(D,x,y)$:
\begin{itemize}
\item on $(xy)$, $\phi=\frac\pi 2 a+b(\pi-\wind(y\rightarrow .))$
\item on $(yx)$, $\phi=-\frac\pi 2 a+b(-\pi+\wind(y\rightarrow.))$
\end{itemize}
where $\wind(y\rightarrow w)$ is the winding of the boundary arc from $y$ to $w$ contained in $(xy)$, $(yx)$ respectively. 
We refer to this set of boundary conditions as {\em $(a,b)$ boundary conditions}. Note that this is in general asymmetric in $x,y$ (jump $+\pi a$ at $x$, $-\pi a-2\pi b$ at $y$). In a configuration, there is a unique harmonic function satisfying the $(a,b)$ boundary conditions. If $c=(\H,x,\infty)$, then 
$h_0(z)=a\arg(z-x)$ and in a general smooth domain $D$, if $\varphi$ is a conformal equivalence $(D,x,y)\rightarrow (\H,0,\infty)$, then
$$h_D=h_\H\circ\varphi-b(\arg(\varphi')-\arg\varphi'(y))$$
Note that $\varphi'$ does not vanish in the simply connected domain $D$, so that there is a single valued branch of $\arg(\varphi')$ in $D$.

One can generalize this to configurations $(D,x_1,\dots,x_n,y)$, with jump $\pi a_i$ at $x_i$, $i\leq n$, and $-\pi\sum_ia_i-2\pi b$ at $y$; we call those {\em (chordal) $(\underline a,b)$ boundary conditions}.

It will be convenient to consider fields with (additive) monodromy. For a domain $D$ with a marked point $y$ in the bulk (puncture), we shall consider the affine space of additively multiply valued functions on $D\setminus\{y\}$ that augment by a fixed quantity (the monodromy) along a counterclockwise circle around $y$ and are locally bounded near $y$.

Given a configuration $c=(D,x_1,\dots,x_n,y)$ consisting of a simply connected domain $D$ with $x_1,\dots,x_n$ marked points on the boundary (in ccwise order), $y$ marked point in the bulk, $\underline a=a_1,\dots,a_n$ a list of parameters, an additively multivalued function $f$ on $D^*=D\setminus\{y\}$ satisfies the {\em $(\underline a,b)$ boundary conditions} if:
\begin{itemize}
\item $f$ increases by $b$ times the winding on the boundary, with  additional jumps of $\pi a_i$ at $x_i$, 
\item $f$ has monodromy $\pi(a_1+\cdots+a_n)+2\pi b$ around $y$,
\item $f(z)=O(1)$ near $y$.
\end{itemize}
There is a unique harmonic function $h_0$ satisfying these conditions, that can be expressed in the unit disk $\D$, $y=0$, as:
$$h_0(z)=b\arg(z)+\sum_i a_i\left(\frac 12\arg(z)-\arg(z-x_i)\right)$$
In a general domain $D$, if $\varphi: D\rightarrow\D$ is a conformal equivalence preserving marked points, $h_D=h_\D\circ\varphi-b\arg\varphi'$ modulo an additive constant.

In the $\SLE$ context, it is necessary to consider domains with rough boundaries. Then the winding of the boundary is no longer defined. However, boundary conditions for the free field intervene only through their harmonic extension. Hence one can use the covariance formula:
$$h_D=h_\D\circ\varphi-b\arg\varphi'$$
where $\varphi: D\rightarrow\D$ is a conformal equivalence preserving marked points, to define $h_D$ in general simply connected domains. 
This is up to an additive constant. If the boundary is rough everywhere and $b\neq 0$, there is no very natural way to fix the constant. On the other hand, it is enough for the boundary to be regular enough in a neighbourhood of, say, $y$ to get an unambiguous definition.

We will mostly concerned with the behaviour of the boundary conditions under deformation of the domain.

\begin{Lem}\label{bcont}
Consider a sequence $(c_n)$ of configurations, $c_n=(D_n,x,x^n_i,y)$, that converges to the configuration $c=(D,x,x_i,y)$ in the following sense: $D_n$ converges to $D$ in the Carath\'eodory topology, $x^n_i$ converges to $x_i$, $i=1,\dots,m$, and there is a smooth boundary arc around $x$ is common to all domains in the sequence. Let $h_c$ be the harmonic extension of $(a,b)$ boundary conditions in $c$. Then:
\begin{enumerate}
\item $h_{c_n}$ converges to $h_c$ uniformly on compact sets of $D$. 
\item If $\mu^{\FF}_c$ is the distribution of the free field in $D$ with $(a,b)$ boundary conditions, $R_U$ the restriction to an open set $U\subset\subset D$, then $(R_U)_*\mu^{\FF}_{c_n}$ converges weakly to $(R_U)_*\mu^{\FF}_{c}$.  
\end{enumerate}
\begin{proof}
Let $\varphi_n$ be the unique conformal equivalence $(\D,0,1)\rightarrow (D_n,y,x)$. Then $h_{c_n}=h_\D\circ\varphi_n-b(\arg\varphi_n'-\arg\varphi'_n(1)+\frac\pi 2)$, where $h_\D$ depends implicitly on the $\phi_n^{-1}(x^n_i)$. Carath\'eodory convergence implies that $\varphi_n$ converges to $\varphi_n$ uniformly on compact sets of $\D$; consequently $\varphi'_n$ also converges uniformly on compact sets. Given the hypothesis on the boundary around $x$, it is not hard to see (using eg the Loewner equations) that $\varphi$ and its derivative converge uniformly in a neighbourhood of $1$. This yields local uniform convergence of $h_{c_n}$.

The weak convergence of $(R_U)_*\mu^{\FF}_{c_n}$ follows from the form of the characteristic functional:
$$\widehat{((R_U)_*\mu^{\FF}_{c_n})}(f)=\int\exp\left(i\langle\phi,f\rangle_{L^2}\right)d\mu^{\FF}_{c_n}(\phi)=\exp\left(i\langle h_{c_n},f\rangle_{L^2}-\frac 12 \langle f,G_{D_n}f\rangle_{L^2}\right)$$
where $f$ runs over $C^\infty_0(U)$.
\end{proof}
\end{Lem}

\subsection{Dirichlet energy}

We study here regularised Dirichlet energy of the harmonic extension of boundary conditions described above. One can think of this as a ground state energy. From the discrete situation, it is natural to define the partition function for the free field in $D$ with boundary condition $\phi_{|\partial D}=\phi_\partial$, for $\phi_\partial$ a smooth (for now) function on $\partial D$ as:
\begin{equation}\label{FFpart}
{\mc Z}^{\FF}_{D,\phi_\partial}={\det}_\zeta(\Lap)^{-\frac 12}\exp\left(-\frac 12\langle m,m\rangle_{\sH^1}\right).
\end{equation} 
The use of the regularized ${\det}_\zeta(\Lap)$ is customary in the physics literature, see eg \cite{Gaw_CFT}; notice that this introduces a metric in addition of the complex structure.

When $\phi_\partial$ is piecewise smooth (with jumps), the Dirichlet energy $\langle m,m\rangle_{\sH^1}$ diverges. We will use another (also customary, see eg \cite{Sonoda}) regularization method, that requires introducing local coordinates (or rather 1-jets) at the marked points where $\phi_\partial$ jumps.

Consider a domain $D$ with smooth boundary, $\phi_\partial$ a piecewise smooth function on $\partial D$ with jumps $\delta_i$ at $x_i$ (say in counterclockwise order, $i=1\dots n$), $m$ its harmonic extension to $D$. Let $z_i$ be an analytic local coordinate at $x_i$ (i.e. $z_i(x_i)=0$, $z_i$ maps a neighbourhood of $x_i$ in $D$ to a neighbourhood of $0$ in $\H$). Define:
$$\langle m,m\rangle_{\sH^1}^{reg}=\lim_{\eps_1\searrow 0,\dots,\eps_n\searrow 0}\left(\int_{\{x\in D, |z_i(x)|\geq \eps\}} |\nabla m|^2dA(x)+\sum_{i=1}^n\frac{\delta_i^2}\pi\log(\eps_i)\right)$$
It is easy to see that this limit exists. There is a simple dependence on the choice of coordinates, that can be expressed by saying that the tensor
$$\exp\left(-\frac 12\langle m,m\rangle_{\sH^1}^{reg}\right)\prod_i(dz_i)^{h_i},$$ 
where $h_i=\delta_i^2/2\pi$, is well defined.

Similarly, for functions with monodromy $2\pi\alpha$ around a bulk point $y$, one can use a local coordinate $w$ at $y$ and define:
$$\langle m,m\rangle_{\sH^1}^{reg}=\lim_{\eps\searrow 0}\left(\int_{\{x\in D, |w(x)|\geq \eps\}} |\nabla m|^2dA(x)+2\pi\alpha^2\log(\eps)\right)$$
so that setting $h_0=\frac{\pi\alpha^2}2$, one defines a tensor:
$$\exp\left(-\frac 12\langle m,m\rangle_{\sH^1}^{reg}\right)|dw|^{2h_0}\prod_i(dz_i)^{h_i}.$$

If $c=(D,x,y)$ is a configuration with smooth boundary, one can define:
$${\mc Z}^{\FF}_{c,(a,b)}={\det}_\zeta(\Lap)^{-\frac 12}\exp\left(-\frac 12\langle m,m\rangle^{reg}_{\sH^1}\right)$$
the partition function for $(a,b)$ boundary conditions. We proceed to evaluating this partition function. 

\begin{Prop}
\begin{enumerate}
\item For a chordal configuration $c=(D,x,y)$, bounded with smooth boundary, we have:
$${\mc Z}^{\FF}_{c,(a,b)}=\lambda{\det}_\zeta(\Lap_D)^{-\frac 12+6\pi b^2}H_D(x,y)^{\frac\pi 2a(2b+a)}$$
where $\lambda$ is a positive constant. More generally, for a chordal configuration $c=(D,x_1,\dots,x_n,y)$, we have:
$${\mc Z}^{\FF}_{c,(\underline a,b)}=\lambda{\det}_\zeta(\Lap_D)^{-\frac 12+6\pi b^2}\prod_{i<j\leq n} H_D(x_i,x_j)^{-\frac{\pi a_ia_j}2}\prod_{i\leq n}H_D(x_i,y)^{\frac\pi 2a_i(2b+\overline a)}$$
where $\overline a=a_1+\cdots+a_n$.
\item For a radial configuration $c=(D,x_1,\dots,x_n,y)$, bounded with smooth boundary, we have:
$${\mc Z}^{\FF}_{c,(\underline a,b)}=\lambda{\det}_\zeta(\Lap_D)^{-\frac 12+6\pi b^2}H_D(y)^{\pi b'(b'-2b)}\prod_iP_D(y,x_i)^{\pi a_ib'}\prod_{i<j}H_D(x_i,x_j)^{-\frac{\pi a_ia_j}2}$$
where $\lambda$ is a positive constant, $\varphi: D\rightarrow \D$, $b'=b+\frac 12\sum_ia_i$.
\end{enumerate}
\end{Prop}
\begin{proof}

1. Let $\varphi:(\D,x_0,y_0)\rightarrow (D,x,y)$ be a conformal equivalence. 
Then it is easy to see that (up to an additive constant):
$$m\circ\varphi=a(\arg(z-y_0)-\arg(z-x_0))+b(2\arg(y_0-z)+\arg(\varphi'))=-a\arg(z-x_0)+(2b+a)\arg(z-y_0)+b\arg(\varphi')$$
We have to compute (taking the natural local coordinates in $\D$):
\begin{align*}
\langle m,m\rangle^{reg}_{\sH^1}
&=
\lim_{\eps\searrow 0}\left(\int_{|z-x|\geq\eps, |z-y|\geq\eps}|\nabla m|^2dA(z)+\pi(a^2+(2b+a)^2)\log(\eps)\right)\\
&=\lim_{\eps\searrow 0}\left(\int_{|z-x_0|\geq\eps/|\varphi'|(x_0), |z-y_0|\geq\eps/|\varphi'|(y_0)}|\nabla m\circ\varphi|^2dA(z)+\pi(a^2+(2b+a)^2)\log(\eps)\right)\\
&=\lim_{\eps\searrow 0}\left(\int_{|z-x_0|\geq\eps, |z-y_0|\geq\eps}|\nabla m\circ\varphi|^2dA(z)+\pi(a^2+(2b+a)^2)\log(\eps)\right)\\
&\hphantom{==}+\pi a^2\log|\varphi'|(x_0)+\pi(2b+a)^2\log|\varphi'|(y_0)
\end{align*}
up to an additive constant. By rotational symmetry, the square terms $\int|\nabla\arg(z-x_0)|^2dA$, $\int|\nabla\arg(z-y_0)|^2dA$ contribute a constant. 
We have:
$$\int_\D(\nabla\arg(z-x_0)).(\nabla\arg(z-y_0))dA(z)=
-\frac 12\int_\D|\nabla(\arg(z-x_0)-\arg(z-y_0))|^2dA(z)+cst$$
(with regularization at $x_0,y_0$). Let $f=\arg(z-x_0)-\arg(z-y_0)$; this is piecewise constant on the boundary, with jumps $\pm\pi$ at $x_0,y_0$. Hence:
$$\int_{|z-x_0|\geq\eps,|z-y_0|\geq\eps}|\nabla f|^2dA=\int_{\dots}|\nabla f^*|^2dA=\int_{\partial(\dots)}f^*\partial_nf^*dl$$
(here $f^*$ is an harmonic conjugate of $f$) and $\partial_nf^*=0$ on the unit circle; besides $f^*=\log|z-x_0|-\log|z-y_0|$, so subtracting divergences leads to:
$$\int_\D(\nabla\arg(z-x_0)).(\nabla\arg(z-y_0))dA(z)=-\pi\log|y_0-x_0|.$$

Besides:
$$\int_\D(\nabla\arg(z-x_0)).(\nabla\arg(\varphi'))dA=\int_\D (\nabla\log|z-x_0|).(\nabla\log|\varphi'|)dA=\lim_{\eps\searrow 0}\int_{\partial(\D\setminus D(x_0,\eps))}\log|\varphi'|\partial_n\log|z-x_0|dl$$
Observe that  on the unit circle $\arg(z-y_0)=-\frac 12\arg(z)+cst$, so that $\partial_n\log|z-x_0|=-\partial_t\arg(z-x_0)=\frac 12$; and $\partial_n\log|z-x_0|=-\eps^{-1}$ on the circle $|z-x_0|=\eps$ (normal derivatives are outward pointing), so that:
$$\int_\D(\nabla\arg(z-x_0)).(\nabla\arg(\varphi'))dA=-\pi\log|\varphi'|(x_0)+\frac 12\int_\U\log|\varphi'|$$
Thus:
$$\langle m,m\rangle^{reg}_{\sH^1}=\pi a(2b+a)(\log|\varphi'|(x_0)+\log|\varphi'|(y_0)+2\log|y_0-x_0|)+b^2(\int_\D|\nabla\log|\varphi'||^2+2\int_\U \log|\varphi'|)+cst
$$

By the Polyakov-Alvarez formula (Proposition \ref{PA}), 
$$\int_\D|\nabla\log|\varphi'||^2+2\int_\U \log|\varphi'|=-12\pi\left(\log{\det}_\zeta(\Lap_D)-\log{\det}_\zeta(\Lap_\D)\right)$$
which concludes the chordal case (with two marked points). The general case is similar.

2. Consider the case of $(\underline a,b)$ boundary conditions on a configuration $(D,x'_i,y')$. Let us consider first the case $D=\D$, with marked point 0. Let $\varphi:(\D,x_i,0)\rightarrow (D,x'_i,y')$ be a conformal equivalence. Then it is easy to see that:
$$m\circ\varphi(z)=-\sum_i a_i\left(\arg(z-x_i)-\frac 12\arg(z)\right)+b(\arg(z)+\arg(\varphi'))=-\sum_i a_i\arg(z-x_i)+b'\arg(z)+b\arg(\varphi')$$
where $b'=b+\frac 12\sum_ia_i$. As before:
$$\langle m,m\rangle_{\sH^1}^{reg}=\langle m\circ\varphi,m\circ\varphi\rangle_{\sH^1}^{reg}+\sum_i \pi a_i^2\log|\varphi'|(x_i)+2\pi (b')^2\log|\varphi'|(0).$$
We have to compute the regularized Dirichlet energy, up to an additive constant. The only new term is:
\begin{align*}
\int_{|z|\geq\eps}(\nabla\arg(z)).(\nabla\arg(\varphi'))dA&=\int_{|z|\geq\eps}(\nabla\log|z|).(\nabla\log|\varphi'|)dA
=\int_{|z|=\eps}\log|\varphi'|\partial_n\log|z|dl+\int_{\U}\log|\varphi'|\partial_n\log|z|dl\\
&=-2\pi\log|\varphi'|(0)+o(1)+\int_{\U}\log|\varphi'|dl
\end{align*}
We get the following expression:
\begin{align*}
\frac 12\langle m,m\rangle^{reg}_{\sH^1}=&-\pi\sum_{i<j}a_ia_j\log|x_i-x_j|+\sum_i  a_ib\pi\log|\varphi'|(x_i)+\frac 12\int_\U\log|\varphi'|dl)\\
&+\pi b'(b'-2b)\log|\varphi'|(0)+\int_\U\log|\varphi'|dl)
+\frac {b^2}2\int_\D|\nabla\log|\varphi'||^2dA+(\sum_i\frac\pi 2a_i^2\log|\varphi'|(x_i))+cst\\
=&-\pi\sum_{i<j}a_ia_j(\log|x_i-x_j|+\frac 12\log|\varphi'|(x_i)+\frac 12\log|\varphi'|(x_j))\\
&+\pi\sum_i  a_ib'\log|\varphi'|(x_i)+\pi b'(b'-2b)\log|\varphi'|(0)
+b^2\left(\int_\U\log|\varphi'|dl+\frac 12\int_\D|\nabla\log|\varphi'||^2dA\right)+cst
\end{align*}
and we conclude by identifying the conformal invariants $P_\D,H_\D$ in the unit disk.

\end{proof}

By comparing with the partition functions of $\SLE$, one obtains the following:

\begin{Thm}\label{PFident}
\begin{enumerate}
\item In the chordal case, the identity 
$${\mc Z}^{\SLE}_{c,\kappa}={\mc Z}^{\FF}_{c,(a,b)}$$
between partition functions of chordal $\SLE_\kappa$ in $c=(D,x,y)$ and the free field in $c$ with $(a,b)$ boundary conditions holds (up to a multiplicative constant), provided that $a=\pm\sqrt{\frac 2{\pi\kappa}}$, $b=a(1-\frac\kappa 4)$. More generally, in a configuration $c=(D,x=x_0,x_1,\dots,x_{n}=y)$, 
$${\mc Z}^{\SLE}_{c,(\kappa,\underline\rho)}={\mc Z}^{\FF}_{c,(\underline a,b)}$$
provided that $a_i=\frac{\eps\rho_i}{\sqrt{2\pi\kappa}}$, $b=\eps\frac{4-\kappa}{\sqrt{8\pi\kappa}}$ for $\eps=\pm 1$ (with the convention $\rho_0=2$).
\item In the radial case, the identity
$${\mc Z}^{\SLE}_{c,(\kappa,\underline\rho)}={\mc Z}^{\FF}_{c,(\underline a,b)}$$
between partition functions of radial $\SLE_\kappa(\underline\rho)$ in $c=(D,z_0,\dots,z_n,y)$ and the free field in $c$ with $(\underline a,b)$ boundary conditions holds (up to a multiplicative constant), provided that $a_i=\frac{\eps\rho_i}{\sqrt{2\pi\kappa}}$, $b=\eps\frac{4-\kappa}{\sqrt{8\pi\kappa}}$ for $\eps=\pm 1$ (with the convention $\rho_0=2$).
\end{enumerate}
\end{Thm}

\begin{proof}
It is merely a matter of matching parameters in the expressions
\begin{align*}
{\mc Z}^{\SLE}_{c,(\kappa,\underline\rho)}&={\det}_\zeta(\Lap_D)^{-\frac\cc 2}\prod_{0\leq i<j}H_D(x_i,x_j)^{-\frac{\rho_i\rho_j}{4\kappa}}\\
{\mc Z}^{\FF}_{c,(\underline a,b)}&={\det}_\zeta(\Lap_D)^{-\frac 12+6\pi b^2}\prod_{0\leq i<j\leq n}H_D(x_i,x_j)^{-\frac{\pi a_ia_j}{2}}\\
\end{align*}
in the chordal case (with the convention $a_n=-(2b+\overline a)$) and
\begin{align*}
{\mc Z}^{\SLE}_{c,(\kappa,\underline\rho)}&={\det}_\zeta(\Lap_D)^{-\frac\cc 2}H_D(y)^{2\alpha}\prod_{i\geq 0}P_D(y,z_i)^{-\frac{\rho\rho_i}{2\kappa}}\prod_{0\leq i<j}H_D(z_i,z_j)^{-\frac{\rho_i\rho_j}{4\kappa}}\\
{\mc Z}^{\FF}_{c,(\underline a,b)}&={\det}_\zeta(\Lap_D)^{-\frac 12+6\pi b^2}H_D(y)^{-2\pi bb'}\prod_iP_D(y,x_i)^{\pi a_ib'}\prod_{i<j}H_D(x_i,x_j)^{-\frac{\pi a_ia_j}2}
\end{align*}
in the radial case, 
where $\rho=(\kappa-6-(\rho_1+\cdots+\rho_n))/2$, $\alpha=\frac\rho {4\kappa}(\rho-\kappa+4)$, $b'=b+\frac 12\sum_ia_i=\frac{4-\kappa+2+\overline\rho}{\sqrt{8\pi\kappa}}=-\frac{\rho}{\sqrt{2\pi\kappa}}$.
\end{proof}

\subsection{Variations of harmonic quantities}

We are considering here a local boundary perturbation of a domain $D$ (growth of a hull at a boundary point) and its effect on various harmonic quantities.

Let $(D_t)_{t\geq 0}$ be a decreasing sequence of domains, $x\in\partial D_0$, so that for any neighbourhood $U$ of $x$, for $t$ small enough, the domains $D_t$ agree outside of $U$. The domains are assumed to be Jordan (the boundary can be parameterized as a continuous, not necessarily simple function); by $x\in\partial D$, we mean a prime end that is a point.

Let $G_t$ be the Green kernel of $D_t$. For any $z,z'$ in $D$, $z,z'\in D_t$ for small $t$, and $G_t(z,z')$ decreases. It follows that $G_0(z,z')-G_t(z,z')$ is positive and harmonic in the two variables in $D_t$. Let $t_n\searrow 0$ and $a_n\nearrow\infty$ such that $a_n(G_0(z,z')-G_n(z,z'))$ has a positive limit for some $x,y\in D$. Then (Harnack principle) $a_n(G_0(.,z')-G_n(.,z'))$ converges to a positive harmonic function in $D$; moreover this function extends continuously to $0$ on the boundary except at $x$. It thus has to be proportional to the Poisson kernel $P_{D,x}$ as a function of $z$; by symmetry, the same is true for the $z'$ variable. Hence for an appropriate choice of $a_n$:
$$\lim a_n (G(z,z')-G_n(z,z'))=P_{D,x}(z)P_{D,x}(z')$$
This argument carries to more general topologies. Let us compute in coordinates for the rest, with the usual $\SLE$ conventions.

In the upper half-plane $\H$, $G(z,z')=-\frac 1{2\pi}\log\frac{z-z'}{z-\overline{z'}}$ and $P_{\H,x}(z)=-\frac 1\pi\Im\frac1{z-x}$. For a family $(\H\setminus K_t)$ corresponding to conformal equivalences $(g_t)$, we get $G_t(z,z')=G(g_t(z),g_t(z'))$ and at $t=0$, for a hull growing at $x$,
$$2\pi\partial_t G(z,z')=\frac 2{(z-x)(z'-x)}-\frac 2{(z-x)(\overline{z'}-x)}=\left(\frac 1{z-x}-\frac 1{\overline{z}-x}\right)\left(\frac 1{z'-x}-\frac 1{\overline{z'}-x}\right)=-4\pi^2P_{\H,x}(z)P_{\H,x}(z')$$

Let $m_t(z)$ be the harmonic function in $(\H\setminus K_t)$ with boundary conditions: $-\frac\pi 2 a+b.(\pi+\wind(\infty\rightarrow.))$ on $(\infty,\gamma_t)$ and $h=\frac\pi 2 a+b.(\pi-\wind(.\rightarrow \infty))$ on $(\gamma_t,\infty)$. Then:
$$m_t(z)=-a\Im\log(g_t(z)-W_t)+b\Im\log(g'_t(z))$$
and
\begin{align*}
d\Im\log(g_t(z)-W_t)=-\pi dP_t(z)&=
\Im\left(\frac 1{g_t(z)-W_t}.\left(\frac 2{g_t(z)-W_t}dt-dW_t\right)\right)-\frac 12 \Im\left(\frac 1{(g_t(z)-W_t)^2}\right)d\langle W\rangle_t\\
&=\pi (\frac\kappa 2-2) P_t'(z)dt+\pi P_t(z)dW_t\\
d\Im\log(g'_t(z))&=-\Im\left(\frac 2{(g_t(z)-W_t)^2}\right)=2\pi P'_t(z)dt
\end{align*}
where $P_t=P_{\H\setminus K_t,\gamma_t}$ and $P'_t=\frac\partial{\partial W_t}P_t$ (this is somewhat dependent on the Loewner convention).

\section{Couplings of SLEs and free fields}

\subsection{Local invariance of the free field under SLE dynamics}

We have obtained partition function identities between (versions of) $\SLE$ on the one hand and the free field (with corresponding boundary conditions) on the other hand. We now show that this implies local identities in distribution between $\SLE$ and the free field, in a way closely analogous to local commutation statements (between two $\SLE$'s with the same partition function) in \cite{Dub_Comm}.

Let $c=(D,x_1,\dots,x_n,y)$ be a configuration ($y$ in the bulk), $(\underline a,b)$ a set of boundary conditions for the free field corresponding to $\SLE_\kappa(\underline\rho)$, as in Theorem \ref{PFident}. (This covers the chordal case, when the monodromy around $y$ is 0). Among the marked points on the boundary, $x_1,\dots,x_m$ are seeds of $\SLE$ ($\rho_i=2$ for $i=1,\dots, m$). We denote by $\mu^{\SLE}_c$, $\mu^{\FF}_c$ the respective distributions of the $\SLE$ system and the free field in $c$.

Up to now, we have considered boundary conditions for the free field up to an additive constant. We now need to fix this constant (in order to compare fields in different domains). For instance, one can require that $\phi(x_n^+)=0$, where $x_n$ is not a seed.

Let $U$ be a connected open subset of $D$, not having any seed $x_i$ on its boundary. Let $\tau_i$ be stopping times for each $\SLE$ such that $\gamma_i^{\tau_i}$ is a.s. at distance at least $\eta>0$ of $U$. Let $s_1,\dots, s_m$ be time parameters for each $\SLE$ (these are somewhat arbitrary, up to bicontinuous time change). Then $c_{\underline s}$ is the configuration $(D\setminus\cup_i\gamma_i^{\tau_i},\gamma_{1,\tau_1},\dots,\gamma_{m,\tau_m},x_{m+1},\dots,y)$.

We have the following:

\begin{Lem}\label{locstab}
In the above situation, the following identity of distributions on $C_0^\infty(U)'$ holds:
$$\int d\mu^{\SLE,\underline\tau}_c(\gamma_1^{\tau_1},\dots,\gamma_m^{\tau_m})(R_U)_*\mu^{\FF}_{c_{\underline\tau}}=(R_U)_*\mu^{\FF}_c$$
where $R_U$ denotes the restriction from a subdomain of $D$ (containing $U$) to $U$.
\end{Lem}

In words: run the $i$-th $\SLE$ strand to time $\tau_i$; this generates a random configuration $c_{\underline\tau}$; sample the free field in $c_{\underline\tau}$, conditionally independently; restrict this field from the random domain $D\setminus(\cup_i\gamma_i^{\tau_i})$ to the fixed domain $U$. Then the resulting mixture of Gaussian fields is again Gaussian, identical in distribution to the restriction of the free field in $D$ to $U$.

\begin{proof}

A distribution $\nu$ on $C_0^\infty(U)'$ is determined by the characteristic functional:
$$\hat\nu(f)=\int\exp\left(i\langle\phi,f\rangle_{L^2}\right)d\nu(\phi),$$
where $f$ runs over $C_0^\infty(U)$. By general Gaussian properties,
$$\widehat{((R_U)_*\mu^{\FF}_{c_{\underline s}})}(f)=\int\exp\left(i\langle\phi,f\rangle_{L^2}\right)d\mu^{\FF}_{c_{\underline s}}(\phi)=\exp\left(i\langle m_{\underline s},f\rangle_{L^2}-\frac 12 \langle f,G_{\underline s}f\rangle_{L^2}\right)$$
where $\mu^{\FF}_{c_{\underline s}}$ has mean $m_{\underline s}$ and covariance $G_{\underline s}$. Therefore we have to prove that:
$$\int \exp\left(i\langle m_{\underline \tau},f\rangle_{L^2}-\frac 12 \langle f,G_{\underline \tau}f\rangle_{L^2}\right)d\mu^{\SLE,\underline\tau}_c(\gamma_1^{\tau_1},\dots,\gamma_n^{\tau_n})=\exp\left(i\langle m_{\underline 0},f\rangle_{L^2}-\frac 12 \langle f,G_{\underline 0}f\rangle_{L^2}\right)$$
for all $f\in C^\infty_0(U)$. Since this has to hold also for all stopping times lesser than $\tau_1,\dots,\tau_m$, we have to prove that:
$$(s_1,\dots,s_m)\longmapsto \exp\left(i\langle m_{\underline s},f\rangle_{L^2}-\frac 12 \langle f,G_{\underline s}f\rangle_{L^2}\right)$$
is a martingale in $s_j$ (stopped at $\tau_j$), the other times being fixed. Due to the (joint) Markov property, we can assume that the other times are 0; we have a single time parameter $t=s_j$. Hence we are left with a stochastic calculus problem.

At this point it is rather convenient to compute in coordinates. So one can assume that $D=\H$, the marked point $x_n$ (where the field is $0$) is at infinity. In this situation, the mean $m_t$ of the field in $c_t$ is given by:
$$m_t(z)=-\sum_i a_i\arg(g_t(z)-g_t(x_i))+b'\left(\arg(g_t(z)-g_t(y))+\arg(g_t(z)-\overline{g_t(y)})\right)-b\arg g'_t(z)$$
and this has to be martingale (stopped at positive distance of $z$), under the $\SLE_\kappa(\underline\rho)$ evolution (conventionally, $g_t(x_j)=W_t$, the driving process). This can be checked directly; we give another argument that avoid computations.

We introduced radial $\SLE_\kappa(\underline\rho)$ by means of the local martingale (w.r.t. the reference measure, that is chordal $\SLE$ in $(\H,0,\infty)$):
$$M_t=\prod_i g'_t(x_i)^{\alpha_i}(g_t(x_i)-W_t)^{\beta_i}\prod_{i<j}|g_t(x_i)-g_t(x_j)|^{\eta_{ij}}$$
where two of the marked points are conjugates $y,\overline y$, for appropriate coefficients. We can perturb the above situation (radial $\SLE_\kappa(\underline\rho)$) by adding a marked point $z$ with weight $\rho_z=\kappa\eps$, say. 
We compute:
\begin{align*}
\left(\frac{\partial_\eps(M^\eps_t)}{M^\eps_t}\right)_{|\eps=0}=&(1-\frac\kappa 4)\log(g'_t(z))+\frac{\kappa-4-2\rho}4\log|g'_t(y)|+\sum_{i\geq 0}\frac{\rho_i}{2}\log(g_t(z)-g_t(x_i))\\
&+\frac{\rho}{2}\left(\log(g_t(z)-g_t(y))+\log(g_t(z)-\overline{g_t(y)})\right)
\end{align*}
which thus a (local) martingale for radial $\SLE_\kappa(\underline\rho)$'s. 
Taking the imaginary part, one gets a process proportional to $m_t(z)$.

This proves that $t\mapsto \langle m_t,f\rangle_{L^2}$ is a martingale (stopped away from the support of $f$). There is only one term in $m_t$ with quadratic variation, so one computes easily:
$$dm_t(z)=-\pi a_jP_{\H}(g_t(z),W_t)\sqrt\kappa dB_t$$
where $P$ is the Poisson kernel. Besides, we have computed that:
$$dG_t(z_1,z_2)=-2\pi P_\H(g_t(z_1),W_t)P_\H(g_t(z_2),W_t)dt$$
so if we define ${\mc E}_t=\exp\left(i\langle m_t,f\rangle_{L^2}-\frac 12\langle f,G_t,f\rangle_{L^2}\right)$, we get:
$$\frac{d{\mc E}_t}{{\mc E}_t}=idm_t(f)-\frac \kappa 2 \pi^2a_j^2 \left(\int _{\H}f(z)P_{\H}(g_t(z),W_t)dA(z)\right)^2dt-\frac 12 \langle f,dG_t f\rangle_{L^2}dt$$
which is thus a martingale, given that $a_j^2=\frac{2}{\pi\kappa}$.
\end{proof}

\begin{Rem}
The (pointwise) first moment martingale $m_t(z)$ was pointed out by Sheffield in the context of the free field, in the chordal case. The proof shows that ${\mc E}$ is, up to a multiplicative constant, the exponential martingale of $m$ (see e.g. \cite{RY}).
\end{Rem}

As in \cite{Dub_dual}, this can used to construct ``local couplings". Define
$$\ell(\gamma^{\underline s},\phi)=\frac{d(R_U)_*\mu^{\FF}_{c_{\underline s}}}{d(R_U)_*\mu^{\FF}_{c_{\underline 0}}}(\phi_{|U}),$$
a function in $\gamma^{\underline s}=\gamma_1^{s_1},\dots,\gamma_m^{s_m}$ ($\SLE$ strands) and $\phi$ (a field). Then:
\begin{align*}
\int \ell(\gamma^{\underline s},\phi)d\mu^{\FF}_c(\phi)=&1&\textrm { for all }\gamma^{\underline s}\\
\int \ell(\gamma^{\underline s},\phi)d\mu^{\SLE}_c(\gamma^{\underline s})=&1&\textrm { a.e. in }\phi
\end{align*}
The first line is obvious while the second one is a rephrasing of the lemma. This shows that $\ell.\mu^{\SLE}_c\otimes\mu^{\FF}_c$ is a coupling of $\mu^{\SLE}_c,\mu^{\FF}_c$ (this builds on the fact that we have a coupling restricted to $\gamma^{\underline\tau}$, $\phi_{|U}$, that extends to a coupling of $\gamma,\phi$ using the Markov property of the $\SLE$ system and of the free field).

Let us analyze the density $\ell$. Let $\delta_i$ be crosscuts around $x_i$, $i=1,\dots,m$, that are at positive distance of each other and of all marked points. Let $\delta=\sqcup_i\delta_i$ be the union of crosscuts and $U$ the connected component of $D\setminus \delta$ with no seed $x_i$, $i=1,\dots,m$, on its boundary. From the Markov property of the free field, it readily appears that:
$$ \frac{d(R_U)_*\mu^{\FF}_{c_{\underline s}}}{d(R_U)_*\mu^{\FF}_{c_{\underline 0}}}(\phi_{|U})= \frac{d(T_\delta)_*\mu^{\FF}_{c_{\underline s}}}{d(T_\delta)_*\mu^{\FF}_{c_{\underline 0}}}(T_\delta\phi)$$
i.e. the density factors through the trace on $\delta$, which we now denote simply by $T$. Given the multiplicative identity:
$$\frac{dT_*\mu^{\FF}_{c_{\underline s}}}{dT_*\mu^{\FF}_{c_{\underline 0}}}=
\frac{dT_*\mu^{\FF}_{c_{s_1,0,\dots}}}{dT_*\mu^{\FF}_{c_{\underline 0}}}\cdot
\frac{dT_*\mu^{\FF}_{c_{s_1,s_2,0,\dots}}}{dT_*\mu^{\FF}_{c_{s_1,0,\dots}}}\cdots\frac{dT_*\mu^{\FF}_{c_{s_1,\dots,s_n}}}{dT_*\mu^{\FF}_{c_{s_1,\dots,s_{n-1},0}}}
$$
we can vary time indices one at a time. Let us focus on $\frac{dT_*\mu^{\FF}_{c_{s,0,\dots}}}{dT_*\mu^{\FF}_{c_{\underline 0}}}$, setting $s_1=s$. Denote by $\bar\mu^{FF}_c$ the zero mean free field in $D$, and $c_s=c_{s,0,\dots,0}$, $c=c_{\underline 0}$. We can decompose:
$$\frac{dT_*\mu^{\FF}_{c_{s}}}{dT_*\mu^{\FF}_{c}}=\frac{dT_*\mu^{\FF}_{c_{s}}}{dT_*\bar\mu^{\FF}_{c_{s}}}\cdot
\frac{dT_*\bar\mu^{\FF}_{c_{s}}}{dT_*\bar\mu^{\FF}_{c}}\cdot
\frac{dT_*\bar\mu^{\FF}_{c}}{dT_*\mu^{\FF}_{c}}$$

Then the Cameron-Martin formula yields:
\begin{align*}
\frac{dT_*\mu^{\FF}_{c_{s}}}{dT_*\bar\mu^{\FF}_{c_{s}}}(w)&=\exp\left(\langle w,N_{s}Tm_{s}\rangle-\frac 12\langle Tm_{s},N_sTm_s\rangle\right)\\
\frac{dT_*\mu^{\FF}_{c}}{dT_*\bar\mu^{\FF}_{c}}(w)&=\exp\left(\langle w,NTm\rangle-\frac 12\langle Tm,NTm\rangle\right)
\end{align*}
The middle term was analyzed in Lemma \ref{denscollar}:
$$\frac{dT_*\bar\mu^{\FF}_{c_{s}}}{dT_*\bar\mu^{\FF}_{c}}(w)=\exp\left(-\frac 12\langle w,(N_s-N)w\rangle+\frac 12m^l(D;\gamma_1^s;\delta)\right)$$
so that:
$$\frac{dT_*\mu^{\FF}_{c_{s}}}{dT_*\mu^{\FF}_{c}}=\exp\left(-\frac 12\langle w,(N_s-N)w\rangle+\langle w,N_sTm_s-NTm\rangle+\frac 12m^l(D;\gamma_1^s;\delta)+\frac 12\langle Tm,NTm\rangle-\frac 12\langle T_sm_s,N_sTm_s\rangle\right)$$

We denote by $N_{\underline s}$ the Neumann jump operator on $\delta$ relative to the configuration $c_{\underline s}$, in such a way that $\langle w,Nw\rangle_{L^2(\delta)}$ is the Dirichlet energy of the harmonic extension of $w$ to $D$ (with zero boundary condition on $\partial D$). Observe that when different time indices evolve, the configuration varies in distinct connected components of $D\setminus\delta$; thus
$$(N_{\underline s}-N)(w_1,\dots,w_m)=\left((N_{s_1}-N)(w_1),\cdots,(N_{s_m}-N)(w_m)\right)$$
where $N=N_{\underline 0}$, $N_{s_i}=N_{0,\dots,0,s_i,0,\dots}$, and $w=(w_1,\dots,w_m)$ identifies $L^2(\delta)$ to $\bigoplus^\perp_iL^2(\delta_i)$. So $N_{\underline s}$ has no cross dependence in the $\underline s$ parameters. We proceed to show that this is also the case for $N_{\underline s}Tm_{\underline s}$.

Observe that $Nw=\sum_i\partial_nP_iw$, where $P_i$ is the harmonic extension to the $i$-th connected component of $D\setminus\delta$ and $\partial_n$ is outward pointing (on $\delta$). Consider $\tilde m$ the harmonic function on $D\setminus\delta$ that agrees with $m$ on $\partial D$ and vanishes on $\delta$. Then $m-\tilde m$ is harmonic on $D\setminus\delta$, agrees with $m$ on $\delta$ and vanishes on $\partial D$. Thus $m-\tilde m=PTm$. Now $\sum_i\partial_n m=0$, since $\partial_n m$ on $\delta_i$ is counted once in each direction (and $m$ is smooth across $\delta$). Thus $NTm=-\sum_i\partial_n \tilde m$. In the varying situation, $\tilde m_{\underline s}$ depends only on $s_i$ in the connected component having $\gamma_i$ in its boundary; this is due to the local character of the boundary condition. Hence $N_{\underline s}Tm_{\underline s}$ has no cross dependence in the $\underline s$ parameters. More precisely, $(NTm)_{s_i}-(NTm)$ vanishes on $\delta_j$, $j\neq i$.

Since $N_{\underline s}$ has no cross dependence in the ${\underline s}$ parameters, one gets:
\begin{align*}
N_{s_1,\dots,s_n}-N&=(N_{s_1}-N)+(N_{s_1,s_2}-N_{s_1})+\cdots+(N_{s_1,\dots,s_{n}}-N_{s_1,\dots,s_{n-1}})\\
&=(N_{s_1}-N)+(N_{s_2}-N)+\cdots+(N_{s_n}-N)
\end{align*}
and a similar identity holds for $(NTm)_{\underline s}$. This shows that:
\begin{align*}
\frac{dT_*\mu^{\FF}_{c_{s_1,\dots,s_n}}}{dT_*\mu^{\FF}_{c}}&=\frac{dT_*\mu^{\FF}_{c_{s_1}}}{dT_*\mu^{\FF}_{c}}\cdot\frac{dT_*\mu^{\FF}_{c_{s_1,s_2}}}{dT_*\mu^{\FF}_{c_{s_1}}}\cdots\frac{dT_*\mu^{\FF}_{c_{s_1,\dots,s_n}}}{dT_*\mu^{\FF}_{c_{s_1,\dots,s_{n-1}}}}\\
&=\lambda\cdot\frac{dT_*\mu^{\FF}_{c_{s_1}}}{dT_*\mu^{\FF}_{c}}\cdot\frac{dT_*\mu^{\FF}_{c_{0,s_2}}}{dT_*\mu^{\FF}_c}\cdots\frac{dT_*\mu^{\FF}_{c_{0,\dots,0,s_n}}}{dT_*\mu^{\FF}_c}\\
\end{align*}
where $\lambda$ depends only on $\gamma$ (not $w$). We now study this term.
Recall that $\langle Tm,NTm\rangle$ is the Dirichlet energy of $PTm$. Consider as above $\tilde m$, which is harmonic on $D\setminus\delta$, agrees with $m$ on $\partial D$ and vanishes on $\delta$. Then:
$$\langle \tilde m,\tilde m\rangle^{reg}_{\sH^1}=
\langle m,m\rangle^{reg}_{\sH^1}+\langle PTm,PTm\rangle^{reg}_{\sH^1}$$
This is easily seen for smooth boundary conditions (vanishing in a neighbourhood of the endpoints of the crosscuts), since $m=\tilde m+PTm$ and 
$$\int_D (\nabla m).(\nabla PTm)dA=\sum_i\int_{\partial D_i}m\partial_nmdl=0$$
where the $D_i$ 's are the connected components of $D\setminus\delta$; $m\partial_nm$ is counted twice with opposite signs on $\partial_i$. By approximation, one gets the result for the sets of boundary conditions under consideration. 

When $s_i$ varies, $\tilde m$ changes only in the connected component having $x_i$ on its boundary. It follows that:
$$\langle Tm_{s_1,s_2,\dots},N_{s_1,s_2,\dots}Tm_{s_1,s_2,\dots}\rangle-\langle Tm_{0,s_2,\dots},N_{0,s_2,\dots}Tm_{0,s_2,\dots}\rangle
+\langle m_{s_1,s_2,\dots},m_{s_1,s_2,\dots}\rangle^{reg}_{\sH^1}-
\langle m_{0,s_2,\dots},m_{0,s_2,\dots}\rangle^{reg}_{\sH^1}
$$
depends only on $s_1$ and not on $s_2,\dots,s_{n-1}$ (it is the variation of the Dirichlet energy of $\tilde m$ in the connected component of $x_1$). 

The loop measure term is handled as follows. We have (Proposition \ref{loop_zeta}):
$$\exp(-m^l(D;\gamma_1^s;\delta))=\frac{\det(\Lap)_D\det(\Lap)_{D\setminus(\gamma_1^s\cup\delta)}}
{\det(\Lap)_{D\setminus\delta}\det(\Lap)_{D\setminus\gamma_1^s}}$$
so that
\begin{align*}
\exp(-m^l(D;\gamma_1^{s_1};\delta)-m^l(D\setminus\gamma_1^{s_1};\gamma_2^{s_2},\delta)-\cdots)
&=\frac{\det(\Lap)_D\det(\Lap)_{D\setminus(\gamma_1^s\cup\delta)}}
{\det(\Lap)_{D\setminus\delta}\det(\Lap)_{D\setminus\gamma_1^{s_1}}}\cdot
\frac{\det(\Lap)_{D\setminus\gamma_1^{s_1}}\det(\Lap)_{D\setminus(\gamma_1^{s_1}\cup\gamma_2^{s_2}\cup\delta)}}
{\det(\Lap)_{D\setminus(\gamma_1^{s_1}\cup\delta)}\det(\Lap)_{D\setminus(\gamma_1^{s_1}\cup\gamma_2^{s_2})}}\cdots
\\
&=\frac{\det(\Lap)_D\det(\Lap)_{D\setminus(\gamma^{\underline s}\cup\delta)}}
{\det(\Lap)_{D\setminus\delta}\det(\Lap)_{D\setminus\gamma^{\underline s}}}
\end{align*}
Recall that:
$${\mc Z}^{\FF}_{\underline s}=\det(\Lap)_{D\setminus\gamma^{\underline s}}\exp(-\frac 12\langle m_{\underline s},m_{\underline s}\rangle^{reg}_{\sH^1})$$
Putting things together, we get:
$$\lambda=\frac{{\mc Z}^{\FF}_0}{{\mc Z}^{\FF}_{\underline s}}\cdot
\frac{{\mc Z}^{\FF}_{s_1}}{{\mc Z}^{\FF}_0}\cdots\frac{{\mc Z}^{\FF}_{s_m}}{{\mc Z}^{\FF}_0}$$
where ${\mc Z}^{\FF}_{s_i}={\mc Z}^{\FF}_{0,\dots,0,s_i,0,\dots}$. Using the identity ${\mc Z}^{\FF}={\mc Z}^{\SLE}$ (Theorem \ref{PFident}), this translates into:
$$\lambda=\frac{{\mc Z}^{\SLE}_0}{{\mc Z}^{\SLE}_{\underline s}}\cdot
\frac{{\mc Z}^{\SLE}_{s_1}}{{\mc Z}^{\SLE}_0}\cdots\frac{{\mc Z}^{\SLE}_{s_m}}{{\mc Z}^{\SLE}_0}$$
It follows that at the stopping times $\tau_i$:
$$\lambda(\gamma_1^{\tau_1},\dots,\gamma_m^{\tau_m})=\frac{d\mu^{\SLE}_c(\gamma_1^{\tau_1})\dots d\mu^{\SLE}_c(\gamma_m^{\tau_m})}{d\mu^{\SLE}_c(\gamma_1^{\tau_1},\dots,\gamma_m^{\tau_m})}$$
The measure in the denominator is the measure induced on the stopped paths $\gamma_i^{\tau_i}$ by the $\SLE$ system; the measure in the numerator has the same $m$ marginals, but these are independent.

\begin{Lem}\label{loccoupl}
The measure $\ell.\mu^{\SLE}_c\otimes\mu^{\FF}_c$, where:
$$\ell(\gamma,\phi)=\ell(\gamma^{\underline\tau},T\phi)=\frac{dT_*\mu^{\FF}_{c_{\underline\tau}}}{dT_*\mu^{\FF}_{c}}(T\phi)$$
is a coupling of $\mu^{\SLE}_c,\mu^{\FF}_c$ such that $\gamma_1^{\tau_1},\dots,\gamma_m^{\tau_m}$ are jointly independent conditionally on $T\phi$. More precisely:
$$\ell(\gamma^{\underline\tau},\phi)d\mu^{\SLE}_c(\gamma^{\underline\tau})=\prod_i \left(\frac{d(T_{\delta_i})_*\mu^{\FF}_{c_{\tau_i}}}{d(T_{\delta_i})_*\mu^{\FF}_{c}}(T_{\delta_i}\phi)\right)d\mu^{\SLE}_c(\gamma_i^{\tau_i})
$$
\end{Lem}

Informally, all the interaction between the different $\SLE$ strands is carried by the field. This depends on the fact that the $\SLE$ strands are occulted from each other by the crosscuts $\delta_i$. Note that this conditional independence property is trivially satisfied in a coupling where the $\SLE$'s are determined by the field. The lemma hinges on and contains the following Gaussian integral evaluation:
$$\int\prod_i \left(\frac{d(T_{\delta_i})_*\mu^{\FF}_{c_{\tau_i}}}{d(T_{\delta_i})_*\mu^{\FF}_{c}}(T_{\delta_i}\phi)\right)d(T_\delta)_*\mu^{\FF}_c(T_\delta\phi)=\frac{{\mc Z}^{\FF}_{\underline \tau}}{{\mc Z}^{\FF}_0}\cdot
\frac{{\mc Z}^{\FF}_0}{{\mc Z}^{\FF}_{\tau_1}}\cdots\frac{{\mc Z}^{\FF}_0}{{\mc Z}^{\FF}_{\tau_m}}$$

In the discrete setting, we gave an elementary version of the computation above, thinking of partition functions as (matrix elements of) transfer operators.

We have established Lemmas \ref{locstab}, \ref{loccoupl} in the case where the $\SLE$ strand is absolutely continuous w.r.t. chordal $\SLE$ in a neighbourhood of its starting point $x$. However, it will be also useful to consider versions where it is absolutely continuous w.r.t. $\SLE_\kappa(\rho^-,\rho^+)$ starting from $x,x^-,x^+$. In the case $\rho^-,\rho^+>-2$ (which is the only one of use here), the $\SLE_\kappa(\rho^-,\rho^+)$ is defined for all times and is driven by a semimartingale (\cite{LSW_restr}, Section 4 in \cite{SS_freefield}). The absolute continuity properties of these processes can be expressed as in Section 3.2; one may note there that the Radon-Nikod\`ym derivative
$$\frac{{\mc Z}(c'_t){\mc Z}(c_0)}{{\mc Z}(c_t){\mc Z}(c'_0)}$$
is still well defined in the case where $x^-$ or $x^+$ are displaced under the evolution, which happens when $\rho^{\pm}<\frac\kappa 2-2$.

Consider a configuration $(D,x_1,\dots,x_m,y)$ as above, $x^-=x=x^+=x_1$. The jump of the field at $x$ is $\pm\frac{\rho^-+2+\rho^+}{\sqrt{2\pi\kappa}}$. The crosscuts $\delta_i$ are defined as before. Then Lemma \ref{locstab} holds, with the same proof. Lemma \ref{loccoupl} also follows.

\subsection{Global coupling}

We have constructed a ``local coupling" between a system of $\SLE$'s and a free field, in the case where the partition functions coincide, for a choice of crosscuts. We now use a limiting argument to construct a global coupling that will enjoy the same properties for any choice of crosscuts. This requires compatibility of the construction with the respective Markov properties of $\SLE$ systems and the free field. 
A ``local to global'' argument is introduced in \cite{DZ_revers} in the context of $\SLE$ reversibility. We broadly follow here the presentation in \cite{Dub_dual}, one difference being in the nature of the Markov properties; there is also here additional structure (conditional independence).

Again,  $c=(D,x_1,\dots,x_n,y)$ is a configuration ($y$ in the bulk), $(\underline a,b)$ a set of boundary conditions for the free field corresponding to $\SLE_\kappa(\underline\rho)$; $x_1,\dots,x_m$ are seeds of $\SLE$ ($\rho_i=2$ for $i=1,\dots, m$). We denote by $\mu^{\SLE}_c$, $\mu^{\FF}_c$ the respective distributions of the $\SLE$ system and the free field in $c$.

The goal is to construct a coupling of $\mu^{\SLE}$ and $\mu^{\FF}$ with natural compatibility with Markov properties, and such that the different $\SLE$ strands are independent conditionally on the field. It is unclear whether it is possible to do this ``in one go". So we shall consider first the coupling of two objects: one $\SLE$ strand (case $m=1$) and a free field.

Let $\eta>0$ be a small parameter, say much smaller than the diameter of $D$ and the distance between marked points. We define a sequence of stopping times for $\gamma$ by $\tau_0=0$, 
$$\tau_{n+1}=\inf\{t\geq\tau_{n}:\dist(\gamma_t,\gamma_{[0,\tau_n]})\geq\eta\}$$
It is easy to see that if $D$ is bounded, there is a fixed $N=N(D)$ such that a.s. $\tau_n=\infty$ for $n> N$.

Let $\partial$ be the union of $\{y\}$ and the smallest connected boundary arc containing all marked points except $x=x_1$. Let $p_0$ be the random integer:
$$p_0=\inf\{p: \dist(K_p,\partial)\leq 3\eta\}$$
where $K_p$ denotes the hull of the $\SLE$ stopped at $\tau_p$. Let $\delta_p$ be the boundary component of $(K_p)^{2\eta}$ that disconnects $K_p$ from other marked points. (One can consider variants that ensure that the crosscut $\delta_p$ is smooth. The important point is that is no closer than, say, $\frac 32\eta$ of $K_p$ and no farther than, say, $2\eta$, and is ${\mc F}_{\tau_p}$ measurable). We define:
$$\ell_{p}(\gamma,\phi)=\frac {d(T_p)_*\mu^{\FF}_{c_{p+1}}}
{d(T_p)_*\mu^{\FF}_{c_{p}}}(T_p\phi)$$
where $c_{p}$ is the configuration sampled at time $\tau_p$ and $T_p=T_{\delta_p}$ (trace on $\delta_p$). This quantity depends on the $\SLE$ strand up to time $\tau_{p+1}$.

Consider the measure $L.\mu^{\SLE}_c\otimes\mu^{\FF}_c$, where the density $L$ is given by:
$$L=L(\gamma,\phi)=\prod_{p\leq p_0}\ell_{p}(\gamma,\phi)$$
First we have to check that this is a coupling of $\mu^{\SLE}_c$, $\mu^{\FF}_c$. For fixed $\phi$, $p\mapsto \prod_{q=0}^{p-1}\ell_q(\gamma,\phi)$ is a discrete time martingale (it is bounded when stopped at $p_0$). This boils down to:
$$\E\left(\frac {d(T_p)_*\mu^{\FF}_{c_{p+1}}}
{d(T_p)_*\mu^{\FF}_{c_p}}(T_p\phi)|\gamma_{[0,\tau_{p}]}\right)=1$$
which follows from the local case (Lemma \ref{locstab}) in the configuration $c_p$ (note that $\delta_p$ is determined by $\gamma_{[0,\tau_p]}$). 
For the other marginal, we need another expression of the density, thinking now of $\gamma$ as fixed. This follows from $\ell_q(\gamma,\phi)=\prod_{q\leq r\leq p_0}\ell_{q,r}(\gamma,\phi)$
where 
$$\ell_{q,r}(\gamma,\phi)=\frac{d(T_r)_*\mu^{\FF}_{c_{q+1}}(.|T_{r+1}\phi)}{d(T_r)_*\mu^{\FF}_{c_{q}}(.|T_{r+1}\phi)}(T_r\phi)$$
for $q\leq r<p_0$ and
$$\ell_{q,p_0}(\gamma,\phi)=\frac{d(T_{p_0})_*\mu^{\FF}_{c_{q+1}}}{d(T_{p_0})_*\mu^{\FF}_{c_{q}}}(T_{p_0}\phi)$$
Indeed, 
$$\prod_{q\leq r\leq p_0}\ell_{q,r}(\gamma,\phi)=\frac{d(T_q,\dots,T_{p_0})_*\mu^{\FF}_{c_{q+1}}}{d(T_q,\dots,T_{p_0})_*\mu^{\FF}_{c_{q}}}(T_q\phi,\dots,T_{p_0}\phi)=\frac{d(T_q)_*\mu^{\FF}_{c_{q+1}}}{d(T_q)_*\mu^{\FF}_{c_{q}}}(T_q\phi)=\ell_q(\gamma,\phi)$$
This density factors through $T_q\phi$ because of the Markov property of the free field (the conditional distributions of $T_{q+i}\phi$ given $T_q\phi$ do not depend on what is on the other side of the crosscut $\delta_q$). Thus:
$$\prod_{q=0}^{p_0}\ell_q(\gamma,\phi)=\prod_{q\leq p_0,q\leq r\leq p_0}\ell_{q,r}=\prod_{0\leq r\leq p_0}\frac{d(T_r)_*\mu^{\FF}_{c_{r+1}}(.|T_{r+1}\phi)}{d(T_r)_*\mu^{\FF}_{c_{0}}(.|T_{r+1}\phi)}(T_r\phi)$$ 
so that integrating $T_0\phi$, then $T_1\phi$, \dots (with the convention that the conditioning by $T_{p_0+1}\phi$ is empty), one gets:
$$\int d\mu_{c_0}^{\FF}(\phi)\prod_{q=0}^{N-1}\ell_q(\gamma,\phi)=1$$ 
for fixed $\gamma$.  
This shows that we have indeed a coupling of $\mu^{\SLE}_c$, $\mu^{\FF}_c$. It is easy to see that in this coupling, $\gamma^{\tau_p}$ is independent of $\phi$ conditionally on $\phi$ inside $(K_p)^{2\eta}$ (more precisely, $T_0\phi,\dots,T_p\phi$). This holds for fixed $p$ or a stopping time for the discrete time filtration $(\sigma(\gamma^{\tau_p}))_{p\geq 0}$.

In the general case, one can proceed in different ways. For simplicity, we can first use a common (discrete) time scale for the $\SLE$ system: at each step, each strand moves at distance $\eta$, synchronously; this yields a sequence of configurations $c_p=c_{\tau^1_p,\dots,\tau^m_p}$. Consider $\delta_p=\delta_p^1\sqcup\dots\sqcup\delta^m_p$. Then we can define similarly to the $m=1$ case:
$$\ell_p(\gamma,\phi)=\frac {d(T_{\delta_p})_*\mu^{\FF}_{c_{p+1}}}
{d(T_{\delta^i_p})_*\mu^{\FF}_{c_p}}(T_{\delta_p}\phi)$$
where $\gamma$ is now an $m$-uplet of paths. As before, the local study ensures that $p\mapsto\prod_{q=0}^{p-1}\ell_q$ is a martingale for fixed $\phi$ (say stopped at $p_0$, the first $p$ such that $(p,\dots,p)\notin G$). Moreover, due to the nested structure of the $\delta_p$'s and the Markov property of the field, we can write:
$$\prod_{q=0}^{p_0-1}\ell_q=\prod_{q=0}^{p_0-1}\frac{d(T_q)_*\mu^{\FF}_{c_{q+1}}(.|T_{q+1}\phi)}{d(T_q)_*\mu^{\FF}_{c_{0}}(.|T_{q+1}\phi)}(T_q\phi)$$
so that for fixed $\gamma$, $\int Ld\mu^{\FF}_{c_0}(\phi)=1$, integrating first $T_0\phi$, then $T_1\phi$ etc, where $L=\prod_{q=0}^{p_0-1}\ell_p$. Thus $L.\mu^{\FF}_{c_0}\otimes\mu^{\SLE}_{c_0}$ defines a coupling of $\mu^{\FF}_{c_0}$, $\mu^{\SLE}_{c_0}$.

Let us go back to the $m=1$ case. We saw that the measure $L.\mu^{\SLE}_c\otimes\mu^{\FF}_c$, where the density $L$ is given by:
$$L=L(\gamma,\phi)=\prod_{p\leq p_0}\ell_{p}(\gamma,\phi)$$
is a coupling of $\mu^{\SLE}_c$, $\mu^{\FF}_c$. Let $n$ be a stopping time in the discrete filtration $(\sigma(\gamma^{\tau_p}))_{p\geq 0}$. Then the measure induced on $(\gamma^{\tau_n},\phi_{|(K_n^{2\eta})^c})$ can be described as follows. The first marginal is just the $\SLE$  strand stopped at $\tau_n$. Conditionally on $\gamma^{\tau_n}$, the distribution of $T_{n-1}\phi$ is that induced by the free field in $c_n$; consequently, the conditional distribution of $\phi_{|(K_n^{2\eta})^c}$ is that of the free field in $c_n$ restricted to $(K_n^{2\eta})^c$. Indeed, the measure induced on $(\gamma^{\tau_n},\phi)$ is simply:
$$\prod_{q=0}^{n(\gamma)-1}\ell_q(\gamma,\phi)d\mu^{\SLE}_c(\gamma^{\tau_n})d\mu^{\FF}_c(\phi)$$
(recall that $p\mapsto\prod_{q=0}^{p-1}\ell_q(\gamma,\phi)$ is a discrete time martingale for fixed $\phi$). Then, reasoning as above (with $n-1$ replacing $p_0$), one gets the expression:
$$\prod_{q=0}^{n-1}\ell_q(\gamma,\phi)=\frac{d(T_{n-1})_*\mu^{\FF}_{c_{n}}}{d(T_{n-1})_*\mu^{\FF}_{c_{0}}}(T_{n-1}\phi)
\prod_{0\leq r<n-1}\frac{d(T_r)_*\mu^{\FF}_{c_{r+1}}(.|T_{r+1}\phi)}{d(T_r)_*\mu^{\FF}_{c_{0}}(.|T_{r+1}\phi)}(T_r\phi)$$
and then one integrates out successively $T_0\phi,\dots,T_{n-2}\phi$.

The construction depends on a small parameter $\eta$, which we now take to 0. The sequence of paired measures $L_\eta.\mu^{\SLE}_c\otimes\mu^{\FF}_c$ has fixed marginals, hence is tight. Thus there exists a subsequential limit $\Theta$, which is again a coupling of $\mu^{\SLE}_c$, $\mu^{\FF}_c$. 

To phrase properties of the coupling, we need to introduce filtrations. First, $({\mc F}^{\SLE}_t)_{t\geq 0}$ is the filtration generated by the $\SLE$ strand. The time scale is arbitrary and the discussion here is invariant under bicontinuous progressive time change (under which the class of stopping times is invariant). A possible time scale is the half-plane capacity of $\psi(\gamma^.)$, where $\psi$ is some conformal equivalence $D\rightarrow\H$. 

Recall that for the free field, we defined ${\mc F}^{\FF}_U=\sigma(\langle \phi,f\rangle_{L^2})_{f\in C^\infty_0(U)}$ for $U$ open subset of $D$, and ${\mc F}^{\FF}_K=\bigcap_{U\supset K}{\mc F}^{\FF}_U$ for $K$ closed. The set of open subsets of $D$ is partially ordered for inclusion, so we can think of $({\mc F}^{\FF}_U)_U$ as a filtration with partially ordered index set (plainly, $U\subset V$ implies ${\mc F}^{\FF}_U\subset{\mc F}^{\FF}_V$). We can phrase now:

\begin{Thm}\label{Tcoupl}
Let $\mu^{\SLE}_c$, $\mu^{\FF}_c$ be the distributions of an $\SLE$ and a free field in a configuration $c$ with common partition functions. Then there exists a coupling $\Theta$ of $\mu^{\SLE}_c$, $\mu^{\FF}_c$ which is maximal in the following sense:
\begin{enumerate}
\item For all ${\mc F}^{\SLE}$-stopping time $\tau$, conditionally on ${\mc F}^{\SLE}_\tau$ the field restricted to $D\setminus\gamma^\tau$ has distribution $\mu^{\FF}_{c_\tau}$.
\item For all open set $U$ having the seed $x$ of the $\SLE$ on its boundary (a continuous arc), the field restricted to $D\setminus U$ is independent of the $\SLE$ stopped upon exiting $U$ conditionally on ${\mc F}^{\FF}_{\partial U}$. Equivalently, the $\SLE$ stopped upon exiting $U$ is independent of the field conditionally on the field restricted to $U$. 
\end{enumerate}
\end{Thm}
\begin{proof}
The limiting arguments here are similar to those in Theorem 6 in \cite{Dub_dual}. Let $\eta_k\searrow 0$ be a sequence along which $L_\eta.\mu^{\SLE}_c\otimes\mu^{\FF}_c$ has a limit $\Theta$.
For the first statement, one can consider a probability space with sample $(\gamma,\phi,\phi_1,\dots,\phi_k,\dots)$ such that $\phi_k\rightarrow\phi$ a.s. (e.g. in the Fr\'echet topology of $C_0^\infty(D)'$) and the marginal $(\gamma,\phi_k)$ has distribution $L_{\eta_k}\mu^{\SLE}_c\otimes\mu^{\FF}_c$. Let us assume first that $\gamma^\tau$ is at uniformly bounded below distance of $\partial$. Consider
$$n(k)=\inf\{n:\tau^k_n\geq\tau\}$$
where the sequence of stopping times $(\tau^k_n)_n$ is from the definition of $L_k=L_{\eta_k}$. Then $\tau^k=\tau^k_{n(k)}$ is a stopping time and $\tau^k\searrow\tau$ a.s. Let $\eps>0$ be fixed. For $k$ large enough (viz. $3\eta_k\leq \eps$), conditionally on ${\mc F}^{\FF}_{\tau^k}$, the field $\phi_k$ restricted to the connected component of $D\setminus(\gamma_\tau)^\eps$ having $\partial$ on its boundary has the distribution of the free field in $c_{\tau^k}$ restricted to that set. One concludes by taking $k\rightarrow \infty$ and then $\eps\searrow 0$.

One obtains the second statement by applying the first statement to $\tau_U=\inf\{t:\gamma_t\notin U\}$ in conjunction with the Markov property of the free field.
\end{proof}

In the situation with several strands, one can rely on the local computation in Lemma \ref{loccoupl} to reduce the problem to one strand.

\begin{Thm}\label{TCoupl2}
Let $\mu^{\SLE}_c$, $\mu^{\FF}_c$ be the distributions of a system of $m$ $\SLE$'s and a free field in a configuration $c$ with common partition functions. Then there exists a coupling $\Theta$ of $\mu^{\SLE}_c$, $\mu^{\FF}_c$ such that:
\begin{enumerate}
\item the marginals $(\gamma_1,\phi)$, \dots, $(\gamma_m,\phi)$ are maximal couplings; 
\item the $\SLE$ strands $\gamma_1,\dots,\gamma_m$ are independent conditionally on the field $\phi$.
\end{enumerate}
\end{Thm}
One can obtain more general stopping statements (involving e.g. sequences of stopping times for the different $\SLE$ strands), which are a bit heavy to formulate and of no direct use here. 

\begin{proof}
For $\eta>0$, consider a coupling
$$L^1_\eta(\gamma_1,\phi)\dots L^m_\eta(\gamma_m,\phi)d\mu^{\SLE}_c(\gamma_1)\dots d\mu^{\SLE}_c(\gamma_m)d\mu^{\FF}_c(\phi)$$
where $L^i_\eta(\gamma_i,\phi)$ is the density we considered above. The marginal distributions are $d\mu^{\SLE}_c(\gamma_i)$ ($\SLE$ system restricted to the $i$-th strand, $i=1,\dots,m$) and $d\mu^{\FF}_c(\phi)$. Moreover, the $\gamma_i$'s are independent conditionally on $\phi$, due to the split form of the density for fixed $\phi$. As $\eta\searrow 0$, the family of measures is tight.

Consider a sequence $\eta_k\searrow 0$ along which these couplings converge to a measure $\Theta$ on $(\gamma_1,\dots,\gamma_m,\phi)$. In particular, the distributions of the marginals $(\gamma_1,\phi)$, \dots, $(\gamma_m,\phi)$ converge. Then the limiting distributions of these paired marginals are maximal couplings, as in the proof of the previous theorem. It is also clear that the conditional independence of the $\SLE$ strands given the field is preserved in the limit. 

What remains to check is that under $\Theta$, $(\gamma_1,\dots,\gamma_m)$ is (jointly) distributed according to $\mu^{\SLE}_c$. Consider disjoint crosscuts $\delta_i$, $i=1\dots m$ separating $x_i$ (the seed of the $i$-th $\SLE$) from all other marked points; more precisely, $D\setminus\delta_i=L_i\sqcup R_i$, with $x_i\in\partial L_i$. The $i$-th $\SLE$ is stopped at time $\tau_i$, when it comes within distance $\eps>0$ of $\delta_i$. For $i\neq j$, $\gamma_i$ is independent from $\gamma_j$ conditionally on $\phi$; besides, $\gamma_i^{\tau_i}$ is depends on $\phi$ only through its restriction to $L_i$. Moreover, the restrictions of the field in the $L_i$'s are independent conditionally on the trace of the field on $\delta=\sqcup_i\delta_i$. Hence the $\gamma_i^{\tau_i}$'s are independent given the trace $T_\delta\phi$. Since $(\gamma_i,\phi)$ is a maximal coupling, the joint distribution of $(\gamma_i^{\tau_i},T_\delta\phi)$ is that of the $\SLE$ stopped at $\tau_i$, and conditionally on $\gamma_i^{\tau_i}$, $\phi$ is distributed as the field in $c_{\tau_i}$. This shows that the distribution under $\Theta$ of $(\gamma_1^{\tau_1},\dots,\gamma_m^{\tau_m},T_\delta\phi)$ is the same as the one in the local coupling of Lemma \ref{loccoupl}. In particular, the joint distribution of $(\gamma_1^{\tau_1},\dots,\gamma_m^{\tau_m})$ is that of the $\SLE$ system $\mu_c^{\SLE}$ with the $i$-th strand stopped at $\tau_i$. Since this is valid for all crosscuts $(\delta_i)$ and all $\eps>0$, the joint distribution of $(\gamma_1,\dots,\gamma_m)$ under $\Theta$ is indeed $\mu_c^{\SLE}$.
\end{proof}

\section{Stochastic ``differential" equations driven by the free field}

In order to build some intuition on the nature of the relationship between $\SLE$ and the free field studied here, it appears rather convenient to draw an analogy with the standard theory of stochastic differential equations (SDEs) driven by (real) Brownian motion (see eg \cite{RY}). However the free field/$\SLE$ situation does not involve a stochastic calculus w.r.t. free field.

\subsection{Definitions}

Let us briefly recall the set-up for stochastic differential equations. Let $(X,B)$ be a pair of adapted processes in a probability space $(\Omega,{\mc F},\P)$, $\sigma_t((X^t_.))$, $b_t((X^t_.))$ progressively measurable functions of the process $X$. The SDE reads:
$$X_t=\int_0^t \sigma_s((X^s_.))dB_s+\int_0^t b_s((X^s_.))ds.$$
A pair $(X,B)$ is a solution of the SDE if $B$ is an ${\mc F}$- Brownian motion and the relation is satisfied (given $B$ and $X$, the RHS is defined as a stochastic integral). It is a strong solution if moreover ${\mc F}$ is generated by $B$. There is uniqueness in law if in all solutions $(X,B)$, the marginal $X$ has the same distribution, and pathwise uniqueness if for any pair of solutions $(X,B)$, $(X',B)$ defined on a common filtered space (with common driving BM's), the processes $X,X'$ are undistinguishable (a.s. equal).

For simplicity, we restrict ourselves to the chordal case: a configuration $c=(D,x,y)$ consists in a simply connected domain $D$ with two marked points $x,y$ on the boundary.

Consider a filtration $({\mc F}_U)_U$ indexed by open neighbourhoods of $x$ in $D$. An ${\mc F}$-free field is a free field such that ${\mc F}^{\FF}_U\subset{\mc F}_U$ and $\phi$ restricted to $D\setminus\overline U$ is independent of ${\mc F}_U$ conditionally on ${\mc F}^{\FF}_{\partial U}$. A stochastic Loewner chain $K_.$ starting from $x$ is ${\mc F}$-adapted if $K$ stopped at first exit of $U$ is ${\mc F}_U$-measurable. We only consider Loewner chains with continuous driving functions. Assume given an assignment: 
$$(K_s)_{0\leq s\leq t}\longmapsto h((K_s)_{0\leq s\leq t})$$
where $h$ is a harmonic function in $D_t=D\setminus K_t$.
(One may also consider the situation where $h$ is defined in $D\setminus\gamma_{[0,t]}$, for a Loewner chain generated by a trace $\gamma$).

We are interested in comparing the boundary values of $h=h((K_s)_{0\leq s\leq t})$ and $\phi$ in $D_t$; the issue is that neither need be defined pointwise on $\partial D_t$. One may proceed as follows: consider a sequence $\delta_n$ of closed, smooth curves converging to $\partial D_t$ (eg equipotentials seen from a bulk point). The $\delta_n$'s depend on the chain but not on the field. Then one requires that the harmonic extension of the trace of the field on $\delta_n$ inside $\delta_n$ (this is a.s. well defined) converges to $h$ uniformly on compact sets of $D_t$. Plainly, this can be checked pathwise. A more compact (if less explicit) formulation in terms of conditional expectation of the field is possible:

\begin{Lem}\label{Lem:eqpwise}
Let $\phi$ be an ${\mc F}$-free field, $K_.$ an ${\mc F}$-adapted Loewner chain. Let $(\delta_n)$ be a $\sigma(K^t_.)$ measurable sequence of nested closed curves approximating $\partial K_t$. Then a.s. conditionally on $K^t_.$, the harmonic extension of $T_{\delta_n}\phi$ inside $\delta_n$ converges uniformly on compact sets to $\E(\phi_{|D_t}|{\mc F}^{\FF}_{\partial D_t})$, the conditional expectation of the field restricted to $D_t$ given ${\mc F}^{\FF}_{\partial D_t}$.
\end{Lem} 
\begin{proof}
The curve $\delta_n$ splits $D_t$ into $L_n$ (that has $\partial K^t_.$ on its boundary) and $R_n$. A compact set $C$ of $D_t$ is contained in $R_n$ for $n$ large enough. From the properties of the field trace, $h_n=\E(\phi_{|C}|{\mc F}^{\FF}_{L_n})=\E(P_{R_n}T_{\delta_n}\phi)_{|C}$. If $m\leq n$, 
$$\E(||h_n-h_m||^2_{L^2(C)})=\E(\E(||h_m-h_n||^2_{L^2(C)}|{\mc F}^{\FF}_{L_n}))$$
and the conditional expectation is constant, since $\delta_m$ is contained in $R_n$. So we can compute it for the free field in $R_n$ with Dirichlet boundary conditions (so that $h_n=0$). Thus:
$$\E(||h_m-h_n||^2_{L^2(C)}|{\mc F}^{\FF}_{L_n})=\E_{R_n}\int_{(\delta_m)^2}\int_CP_{R_m}(x,z)(T_{\delta_m}\phi)(x)P_{R_m}(y,z)(T_{\delta_m}\phi)(y)dl(x)dl(y)dA(z)$$
This can be exactly evaluated by general Gaussian arguments. Under $\mu^{\FF}_{R_n}$, $w=T_{\delta_m}\phi$ has covariance $(G_{R_n})_{|\delta_m}$. For a symmetric kernel $B(x,y)=\sum_{i,j} f_i(x)f_j(y)$ on $(\delta_m)^2$, one gets:
\begin{align*}
\E_{R_n}(\int_{(\delta_m)^2}w(x)B(x,y)w(y)dl(x)dl(y))&=\sum_{i,j}\E_{R_n}(\int_{\delta_m} f_iwdl\int_{\delta_m} f_jwdl)\\
&=\sum_{i,j}\int_{(\delta_m)^2}f_i(x)G_{R_n}(x,y)f_j(y)dl(x)dl(y)=\Tr_{L^2(\delta_m)}(BG_{R_n})
\end{align*}
This is valid for finite rank kernels, and by approximation applies to the trace class kernel:
$$B(x,y)=\int_C P_{R_m}(x,z)P_{R_m}(y,z)dA(z)$$
It follows that:
$$\E(||h_m-h_n||^2_{L^2(C)}|{\mc F}^{\FF}_{L_n})=\int_{(\delta_m)^2}\int_CP_{R_m}(x,z)G_{R_n}(x,y)P_{R_m}(y,z)dl(x)dl(y)dA(z)$$
In terms of path decompositions, this corresponds to a Brownian loop in $R_n$ starting and ending at $z$ and decomposed w.r.t. its first and last visit to $\delta_m$. Given that the transition kernels in $R_n$ have uniform exponential decay and that as $n\rightarrow\infty$, the transition kernel in $R_m$ converges uniformly to that of $D_t$ uniformly on $C\times C\times [0,t]$, it is easy to see that $\E(||h_m-h_n||^2_{L^2(C)})$ converges to 0 as $m\rightarrow\infty$, uniformly in $n\geq m$. One can refine this (using eg the Harnack inequality to control derivatives of the Poisson kernel) to get that for any $k\geq 0$, $\E(||h_m-h_n||^2_{\sH^k(C)})$ converges to 0 as $m\rightarrow\infty$, uniformly in $n\geq m$. It follows that $(h_m)$ converges a.s. in any $\sH^k(C)$, and consequently (Sobolev imbedding) converges uniformly in $C$ (ie in $(C_0(C),||.||_\infty)$).

From the free field Markov property, we have $h=\E(\phi_{|C}|{\mc F}^{\FF}_{\partial D_t})=\E(\phi_{|C}|{\mc F}^{\FF}_{K_t\cup\partial D})$. By definition of ${\mc F}^{\FF}$ on closed sets, we have ${\mc F}^{\FF}_{K_t\cup\partial D}=\cap_{n>0}{\mc F}^{\FF}_{L_n}$.
For any $m>0$, we have:
$$h=\E(\phi_{|C}|\cap_n{\mc F}^{\FF}_{L_n})=\E(\E(\phi_{|C}|{\mc F}^{\FF}_{L_m})|\cap_n{\mc F}^{\FF}_{L_n})=\E(h_m|\cap_n{\mc F}^{\FF}_{L_n})$$
Since $\lim_m h_m$ exists a.s. and is consequently $(\cap_n{\mc F}^{\FF}_{L_n})$-measurable, we get that $h=\lim_m h_m$, which concludes the proof.
\end{proof}

Consider the following problem, given the data of $h=h((K_s)_{0\leq s\leq t})$: find a  probability space with filtration $({\mc F}_U)_U$ on which are defined a field $\phi$ and a stochastic Loewner chain $(K_t)_{t\geq 0}$ such that:
\begin{enumerate}
\item $\phi$ is an ${\mc F}$-free field,
\item $(K_.)$ is ${\mc F}$-adapted,
\item 
For all $t\geq 0$, $\E(\phi_{|D_t}|{\mc F}^{\FF}_{\partial D_t})=h((K^t_.))$.
\end{enumerate}

This imposes some compatibility conditions on $h$ under stopping of the Loewner chain: if $s\leq t$, $h_s=h((K^s_.))$ and $h_t=h((K^t_.))$ agree on $\partial D_s$ in the sense that if $(\delta_n)$ is a sequence of closed curves approaching $\partial D_s$, the harmonic extension of the restriction of $h_t$ to $\delta_n$ converges to $h_s$ locally uniformly in $D_s$. 
Under continuity assumptions (as in Lemma \ref{bcont}), it is enough to check the condition $\E(\phi_{|D_t}|{\mc F}^{\FF}_{\partial D_t})=h((K^t_.))$ for a countable dense set of times $\{t_i\}$. 

By analogy with the SDE framework, one can phrase:

\begin{Def}
The stochastic equation 
$$\E(\phi_{|D_t}|{\mc F}^{\FF}_{\partial D_t})=h((K^t_.))\quad \forall t\geq 0$$
has a solution if there exists a filtered probability space $(\Omega,{\mc F},({\mc F}_U),\P)$ on which are defined an ${\mc F}_.$-free field $\phi$ and an ${\mc F}_.$-adapted stochastic Loewner chain $K_.$ satisfying the equation.

The solution is strong if moreover ${\mc F}_.={\mc F}^{\FF}_.$.

There is uniqueness in law if for any two solutions $(\phi,K)$, $(\phi',K')$, the marginal distributions of the Loewner chain are identical.

There is pathwise uniqueness if for any filtered space on which are defined a field $\phi$ and two chains $K$, $\tilde K$ such that $(\phi,K)$ and $(\phi,\tilde K)$ are solutions, the Loewner chains are a.s. equal.
\end{Def}

\subsection{Existence and uniqueness in law}

We have considered different types of boundary conditions, in particular chordal $(a,b)$ boundary conditions in a configuration $c=(D,x,y)$. This defines an assignment $h_{a,b}=h_{a,b}(K^t_.)$, provided that the domain $D$ is regular enough ($C^1$) around $y$ and stopped chains $K^t_.$ stay away from $y$.

\begin{Thm}\label{Thm:SDEweak}
Let $\kappa>0$, $a=\pm\sqrt{\frac 2{\pi\kappa}}$, $b=a(1-\frac\kappa 4)$. Then the stochastic equation in $(\phi,K_.)$:
$$\E(\phi_{|D_t}|{\mc F}^{\FF}_{\partial D_t})=h_{a,b}((K^t_.))\quad \forall t\geq 0$$
has a solution. It is unique in law and distributed as chordal $\SLE_\kappa$ in $(D,x,y)$.
\end{Thm}
\begin{proof}
Existence. Take a maximal coupling of a free field with $(a,b)$ boundary conditions and a chordal $\SLE_\kappa$ in $c$, which exists by Theorem \ref{Tcoupl}. Define ${\mc F}_U={\mc F}^{\SLE}_{\tau_U}\vee{\mc F}^{\FF}_U$, where $\tau_U$ is the time of first exit of $U$ by the $\SLE$. By definition, the $\SLE$ is ${\mc F}$-adapted; and $\phi$ is an ${\mc F}$-free field by Theorem \ref{Tcoupl}. Besides, for a time $t$, $\phi_{|D_t}$ is distributed as an $(a,b)$ free field in $D_t$ conditionally on ${\mc F}_{\partial D_t}$. It follows that $\E(\phi_{|D_t}|{\mc F}^{\FF}_{\partial D_t})=h_{a,b}((K^t_.))$ is a.s. satisfied at $t$; consequently it is a.s. satisfied for $t$ in a dense countable set of times $\{t_i\}$. It is then easy to see that the equation is satisfied for all times, a.s.

Uniqueness. We reason as in Lemma \ref{locstab}, in reverse (a standard argument, see eg \cite{LSW_LERW}). Consider a solution $(\phi,K_.)$, with filtration ${\mc F}$; denote ${\mc G}_t={\mc F}_{K_t}$. By the Markov property of the field, the distribution of $\phi_{|D_t}$ conditionally on ${\mc G}_t$ is that of a free field in $D_t$ with mean $h_t=h_{a,b}(K^t_.)$. Consequently, if $f\in C^\infty_0(D)$ is a test function 
and $\tau=\inf\{t\geq 0: \dist(K^t_.,supp(f))\leq\eps\}$, we have:
$$M_{t\wedge\tau}\stackrel{def}{=}\E(\langle \phi,f\rangle_{L^2}|{\mc G}_{t\wedge\tau})=\langle h_{t\wedge\tau},f\rangle_{L^2}$$ 
and $M$ is by construction a bounded ${\mc G}$-martingale, and is continuous (the Loewner chain is also assumed to be generated by a continuous process). For simplicity, map $(D,x,y)$ to $(\H,0,\infty)$. Then $h_.$ is a martingale taking values in the space of harmonic functions in a neighbourhood of infinity (the boundary condition is local, hence it is fixed in a neighbourhood of infinity). The harmonic conjugation (with condition $h^*(\infty)=0$) is a linear operation. Consequently, $(h+ih^*)_t(z)$ is a complex valued martingale for $z$ in a neighbourhood of $y=\infty$. With usual notations, this means that if $Z_t=g_t(z)-W_t$, 
$$m_t(z)=-a\log(Z_t)-b\log g'_t(z)$$
is a martingale (bounded if stopped upon exiting $D(0,M)$, $|z|>M$). Since $\log g'_t(z)=-\int_0^t\frac{2ds}{Z_s^2}$, we have:
$$Z_t=\exp(-a^{-1}(m_t(z)+b\int_0^t\frac{2ds}{Z_s^2}))$$
so that $Z$ is a semimartingale; then $W_t=-Z_t+z+\int_0^t\frac{2ds}{Z_s}$ is also a semimartingale. Thus one can write $dW_t=\sigma_t dB_t+db_t$; plugging this back in 
$$-dm_t(z)=\frac{a}{Z_t}\left(\frac {2}{Z_t}dt-dW_t\right)-\frac a{2Z_t^2}d\langle W\rangle_t+\frac{2b}{Z_t^2}dt$$
evaluated at two distinct $z$ points, one gets $\sigma_t\equiv \sqrt\kappa$ and $db_t=0$. Thus the Loewner chain is distributed as chordal $\SLE_\kappa$.
\end{proof}

We restricted to the chordal case for simplicity; however it is clear that the result applies whenever an identity of partition functions as in Theorem \ref{PFident} holds. Following the discussion at the end of Section 6.1, it also applies when the $\SLE$ strand is absolutely continuous w.r.t. an $\SLE_\kappa(\rho^-,\rho^+)$, $\rho^\pm>-2$, near its start at $x^-=x=x^+$.

\subsection{Pathwise uniqueness}

In this subsection, we are considering the question of pathwise uniqueness in the chordal case for general $\kappa>0$. Pathwise uniqueness combined with the already established existence of weak solutions implies existence of strong solutions (in which the $\SLE$ path is a function of the field).

The general strategy consists in starting from a weak solution $(\phi,K)$ to construct a triplet $(K,\phi,\hat K)$ where $\hat K$ is a dual $\SLE$ path (or collection of such paths) such that $K$, $\hat K$ are independent conditionally on the field and $\hat K$ determines $K$. This implies that $K$ is actually a strong solution. Moreover, if $K$, $\tilde K$ are two weak solutions defined on the same probability space (common field), then $K$, $\tilde K$ are equal since they are determined by the common auxiliary $\hat K$; this yields pathwise uniqueness.

The construction of the auxiliary path (or collection of paths) depends on $\kappa$; we will consider separately the cases $\kappa=4$, $\kappa<4$, $\kappa\geq 8$, $4<\kappa<8$.

{\bf Case $\kappa=4$}

In the cases $\kappa=4,8$, the corresponding free field boundary conditions have symmetries compatible with reversibility. We now exploit this fact, in conjunction with the simplicity of the trace, for $\kappa=4$.

We have already proved the existence of a solution. It is enough to prove that if $(\phi,K)$, $(\phi,\tilde K)$ are two solutions defined on the same filtered space, then $K=\tilde K$. It implies in particular that all solutions are strong (as in the case of SDEs).

To be able to use densities, we prove a different version. Namely, consider a solution $(\phi,K)$ of the problem in $(D,x,y)$ ($(a,0)$ boundary conditions) and $(\phi,\hat K)$ a solution in $(D,y,x)$ ($(-a,0)$ boundary conditions), coupled so that the fields agree (they have the same boundary conditions) and the chains are independent conditionally on the field. Then we claim that $K$, $\hat K$ have the same range. As these are simple paths, this determines the chain completely. Applying twice this result (take $(\phi,\hat K)$ a solution of the problem in $(D,y,x)$, independent of $K,\tilde K$ conditionally on $\phi$; then $K$ and $\tilde K$ are the reverse of $\hat K$), one gets pathwise uniqueness.

Hence we consider a triplet $(K,\phi,\hat K)$ such that $(K,\phi)$, $(\phi,\hat K)$ are solutions in $(D,x,y)$, $(D,y,x)$ respectively, and $K,\hat K$ are independent conditionally on $\phi$. We reason now as in Theorem \ref{TCoupl2}. Consider a crosscut $\delta$ splitting $D$ in $L,R$ ($x$ on the boundary of $L$, $y$ on the boundary of $R$). The chains $K,\hat K$ are stopped at $\tau,\hat\tau$ when they come within distance $\eps>0$ of the crosscut $\delta$.
Then $K^\tau$ is independent of ${\mc F}^{\FF}_R$ conditionally on ${\mc F}^{\FF}_\delta$; $\hat K^{\hat\tau}$ is independent of ${\mc F}^{\FF}_L$ conditionally on ${\mc F}^{\FF}_\delta$; $K^\tau$ is independent of $\hat K^{\hat\tau}$ conditionally on the field; ${\mc F}^{\FF}_L$, ${\mc F}^{\FF}_R$ are independent conditionally on ${\mc F}^{\FF}_\delta$. It follows that $K^\tau$, $\hat K^{\hat\tau}$ are independent conditionally on ${\mc F}^{\FF}_\delta$. 

Besides, the marginal distributions of $(K^\tau,T_\delta\phi)$, $(\hat K^{\hat\tau},T_\delta\phi)$ are fixed:
$$d\mu^{\SLE}_c(K^\tau)d(T_\delta)_*\mu_{c_\tau}^{\FF}(T_\delta\phi)$$
and symmetrically for $(\hat K^{\hat\tau},T_\delta\phi)$. The joint distribution of $(K,T_\delta\phi,\hat K^{\hat\tau})$ is thus:
$$\frac{d(T_\delta)_*\mu_{c_\tau}^{\FF}}{d(T_\delta)_*\mu_{c}^{\FF}}(T_\delta\phi)\cdot\frac{d(T_\delta)_*\mu_{{\hat c}_{\hat\tau}}^{\FF}}{d(T_\delta)_*\mu_{c}^{\FF}}(T_\delta\phi)\cdot d\mu^{\SLE}_c(K^\tau)d\mu^{\SLE}_{\hat c}(\hat K^{\hat\tau})d(T_\delta)_*\mu_{c}^{\FF}(T_\delta\phi)$$
To obtain the joint distribution of $(K^\tau,\hat K^{\hat\tau})$, one integrates out $T_\delta\phi$; as in Lemma \ref{loccoupl}, this yields the joint distribution:
$$d\mu_c^{\SLE}(K^\tau)d\mu^{\SLE}_{\hat c_{\tau}}(\hat K^{\hat\tau})
=d\mu^{\SLE}_{\hat c}(\hat K^{\hat\tau})d\mu_{c_{\hat\tau}}^{\SLE}(K^\tau)$$
ie the same distribution as when $\hat K$ is the reverse of $K$. Since this holds for all crosscuts $\delta$ and all $\eps>0$, it is easy to see 
that in this coupling $\hat K$ is the reverse of $K$.

{\bf Case $\kappa\in (0,4)$}

When $\kappa\notin\{4,8\}$, the coupling of the free field and chordal $\SLE$ is not compatible with $\SLE$ reversibility (at least, not in an obvious way). But it is still compatible with some duality identities (eg \cite{Dub_dual}, in particular Proposition 10), which will be enough for our purposes.

So consider a solution $(K,\phi)$ of the stochastic equation relative to $(a,b)$ boundary conditions, $a=\pm\sqrt\frac{2}{\pi\kappa}$, $b=a(1-\frac\kappa 4)$ in a domain $(D,x,y)$. It is clearer to begin with a regular version with additional marked points: in $(D,z_1,x,z_2,y)$ (marked points in this order on the boundary), consider $(\underline a,b)$ boundary conditions with $\underline a=\frac{\kappa-4}{\sqrt{2\pi\kappa}},\frac{2}{\sqrt{2\pi\kappa}},\frac{\kappa-4}{2\sqrt{2\pi\kappa}}$. Via Theorem \ref{PFident}, this corresponds to an $\SLE_\kappa(\rho_1,\rho_2)$ from $x$ to $y$ in $(D,x,z_1,z_2,y)$ ($\rho_1=(\kappa-4)$ at $z_1$, $\rho_2=(\kappa-4)/2$ at $z_2$). Eventually, $z_1,z_2$ will collapse on $y$, yielding simply chordal $\SLE_\kappa$. 

The parameters are chosen so that there is a dual chain $\hat K$ which is an $\SLE_{\hat\kappa}(\hat\rho_1,\hat\rho_2)$ from $y$ to $x$, $\hat\kappa=16/\kappa$, $\hat\rho_1=\hat\kappa-4$, $\hat\rho_2=(\hat\kappa-4)/2$, with the same partition function. For instance
$$\frac{\kappa-4}{\sqrt{2\pi\kappa}}=-\frac{\hat\kappa-4}{\sqrt{2\pi\hat\kappa}}$$ 
so that the fields associated to $K$, $\hat K$ share the same boundary conditions, up to a global sign.

Given a weak solution $(K,\phi)$, one can thus construct a triplet $(K,\phi,\hat K)$ where $\hat K$ is independent of $K$ conditionally on the field and $(\hat K,\phi)$ is a solution of the dual equation. 
We study the joint distribution $(K,\hat K)$. Consider two disjoint crosscuts $\delta_1,\delta_2$ disconnecting $x$ (resp. $y$) from other marked points; the chains $K$, $\hat K$ are stopped at $\tau,\hat\tau$ when they come within distance $\eps>0$ of $\delta=\delta_1\sqcup\delta_2$. 
Arguing as in the $\kappa=4$ case, we see that $K^\tau$, $\hat K^{\hat\tau}$ are independent conditionally on $T_\delta\phi$, and consequently the distribution of the triplet $(K^\tau,T_\delta\phi, \hat K^{\hat\tau})$ is as in Lemma \ref{loccoupl}. 
Integrating out $T_\delta\phi$ shows that $(K,\hat K)$ is a maximal coupling of $K,\hat K$. In such a coupling (\cite{Dub_dual}, Proposition 10), $K$ (stopped when it hits the boundary arc $(z_2,z_1)$) is the (say, left) boundary of the range of $\hat K$ (stopped when it hits $(z_1,z_2)$ at $x$). Thus $K$ is determined by $\hat K$; one concludes as in the $\kappa=4$ case.

{\bf Case $\kappa\geq 8$}

The argument here is best understood in terms of Uniform Spanning Trees (UST). Chordal $\SLE_8$ is the scaling limit of the Peano path of a UST with Dirichlet/Neumann boundary conditions, \cite{LSW_LERW}. The auxiliary object $\hat K$ we are using is an arbitrarily fine subtree (and dual subtree) with finitely many branches, that are $\SLE_2$-type curves.

The following lemma provides path decompositions for some versions of $\SLE_\kappa$, $\kappa\geq 8$ (notice that the statements are simpler in the case $\kappa=8$).

\begin{Lem}\label{Lem:decomp8}
In a configuration $c=(D,x,z_1,y,z_2)$, consider an $\SLE_\kappa(\underline\rho)$ chain $K$, $\kappa\geq 8$, $\rho=\frac\kappa 2-4,2,\frac\kappa 2-4$ at $z_1,y,z_2$ ($x,z_1,y,z_2$ in this order on the boundary), coupled with a free field $\phi$; let $\mu_c$ be the law of that $\SLE$. Let $z$ be a point on the boundary arc $(xy)$; $D_l$ the random subdomain swallowed when the trace hits $z$, $\hat K$ the boundary arc $\partial D_l\cap D$; $D_r=D\setminus \overline{D_l}$; the endpoints of $\hat K$ are $z$ and a random point $z'$ on $(yx)$. The dual path $\hat K$ is determined by the field and the restrictions of $(K,\phi)$ to $D_l$, $D_r$ respectively are independent conditionally on $\hat K$. The marginal distributions are ($\hat\kappa=16/\kappa$):
\begin{enumerate}
\item If $z\in (x,z_1)$, $\hat K$ is an $\SLE_{\hat\kappa}(\underline\rho)$ starting from $z$ in $D$, $\rho=-\frac{\hat\kappa}2,{\hat\kappa}-2,-\frac{\hat\kappa}2,{\hat\kappa}-2,-\frac{\hat\kappa}2,\frac{\hat\kappa}2-2$ at $z^+,z_1,y,z_2,x,z^-$ (in this order), stopped when it hits $(yx)$. If $z\in (z_1,y)$, $\hat K$ is an $\SLE_{\hat\kappa}(\underline\rho)$ starting from $z$ in $D$, $\rho=\frac{\hat\kappa}2-2,-\frac{\hat\kappa}2,{\hat\kappa}-2,-\frac{\hat\kappa}2,{\hat\kappa}-2,-\frac{\hat\kappa}2$ at $z^+,y,z_2,x,z_1,z^-$ (in this order). 
\item Conditionally on $\hat K$, $K^{\tau_z}$ has distribution $\mu_{c_l}$ the configuration $c_l=(D_l,\tilde x=x,\tilde z_1=z_1\wedge z,\tilde y=z,\tilde z_2=z_2\wedge z',\tilde z=z')$ (the boundary arcs $(xy)$, $(yx)$ are ordered from $x$ to $y$).
\item Conditionally on $\hat K$, $K$ after $\tau_z$ has distribution $\mu_{c_r}$ in the configuration $c_r=(D_r,\tilde x=x,\tilde z_1=z_1\vee z,\tilde y=z,\tilde z_2=z_2\vee z',\tilde z=z')$.
\end{enumerate}
\end{Lem}

\begin{proof}
Given a solution $(K,\phi)$ in the configuration $c$, consider a solution $(\hat K,\phi)$ of the dual problem, as summarized in Table \ref{Tab:a}, \ref{Tab:b} depending on the position of $z$; $\hat K$ is taken independent of $K$ conditionally on the field. 
\begin{table}[htdp]
\caption{$z\in(x,z_1)$}
\begin{center}
\begin{tabular}{|c|c|c|c|c|c|c|}
\hline
$x$&$z^-$&$z$&$z^+$&$z_1$&$y$&$z_2$\\
\hline
$[\kappa]$&$\frac\kappa 2-2$&$-\frac\kappa 2$&$2$&$\frac\kappa 2-4$&$2$&$\frac\kappa 2-4$\\
\hline
$-\frac{\hat\kappa}2$&$\frac{\hat\kappa}2-2$&$[\hat\kappa]$&$-\frac{\hat\kappa}2$&$\hat\kappa-2$&$-\frac{\hat\kappa}2$&$\hat\kappa-2$\\
\hline
\end{tabular}
\end{center}
\label{Tab:a}
\end{table}%
\begin{table}[htdp]
\caption{$z\in(z_1,y)$}
\begin{center}
\begin{tabular}{|c|c|c|c|c|c|c|}
\hline
$x$&$z_1$&$z^-$&$z$&$z^+$&$y$&$z_2$\\
\hline
$[\kappa]$&$\frac\kappa 2-4$&$2$&$-\frac\kappa 2$&$\frac\kappa 2-2$&$2$&$\frac\kappa 2-4$\\
\hline
$-\frac{\hat\kappa}2$&$\hat\kappa-2$&$-\frac{\hat\kappa}2$&$[\hat\kappa]$&$\frac{\hat\kappa}2-2$&$-\frac{\hat\kappa}2$&$\hat\kappa-2$\\
\hline
\end{tabular}
\end{center}
\label{Tab:b}
\end{table}%
In the tables, entries in a row are $\rho$ parameters, except $[\kappa]$ that designates the starting point of the $\SLE_\kappa$ under consideration. Reasoning as in the case $\kappa=4$ (see Lemma \ref{loccoupl}), we see that $(K,\hat K)$ is a maximal coupling. The $\rho$ coefficients at $z^\pm$ are chosen so that $\hat K$ is the boundary of $K$ stopped when it hits $z$; this is a duality identity of the type considered in \cite{Dub_dual}, \cite{Zhan_dual}. Since $\hat K$ is determined by $K$ and independent of it conditionally on $\phi$, it is determined by $\phi$. The situation in $c_r$ is the same as in $c$, by the Markov property and the fact that $(K,\phi)$ is a solution.

The chain $K$ stays in $\overline D_l$ until it reaches $z$ at time $\tau_z$. We have to determine the distribution of $K$ up to $\tau_z$ conditionally on $\hat K$. By construction, $(K,\phi)$ and $(\hat K,\phi)$ are solutions of dual problems in $D$. 
By Lemma \ref{bcont} (as $\hat K$ gets closer to its random endpoint $z$), conditionally on $\hat K$, $\phi$ restricted to $D_l$ is a free field with $(\underline a,b)$ boundary conditions, with jumps at $x,z_1\wedge z,z_2\wedge z',z'$ (the jumps at $z^+,y$ in $D$ agglomerate in a jump at $z'$). 

The fact that $(K,\phi)$ is a solution in $D$ is a pathwise, local condition (Lemma \ref{Lem:eqpwise}). It follows that $(K^{\tau_z},\phi_{|D_l})$ is a solution in $D_l$. Uniqueness in law then determines the distribution of $K^{\tau_z}$ conditionally on $\hat K$.

\end{proof}

One can use the previous lemma to reconstruct a chordal $\SLE_\kappa$, $\kappa\geq 8$ from its dual branches as follows (see also \cite{Sheff_explor} for related considerations).

\begin{figure}[htb!]
\begin{center}
\centerline{\scalebox{.6}{\rotatebox{270}{\psfig{file=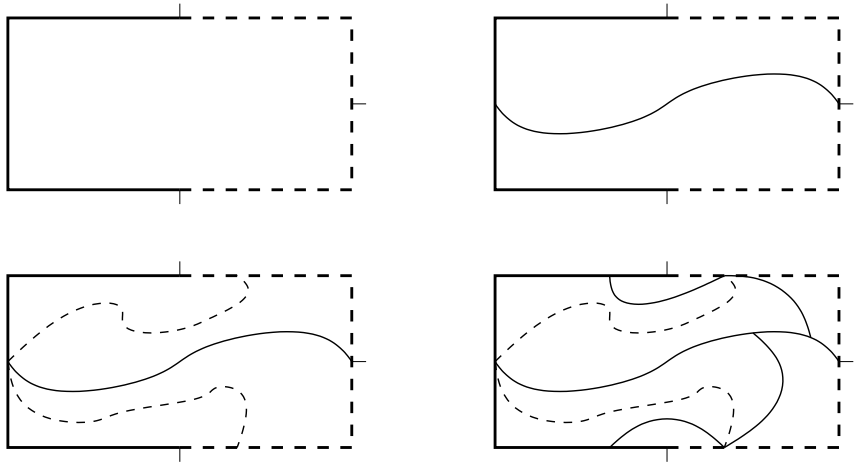}}}}
\end{center}
\caption{Iterative splitting of a domain by dual arcs ($\kappa\geq 8$)}
\label{fig:bintree}
\end{figure}

Start from a chordal $\SLE_\kappa$ in a domain $(D,x,y)$ coupled with a free field; pick a point $z$ on the boundary. This is a particular case of Lemma \ref{Lem:decomp8} with $z_1=y=z_2$. The branch $\hat K$ starting from $z$ is determined by the field; it splits $D$ into $D_l$, $D_r$. In $D_l$, the branch $\hat K_l$ starting from $z'$ is determined by the restricted field; similarly, $\hat K_r$ is the branch starting from $z'$ in $D_r$. The branch $\hat K_l$ splits $D_l$ into $D_{ll}$ and $D_{lr}$. Recursively, every subdomain is dissected in two by a branch determined by the field (see Figure \ref{fig:bintree}). All the branches are boundary arcs of the chain $K$ at some time; and the cells are visited in a prescribed order (eg at level 2, $D_{ll}$, $D_{lr}$, $D_{rl}$, $D_{rr}$). 

We can extract information on $K$ from the field in this form of the branches $\hat K$. We need to prove that all the information on $K$ can be obtained by this countable set of branches; informally, the splitting procedure yields information on $K$ at arbitrary small scales, everywhere in $D$.

\begin{Lem}\label{Lem:splitting}
Let $w_1,w_2$ be distinct interior points of $D$. In the iterative splitting of $(D,x,y)$, $w_1$ and $w_2$ are eventually in distinct cells a.s.
\end{Lem}
\begin{proof}
Without loss of generality, assume that $D$ is bounded with, say, Jordan boundary. Assume by contradiction that $w_1,w_2$ are in the same cell at any level of the splitting. Note that a.s. they are not on any branch $\hat K$ (that have zero Lebesgue measure). The splitting start from $z=z_0$ in $(xy)$ (distinct from $x,y$). The endpoint of the first branch is $z_1$; at the next level, the cell containing $w_1,w_2$ is dissected by a branch from $z_1$ to some $z_2$, and so on. The branches $\hat K_n$ from $z_n$ to $z_{n+1}$ are simple and disjoint except at their endpoints, and of ``alternating colors" (ie alternately left and right boundary arcs of $K$).

There are two possibilities: either the concatenation of the $\hat K_n$ contains infinitely many simple disjoint cycles circling around $w_1,w_2$, or eventually the successive cells containing $w_1,w_2$ have a common boundary point, and the $\hat K_n$ are arranged in a zigzag.

In the first case, the diameter of each simple cycle circling around $w_1,w_2$ is bounded away from 0. All the points on these cycles are visited by the trace in prescribed order. This contradicts the continuity of the trace of $K$ (\cite{RS01,LSW_LERW}).

In the second case, consider the harmonic measure $h_n$ of the branch $\hat K_n$ in the cell $D_n$, seen from $w_1$ or $w_2$. Then $(h_{2n})_n$ and $(h_{2n+1})_n$ are eventually increasing (and never zero), hence bounded away from 0. Since the harmonic measure of the connected set $\hat K_n$ is bounded away from 0 seen from two distinct points $w_1,w_2$, the diameter of $\hat K_n$ is also bounded away from 0; one concludes as in the first case. 

\end{proof}

One can now conclude that in a solution $(K,\phi)$, the chain $K$ is determined by the field. Indeed, the splitting of the domain is determined by the field. 
Enumerate a dense sequence of points $w_n$ in $D$; they are hit by $K$ at times $\tau_n$, which constitute a dense sequence of stopping times. Enumerate also the cells of the splitting at all levels; let $\sigma_m$ be the random time at which the trace enters the $m$-th cell, which it does at a point $x_m$ determined by the field. By the previous lemma, for $i\neq j$, the times $\tau_i,\tau_j$ are a.s. separated by a random time $\sigma_m$. Hence the family of times $\sigma_m$ is a.s. dense. The position of the continuous trace of $K$ is prescribed on an a.s. dense set of times. Thus we get pathwise uniqueness of the solution $K$.

{\bf Case $\kappa\in(4,8)$}

The argument is similar to the case $\kappa\geq 8$, however a bit more involved. Again, the $\SLE$ trace can be recovered from a tree of dual arcs which is independent conditionally on the field.

We begin with a path decomposition, analogous to Lemma \ref{Lem:decomp8}; the formal commutation identities are the same, but the geometry of paths is different. Recall in particular that for $\kappa\in (4,8)$, $\SLE$ develops cutpoints (\cite{Bef_SLE6,Dub_exc}). Hence the complement of the boundary of the $\SLE$ hull stopped at a finite time has countably many connected components (rather than just two in the $\kappa\geq 8$ case).

\begin{Lem}\label{Lem:decomp48}
In a configuration $c=(D,x,z_1,y,z_2)$, consider an $\SLE_\kappa(\underline\rho)$ chain $K$, $\kappa\in (4,8)$, $\rho=\frac\kappa 2-4,2,\frac\kappa 2-4$ at $z_1,y,z_2$ ($x,z_1,y,z_2$ in this order on the boundary), coupled with a free field $\phi$; let $\mu_c$ be the law of that $\SLE$. Let $z$ be a point on the boundary arc $(xy)$; $D_r=D\setminus K_{\tau_z}$ a random simply connected subdomain, $\tilde K=\overline{\partial D_r\cap D}$. Let $\hat K$ be a solution of the dual problem starting from $z$ (see Tables \ref{Tab:a}, \ref{Tab:b}), independent of $K$ conditionally on the field, stopped when it hits $(yx)$ at $z'$; its last visit on $(xy)$ before hitting $(yx)$ is at $z''$. Conditionally on $\hat K$, the restriction of $K$ to different connected components of $D\setminus\hat K$ are independent.
\begin{enumerate}
\item If $z\in (x,z_1)$, $\tilde K$ is the first excursion of $\hat K$ from $(xy)$ to $(yx)$. Let $D_l$ be the component of $D\setminus\hat K$ with $x$ on its boundary. Conditionally on $\hat K$, $K^{\tau_z}$ in $D_l$ has distribution $\mu_{c_l}$ where $c_l=(D_l,x,z,z,z_2\wedge z')$, stopped when it hits $z''$. After $\tau_{z''}$, the distribution of $K$ is $\mu_{c_r}$, $c_r=(D_r,z'',z_1,y,z'\vee z_2)$.
\item If $z\in (z_1,y)$, $\tilde K=\hat K$. Let $D_l=D\setminus \overline{D_r}$, $D'_l$ the connected component of $D_l$ with $x$ on its boundary. Conditionally on $\hat K$, the distribution of $K^{\tau_{z''}}$ is $\mu_{c'_l}$ in $c'_l=(D'_l,x,z_1,z'',z_2\wedge z')$; the distribution of $K$ in another connected component $D''_l$ of $D_l$ corresponding to an excursion of $\hat K$ from $y''$ to $x''$ on $(xy)$ is $\mu_{c''_l}$, $c''_l=(D''_l,x'',x'',y'',x'')$. After $\tau_z$, the distribution of $K$ is $\mu_{c_r}$, $c_r=(D_r,z,z,y,z'\vee z_2)$
\end{enumerate}
\end{Lem}
\begin{proof}
The general argument is as in Lemma \ref{Lem:decomp8}, based on the same commutation relations (Tables \ref{Tab:a}, \ref{Tab:b}), the difference being in the geometric interpretation.

In the case $z\in (x,z_1)$, consider $\hat K$ a solution of the dual problem (Table \ref{Tab:a}), independent of $K$ conditionally on the field. The path $\hat K$ hits $(xy)$ between $z$ and $z_1$ (and not in $(xz)$ or $(z_1y)$) before it first hits $(yx)$ at $z'$, where it is stopped. 
Then Lemma \ref{loccoupl} shows that $(K,\hat K)$ is a maximal coupling in the sense of \cite{Dub_dual}, for $K$ stopped at $\tau_z$ and $\hat K$ stopped at $\hat\tau_{z'}$.

For any stopping time $\hat\tau$ for $\hat K$ lesser than $\hat\tau_{z'}$, stop $K$ the first time it hits $\hat K^{\hat\tau}$ or disconnect it from $y$. Given the values of the $\rho$ parameters, it can hit $\hat K^{\hat\tau}$ only at the tip $\hat K_{\hat\tau}$; if it does not, the path $\hat K$ goes back to $(xy)$ at the point hit by $K$ at disconnection time. Hence any point of $\hat K$ on the first excursion $\tilde K$ from $(xy)$ to $(yx)$ is on $K$, while points on excursions of $\hat K$ from $(xy)$ to $(xy)$ are not. Moreover, $\hat K$ can hit $K$ only on its right boundary. Reasoning as eg in \cite{Dub_dual}, one concludes that the right boundary of $K^{\tau_z}$ is the first excursion of $\hat K$ from $(xy)$ to $(yx)$. 

Given that $(\hat K,\phi)$ is a solution, conditionally on $\hat K$, $\phi$ restricted to different connected components of $D\setminus\hat K$ is a free field with prescribed $(\underline a,b)$-type boundary conditions. Since $(K,\phi)$ is also a solution, and this is a local property (Lemma \ref{Lem:eqpwise}), this determines the distribution of $K$ in the connected components of $D\setminus\hat K$ it visits.

The case $z\in (z_1,x)$ is similar (and simpler). There (see Table \ref{Tab:b}), $\hat K$ hits $(xy)$ between $z_1$ and $z$ (and not in $(xz)$ or $(zy)$) before it first hits $(yx)$ at $z'$, where it is stopped. The chain hits $z$ a.s.; the right boundary $\tilde K$ of $K^{\tau_z}$ intersects $(xy)$ between $z_1$ and $z$. Reasoning as before, we see that $(K,\hat K)$ is a maximal coupling for $K$ stopped at $\tau_z$ and $\hat K$ stopped at $\hat\tau_{z'}$; given the choice of parameters, this implies that $\tilde K=\hat K$.
\end{proof}

As in the $\kappa\geq 8$ case, this path decomposition result can be used recursively to describe a chordal $\SLE_\kappa$ by a collection of dual paths determined by the field.
\begin{figure}[htb!]
\begin{center}
\centerline{\scalebox{.6}{\rotatebox{270}{\psfig{file=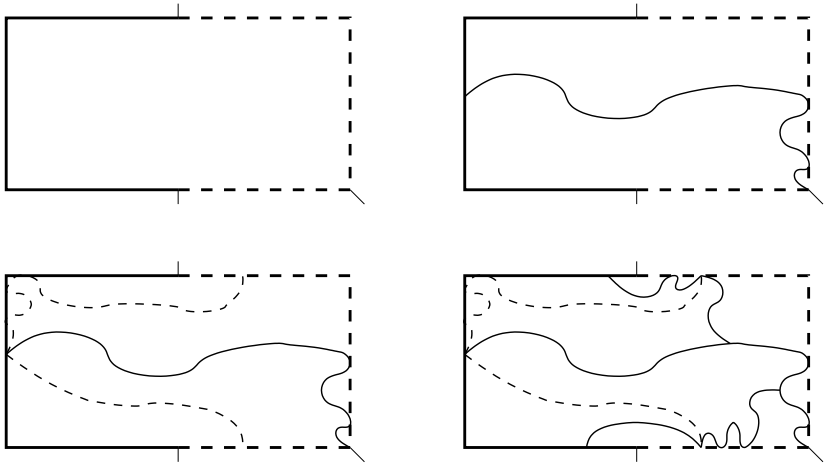}}}}
\end{center}
\caption{Iterative splitting of a domain by dual arcs ($\kappa\in (4,8)$)}
\label{fig:bintree48}
\end{figure}
We note that in this context it is rather natural to consider not simply chordal $\SLE_\kappa$, but a fuller version, such as branching $\SLE_\kappa$ (\cite{CN_full, Sheff_explor}).

Let us start with a chordal $\SLE_\kappa$ in $(D,x,y)$, $4<\kappa<8$, coupled with a free field $\phi$; pick a point $z$ on $(xy)$. This is the situation of Lemma \ref{Lem:decomp48} with $z_1=y=z_2$. The original domain $D$ is split in subdomains by $\hat K$; the trace of $K$ is contained in $\overline D_l$, $\overline D_r$. The distribution of $(K,\phi)$ in these subdomains is of the general type considered in Lemma \ref{Lem:decomp48}, and we can take $z'$ to play the r\^ole of $z$. Cells intersecting $K$ are split recursively. Note that when case 2 of Lemma \ref{Lem:decomp48} applies, one gets countably many subcells (and new marked points have to be picked in these cells). The collection of dual branches thus produced is independent of $K$ conditionally on the field.

There remains to check that points separated by the trace of $K$ are also separated by the dual branches. For this purpose, it seems practical to use a slightly different subdivision scheme. One proceeds as described above; the difference being that when case 1 in Lemma \ref{Lem:decomp48} applies, one takes $z''$ (rather than $z'$) as the new  marked point in $D_l$. Also, one does not further divide cells the interior of which is not visited by the trace of $K$ (these are determined by the tree structure). This ensures that the successive arcs $\tilde K$ are boundary arcs of the original $\SLE_\kappa$ (rather than some branching version of it).

Consider $w_1,w_2$ distinct points of $D$. If at some level $n$, there is a cell containing both $w_1,w_2$ which is not further divided, this means that $w_1,w_2$ are swallowed at the same time by $K$. Reasoning as in Lemma \ref{Lem:splitting} shows that $w_1,w_2$ cannot belong both to cells in an infinite strictly decreasing sequence, as that would violate the continuity of the trace of $K$. Therefore points $w_1,w_2$ that are swallowed at distinct times by the trace are separated by the collection of dual arcs constructed above. For $(w_n)$ a dense sequence of points in $D$, the sequence of stopping times $(\tau_{w_n})$ is a.s. dense, and one concludes as in the case $\kappa\geq 8$.

By considering successively the cases $\kappa=4$, $\kappa\in (0,4)$, $\kappa\in [8,\infty)$, $\kappa\in (4,8)$, we have established the following theorem:

\begin{Thm}\label{Thm:SDEstrong}
Let $\kappa>0$, $a=\pm\sqrt\frac{2}{\pi\kappa}$, $b=a(1-\frac\kappa 4)$. Then the stochastic equation in $(\phi,K_.)$:
$$\E(\phi_{|D_t}|{\mc F}^{\FF}_{\partial D_t})=h_{a,b}((K^t_.))\quad \forall t\geq 0$$
has a strong solution. It is pathwise unique and distributed as chordal $\SLE_\kappa$ in $(D,x,y)$.
\end{Thm}

The boundary condition for the field involves the embedding of the domain in the plane. From pathwise uniqueness, one can deduce the following covariance result:

\begin{Cor}
Let $\psi:(D,x,y)\rightarrow (\tilde D,\tilde x,\tilde y)$ be an equivalence of configuration. Let $(K,\phi)$, $(\tilde K,\tilde\phi)$ be solutions of the equation as in Theorem \ref{Thm:SDEstrong}. If the fields are coupled so that 
$$\phi=\tilde\phi\circ\psi-b\arg\psi'$$
then $\psi(K_.)=\tilde K_.$ a.s., up to time parameterization.
\end{Cor}
Note that $\psi$ is in particular a diffeomorphism, hence operates on distributions, which is what is meant by $\circ$.
\begin{proof}
This follows immediately from pathwise uniqueness and the fact that both members of the equation in Theorem \ref{Thm:SDEstrong} have the same covariance rule.
\end{proof}

The use of the chordal model is essentially conventional and for simplicity. However we observe that existence of strong solutions and pathwise uniqueness are local properties, hence they hold in more general settings.

To illustrate this, consider the following situation. Let $D$ be a planar  simply connected domain. Points $x_1,\dots,x_m$ are marked on the boundary and a point $y$ is marked in the bulk. (One could mark more points in the bulk). Consider a field with $(\underline a,b)$ boundary conditions, where $\kappa>0$, $a_1=\pm\sqrt\frac{2}{\pi\kappa}$, $b=a_1(1-\frac\kappa 4)$. Let $\partial$ be the smallest boundary arc containing all marked points except $x_1$. Then:

\begin{Cor}\label{Cor:SDEstrong}
The stochastic equation in $(\phi,K_.)$:
$$\E(\phi_{|D_t}|{\mc F}^{\FF}_{\partial D_t})=h_{\underline a,b}((K^t_.))\quad \forall t\geq 0$$
on chains stopped when they first hit $\partial\cup\{y\}$ or disconnect $y$ has a strong solution. It is pathwise unique and distributed as $\SLE_\kappa({\mc Z}_{\underline a,b})$ in $(D,x,y)$.
\end{Cor}
The partition function ${\mc Z}_{\underline a,b}$ is as in Theorem \ref{PFident}. The resulting $\SLE$ is a radial $\SLE_\kappa(\underline\rho)$ (one can also omit $y$ to get a chordal $\SLE_\kappa(\underline\rho)$).

\begin{proof}
Let $U$ be a subdomain of $D$ having on its boundary a boundary
an arc of $\partial D$ containing $x_1$; assume also that $U$ is simply connected and at positive distance of other marked points. 

Consider a solution $(K,\phi)$ of the stochastic equation with $(\underline a,b)$ boundary conditions and restrict it to $(K^{\tau_U},\phi_{|U})$ ($\tau_U$ time of first exit of $U$ by $K$). One can construct a solution $(\tilde K,\tilde\phi)$ of the chordal problem in $(D,x_1,x_2)$ by applying the density of $\frac{d\mu^{\FF}_{(a,b)}}{d\mu^{\FF}_{(\underline a,b)}}(\phi_{|U})$ (such a density exists by the Cameron-Martin formula) to the restricted pair $(K^{\tau_U},\phi_{|U})$ and extending it to a solution of the chordal problem in $D\setminus K^{\tau_U})$ (by Theorem \ref{Thm:SDEweak}). 

Considering two solutions $(K_1,\phi)$, $(K_2,\phi)$, one can similarly extend them (after $\tau_U$ for the chain, outside of $U$ for the field) to solutions $(\tilde K_1,\tilde\phi)$, $(\tilde K_2,\tilde\phi)$ of the chordal problem in $(D,x_1,x_2)$. Pathwise uniqueness holds for that problem, so that $\tilde K_1,\tilde K_2$ are a.s. equal. Since $K_i,\tilde K_i$ agree before $\tau_U$ and the measures are mutually absolutely continuous, it follows that $K_1,K_2$ agree a.s. before $\tau_U$.

By considering a countable set of such $U$'s, it follows that $K_1,K_2$ agree a.s. until they hit $\partial\cup\{y\}$ or disconnect $y$. Thus we have pathwise uniqueness. The distribution of the solution was identified in Theorem \ref{Thm:SDEweak}.
\end{proof}

\section{Some consequences}

In this section we gather some consequences of the previous constructions. We first describe a continuous version of Temperley's bijection. Some path decompositions of $\SLE$ are then listed.

\subsection{Temperley's bijection in the continuum}

In \cite{Ken_domino_conformal,Ken_domino_GFF}, Kenyon proves the convergence of the height function of a domino tiling in a simply connected Jordan domain with smooth boundary to the massless free field; the boundary condition is of $(a,b)$ type, with $b=\frac{2}\pi.\frac{\sqrt\pi}4$, which is coherent with the expression $b=\pm\frac{\kappa-4}{\sqrt{8\pi\kappa}}$ for $\kappa\in\{2,8\}$. The jump ($a=-2\pi b$) is at a marked boundary point corresponding to the root of the associated spanning tree.

In the discrete setting, there is a bijection between height functions (satisfying appropriate local conditions) and spanning trees. The goal here is to prove that the correspondence still holds in the continuum, ie the tree can be recovered from the field and vice versa.

The notion of scaling limit of a uniform spanning tree is analyzed in \cite{Sch99}. In \cite{LSW_LERW}, Lawler, Schramm and Werner prove convergence of the Peano path of a UST to $\SLE_8$ for appropriate boundary conditions. The approach here will be closer to the one in \cite{Sch99}.

Start from a free field $\phi$ in $(D,x)$ with $(a,b)$ boundary conditions, $D$ a simply connected domain (say with smooth Jordan boundary), $b=\frac{1}{2\sqrt\pi}$, $a=-2\pi b$. Pick a boundary point $z\in\partial D$ distinct from $x$. Take a path $\gamma$ coupled with the field which is an $\SLE_2(-1,-1)$ started at $z,z^-,z^+$, aiming at $x$. There is pathwise uniqueness in this situation, since this is a local property (as in Corollary \ref{Cor:SDEstrong}) and one can reason from the duality identities as in Table \ref{Tab:a}.

This splits the original domain $D$ into two subdomains $D_l,D_r$. Pick another boundary point $z_l,z_r$ in each of these domains, distinct from all marked points. Then $D_l,D_r$ can be split iteratively as in the proof of pathwise uniqueness for $\kappa\geq 8$ (see Figure \ref{fig:bintree}). The branches obtained in this fashion determine the scaling limit of the tree (either in the Peano path or collection of branches formalism).

Conversely, assume given the tree, in the form of the countable collection of branches described above (this can be deduced from the limiting objects considered in \cite{Sch99,LSW_LERW}). Let $\phi_0$ be the mean of the field with $(a,b)$ boundary conditions as above, a harmonic function in $D$. At level 1, the domain $D$ is split in $D_l,D_r$ by a curve $\gamma$. Let $\phi_1$ be the function which is harmonic in $D_l$, $D_r$ with $(\underline a,b)$ boundary conditions as in Lemma \ref{Lem:decomp8}; if $\gamma$ is coupled with a field $\phi$ as above, $\phi_1=\E(\phi|\gamma)$. The function $\phi_n$ is defined recursively by the tree at level $n$; it is harmonic on the complement of the branches; in a coupling with the a free field $\phi$, it is the expected value of the field conditionally on the tree at level $n$. For two nested cells $D_n\subset D_m$ at levels $n\geq m$ with a common boundary arc, $\phi_n-\phi_m$ converges to zero (eg in the sense of Lemma \ref{Lem:eqpwise}) at the common boundary arc. This fixes the offset of $\phi_n$ in cells with rough boundaries.

This defines from the continuous tree a sequence of a.e. harmonic functions $\phi_n$. The point is to prove that this converges to a free field. Let $\phi$ be a free field coupled with the tree as above, so that $\phi_n$ is the expected value of the field given the branches at level $\leq n$. Let ${\mc F}_n$ be the $\sigma$-algebra generated by the tree at level $n$ and $G_n$ be the Green kernel in the complement of the tree at level $n$, with Dirichlet boundary conditions. Then if $f\in C_0^\infty(D)$, 
$$\E(\langle\phi-\phi_n,f\rangle^2)=\E(\E(\langle\phi-\phi_n,f\rangle^2|{\mc F_n}))=\E(\int_{D^2} f(x)G_n(x,y)f(y)dA(x)dA(y)).$$
Notice that $G_n(x,y)$ is nonnegative and decreasing in $n$ for fixed $x,y$ (domain monotonicity). Moreover, by Lemma \ref{Lem:splitting}, $G_n(x,y)=0$ eventually for fixed $x\neq y$. It follows that $\langle\phi-\phi_n,f\rangle$ converges to 0 in $L^2$. This implies that $\langle \phi,f\rangle$ is $\mc F_\infty=\sigma({\mc F}_1,{\mc F}_2,\dots)$- measurable. As this holds for all $f\in C^\infty_0(D)$, $\phi$ itself is ${\mc F}_\infty$-measurable.

Note that there is a similar correspondence between fields and chordal $\SLE$ for any $\kappa\geq 8$. For $\kappa\in(4,8)$, the data of one chordal $\SLE$ path is not sufficient to reconstruct the field; it can be expected that reconstruction is possible from a ``fuller" version, such as branching $\SLE_\kappa$ (\cite{CN_full,Sheff_explor}), with similar arguments.

\subsection{Strong duality identities}

Duality for $\SLE$, conjectured by Duplantier, states that boundary arcs of $\SLE_\kappa$, $\kappa>4$, can be described as (versions of) $\SLE_{\hat\kappa}$, $\hat\kappa=16/\kappa$. Various such identities are established in \cite{Dub_dual,Zhan_dual}. In \cite{Dub_kapparho}, ``strong" duality identities are conjectured; these bear on the joint distribution of an $\SLE_\kappa$ and its boundary (rather than just the marginal distribution of the boundary) and are based on computations that can be understood in terms of partition function identities.

We note that such identities have been established en route to proving  pathwise uniqueness in the $\kappa\geq 8$, $\kappa\in (4,8)$. For clarity, let us make separate statements for the two cases. Many variants are possible.

\begin{Prop}\label{Prop:strongdual}
Consider $K$ a chordal $\SLE_\kappa$ from $0$ to $\infty$ in $\H$, $\tau_x$ the time at which $x>0$ is swallowed ($\kappa>4$, $\hat\kappa=16/\kappa$).
\begin{enumerate}
\item If $\kappa\geq 8$, let $\hat K=\partial K_{\tau_x}\cap\H$, $D_0$ the connected component of $\H\setminus \hat K$ that has $0$ on its boundary. Then $\hat K$ is an $\SLE_{\hat\kappa}(\underline{\hat\rho})$ from $x$ to infinity stopped upon hitting $\R^-$ at $y$, $\hat\rho=-\frac{\hat\kappa} 2, \frac{\hat\kappa} 2-2,-\frac{\hat\kappa} 2$ at $0,x^-,x^+$. Conditionally on $\hat K$, $K^{\tau_x}$ is distributed as an $\SLE_\kappa(\rho)$ from $0$ to $x$ in $D_0$, $\rho=\frac{\kappa}2-4$ at $y$.
\item If $\kappa\in (4,8)$, let $\hat K$ be the boundary arc of $K$ straddling $x$; let $d$ be its right endpoint, and $D_0$ be the connected component of $\H\setminus\hat K$ that has $0,\infty$ on its boundary. Conditionally on $d$, $\hat K$ is an $\SLE_{\hat\kappa}(\underline{\hat\rho})$ from $d$ to infinity stopped upon hitting $(0,x)$ at $g$, $\hat\rho=-\frac{\hat\kappa}2,\hat\kappa-4,\hat\kappa-2$ at $0,x,d^+$.  Conditionally on $\hat K$, $K^{\tau_x}$ is distributed as an $\SLE_\kappa(\underline\rho)$ from $0$ to $\infty$ in $D_0$, $\rho=\frac{\kappa}2-4,\frac{\kappa}2-4$ at $g,d$, stopped when it hits $d$.
\end{enumerate}
\end{Prop}
\begin{proof}
The case $\kappa\geq 8$ is part of the pathwise uniqueness proof. Let us briefly discuss the modification for the case $\kappa\in (4,8)$. Conditionally on $d$, $K$ is an $\SLE_\kappa(\underline\rho)$, $\rho=\kappa-4,-4$ at $x,d$. The relevant parameters are summarized in Table \ref{Tab:c}. 
\begin{table}[htdp]
\caption{strong duality - one sided}
\begin{center}
\begin{tabular}{|c|c|c|c|c|}
\hline
$0$&$x$&$d$&$d^+$&$\infty$\\
\hline
$[\kappa]$&$\kappa-4$&$-\frac\kappa 2$&$\frac\kappa 2-4$&$2$\\
\hline
$-\frac{\hat\kappa}2$&${\hat\kappa}-4$&$[\hat\kappa]$&$\hat\kappa-2$&$-\frac{\hat\kappa}2$\\
\hline
\end{tabular}
\end{center}
\label{Tab:c}
\end{table}%
One can couple $K$, under the conditional measure and stopped at $\tau_x$, with a free field. Taking $\hat K$ a solution of the dual problem with the same field, starting from $d$ and stopped when it hits $(0,x)$, we see reasoning as before that $\hat K$ is the boundary arc of $K$ straddling $x$ (see Theorem 1 in \cite{Dub_dual}). This gives the conditional distribution of the field in $D_0$. Considering $K^{\tau_x}$ in $D_0$ shows that it is a solution of a stochastic equation there (since this is a local condition, see Lemma \ref{Lem:eqpwise}), which determines its distribution by weak uniqueness.
\end{proof}

Similarly, one can consider two-sided situations, ie versions of $\SLE_\kappa$ conditioned on both left and right boundary arcs. In \cite{Dub_exc}, properties of $\SLE_\kappa(\underline\rho)$ in $(D,x,y)$, $\rho=\kappa-4,\kappa-4$ at $x^-,x^+$ are studied; in particular, for $\kappa\in (4,8)$, it is a chain of iid ``beads". We will now briefly discuss how to identify the distribution of these beads conditionally on their boundary.

The relevant system of commuting $\SLE$'s is summarized in Table \ref{Tab:d}. 
\begin{table}[htdp]
\caption{strong duality - two sided}
\begin{center}
\begin{tabular}{|c|c|c|c|c|c|}
\hline
$x^-$&$x$&$x^+$&$y^-$&$y$&$y^+$\\
\hline
$\kappa-4$&$[\kappa]$&$\kappa-4$&$-\frac\kappa 2$&$2$&$-\frac\kappa 2$\\
\hline
$\hat\kappa-4$&$-\frac{\hat\kappa}2$&${\hat\kappa}-4$&$[\hat\kappa]$&$-\frac{\hat\kappa}2$&$2$\\
\hline
$\hat\kappa-4$&$-\frac{\hat\kappa}2$&${\hat\kappa}-4$&$2$&$-\frac{\hat\kappa}2$&$[\hat\kappa]$\\
\hline
\end{tabular}
\end{center}
\label{Tab:d}
\end{table}%
Consider three chains $K$, $\hat K_l$, $\hat K_r$ corresponding to the lines of Table \ref{Tab:d} coupled with a common field $\phi$. Reasoning on $K,\hat K_l$ shows that $\hat K_l$ is the left boundary of $K$; symmetrically, $\hat K_r$ is its right boundary. This entails pathwise uniqueness for $\hat K_l, \hat K_r$ and the fact that they do not cross (however they intersect at the cutpoints of $K$ if $\kappa\in (4,8)$). This determines the distribution of the field right of $\hat K_l$ and also the distribution of $\hat K_r$ limited to the domain right of $\hat K_l$. Consequently, one gets the distribution of the field between $\hat K_l$ and $\hat K_r$, and finally the distribution of $K$ in the domain (or chain of domains if $\kappa\in (4,8)$) delimited by $\hat K_l$, $\hat K_r$. The conclusion is that the distribution of $K$ restricted to a bead $D$ (between consecutive cutpoints $X,Y$ of $K$) is that of an $\SLE_\kappa(\underline\rho)$, from $X$ to $Y$ in $D$, $\rho=\frac{\kappa}2-4,\frac{\kappa}2-4$ at $X^-,X^+$.

{\bf Acknowledgments.} I wish to thank Oded Schramm for very interesting conversations, and Yves Le Jan for useful references and conversations; I also thank Scott Sheffield and Roland Friedrich for discussions on an earlier version of this article.

\bibliographystyle{abbrv}
\bibliography{biblio}

\begin{thebibliography}{10}

\bibitem{Alv_boundary}
O.~Alvarez.
\newblock Theory of strings with boundaries: fluctuations, topology and quantum
  geometry.
\newblock {\em Nuclear Phys. B}, 216(1):125--184, 1983.

\bibitem{BB_review}
M.~Bauer and D.~Bernard.
\newblock 2{D} growth processes: {SLE} and {L}oewner chains.
\newblock {\em Phys. Rep.}, 432(3-4):115--221, 2006.

\bibitem{Bef_SLE6}
V.~Beffara.
\newblock Hausdorff dimensions for {$\rm SLE\sb 6$}.
\newblock {\em Ann. Probab.}, 32(3B):2606--2629, 2004.

\bibitem{CN_full}
F.~Camia and C.~M. Newman.
\newblock {The Full Scaling Limit of Two-Dimensional Critical Percolation}.
\newblock {\em preprint, arXiv:math.PR/0504036}, 2005.

\bibitem{DaPrato_LMS}
G.~Da~Prato and J.~Zabczyk.
\newblock {\em Second order partial differential equations in {H}ilbert
  spaces}, volume 293 of {\em London Mathematical Society Lecture Note Series}.
\newblock Cambridge University Press, Cambridge, 2002.

\bibitem{Dub_kapparho}
J.~Dub{\'e}dat.
\newblock {${\rm SLE}(\kappa,\rho)$} martingales and duality.
\newblock {\em Ann. Probab.}, 33(1):223--243, 2005.

\bibitem{Dub_Euler}
J.~Dub{\'e}dat.
\newblock Euler integrals for commuting {SLE}s.
\newblock {\em Journal Statist. Phys.}, 123(6):1183--1218, 2006.

\bibitem{Dub_exc}
J.~Dub{\'e}dat.
\newblock Excursion decompositions for {SLE} and {W}atts' crossing formula.
\newblock {\em Probab. Theory Related Fields}, 134(3):453--488, 2006.

\bibitem{Dub_Comm}
J.~Dub{\'e}dat.
\newblock Commutation relations for {SLE}.
\newblock {\em Comm. Pure Applied Math.}, 60(12):1792--1847, 2007.

\bibitem{Dub_dual}
J.~Dub\'edat.
\newblock {Duality of Schramm-Loewner Evolutions}.
\newblock {\em preprint, arXiv:0711.1884}, 2007.

\bibitem{Dub_Vir}
J.~Dub{\'e}dat.
\newblock S{LE} partition functions, $\zeta$-regularization and {V}irasoro
  representations.
\newblock {\em in preparation}, 2007.

\bibitem{FriKal}
R.~Friedrich and J.~Kalkkinen.
\newblock On conformal field theory and stochastic {L}oewner evolution.
\newblock {\em Nuclear Phys. B}, 687(3):279--302, 2004.

\bibitem{Gaw_CFT}
K.~Gaw{\c{e}}dzki.
\newblock Lectures on conformal field theory.
\newblock In {\em Quantum fields and strings: a course for mathematicians, Vol.
  1, 2 (Princeton, NJ, 1996/1997)}, pages 727--805. Amer. Math. Soc.,
  Providence, RI, 1999.

\bibitem{Glimm_Jaffe}
J.~Glimm and A.~Jaffe.
\newblock {\em Quantum physics}.
\newblock Springer-Verlag, New York, second edition, 1987.
\newblock A functional integral point of view.

\bibitem{Janson}
S.~Janson.
\newblock {\em Gaussian {H}ilbert spaces}, volume 129 of {\em Cambridge Tracts
  in Mathematics}.
\newblock Cambridge University Press, Cambridge, 1997.

\bibitem{Ken_domino_conformal}
R.~Kenyon.
\newblock Conformal invariance of domino tiling.
\newblock {\em Ann. Probab.}, 28(2):759--795, 2000.

\bibitem{Ken_domino_GFF}
R.~Kenyon.
\newblock Dominos and the {G}aussian free field.
\newblock {\em Ann. Probab.}, 29(3):1128--1137, 2001.

\bibitem{KenWil_grove}
R.~Kenyon and D.~Wilson.
\newblock {Boundary Partitions in Trees and Dimers}.
\newblock {\em preprint, arXiv:math.PR/0608422}, 2006.

\bibitem{KPW}
R.~W. Kenyon, J.~G. Propp, and D.~B. Wilson.
\newblock Trees and matchings.
\newblock {\em Electron. J. Combin.}, 7:Research Paper 25, 34 pp. (electronic),
  2000.

\bibitem{Kont_arbeit}
M.~Kontsevich.
\newblock {SLE, CFT, and phase boundaries}.
\newblock {\em Arbeitstagung 2003, preprint, MPI 2003 (60)}.

\bibitem{KontSuh}
M.~Kontsevich and Y.~Suhov.
\newblock On {M}alliavin measures, {SLE} and {CFT}.
\newblock {\em preprint, arXiv:math-ph/0609056}, 2006.

\bibitem{LSW_restr}
G.~Lawler, O.~Schramm, and W.~Werner.
\newblock Conformal restriction: the chordal case.
\newblock {\em J. Amer. Math. Soc.}, 16(4):917--955 (electronic), 2003.

\bibitem{Law}
G.~F. Lawler.
\newblock {\em Conformally invariant processes in the plane}, volume 114 of
  {\em Mathematical Surveys and Monographs}.
\newblock American Mathematical Society, Providence, RI, 2005.

\bibitem{LSW_LERW}
G.~F. Lawler, O.~Schramm, and W.~Werner.
\newblock Conformal invariance of planar loop-erased random walks and uniform
  spanning trees.
\newblock {\em Ann. Probab.}, 32(1B):939--995, 2004.

\bibitem{Law_Fer_RWLS}
G.~F. Lawler and J.~A. Trujillo~Ferreras.
\newblock Random walk loop soup.
\newblock {\em Trans. Amer. Math. Soc.}, 359(2):767--787 (electronic), 2007.

\bibitem{LW}
G.~F. Lawler and W.~Werner.
\newblock The {B}rownian loop soup.
\newblock {\em Probab. Theory Related Fields}, 128(4):565--588, 2004.

\bibitem{LJ_loops}
Y.~{Le Jan}.
\newblock {Markov loops, determinants and Gaussian fields}.
\newblock {\em {arxiv:math.PR/0612112}}, 2006.

\bibitem{OPS_extr}
B.~Osgood, R.~Phillips, and P.~Sarnak.
\newblock Extremals of determinants of {L}aplacians.
\newblock {\em J. Funct. Anal.}, 80(1):148--211, 1988.

\bibitem{Pol_bosonic}
A.~M. Polyakov.
\newblock Quantum geometry of bosonic strings.
\newblock {\em Phys. Lett. B}, 103(3):207--210, 1981.

\bibitem{RY}
D.~Revuz and M.~Yor.
\newblock {\em Continuous martingales and {B}rownian motion}, volume 293 of
  {\em Grundlehren der Mathematischen Wissenschaften}.
\newblock Springer-Verlag, Berlin, third edition, 1999.

\bibitem{RS01}
S.~Rohde and O.~Schramm.
\newblock Basic properties of {SLE}.
\newblock {\em Ann. of Math. (2)}, 161(2):883--924, 2005.

\bibitem{Sch99}
O.~Schramm.
\newblock Scaling limits of loop-erased random walks and uniform spanning
  trees.
\newblock {\em Israel J. Math.}, 118:221--288, 2000.

\bibitem{SS_freefield}
O.~Schramm and S.~Sheffield.
\newblock {Contour lines of the two-dimensional discrete Gaussian free field}.
\newblock {\em preprint, arXiv:math.PR/0605337}, 2006.

\bibitem{SS_FF2}
O.~Schramm and S.~Sheffield.
\newblock {In preparation}.
\newblock 2007.

\bibitem{Sheff_explor}
S.~Sheffield.
\newblock {SLE exploration trees and conformal loop ensembles}.
\newblock {\em preprint, arXiv:math.PR/0609167}, 2006.

\bibitem{Sheff_GFF}
S.~Sheffield.
\newblock Gaussian free fields for mathematicians.
\newblock {\em Probab. Theory Related Fields}, 139(3-4):521--541, 2007.

\bibitem{Simon_Pphi}
B.~Simon.
\newblock {\em The {$P(\phi )\sb{2}$} {E}uclidean (quantum) field theory}.
\newblock Princeton University Press, Princeton, N.J., 1974.
\newblock Princeton Series in Physics.

\bibitem{Simon_Trace}
B.~Simon.
\newblock {\em Trace ideals and their applications}, volume 120 of {\em
  Mathematical Surveys and Monographs}.
\newblock American Mathematical Society, Providence, RI, second edition, 2005.

\bibitem{Sonoda}
H.~Sonoda.
\newblock Functional determinants on punctured {R}iemann surfaces and their
  application to string theory.
\newblock {\em Nuclear Phys. B}, 294(1):157--192, 1987.

\bibitem{W1}
W.~Werner.
\newblock Random planar curves and {S}chramm-{L}oewner evolutions.
\newblock In {\em Lectures on probability theory and statistics}, volume 1840
  of {\em Lecture Notes in Math.}, pages 107--195. Springer, Berlin, 2004.

\bibitem{Zhan_dual}
D.~Zhan.
\newblock Duality of chordal {SLE}.
\newblock {\em preprint, arXiv:0712.0332}, 2007.

\bibitem{DZ_revers}
D.~Zhan.
\newblock Reversibility of chordal {SLE}.
\newblock {\em to appear, Ann. Probab.; arXiv:math.PR/0705.1852}, 2007.

\end{thebibliography}

-----------------------

Department of Mathematics

The University of Chicago

\end{document}